\documentclass{amsart}
\usepackage{amscd}
\usepackage{amssymb}
\usepackage{hyperref}
\usepackage{cleveref}
\usepackage{url}

\RequirePackage{tikz-cd}
\RequirePackage{amssymb}
\usetikzlibrary{calc}
\usetikzlibrary{decorations.pathmorphing}

\RequirePackage{tikz}

\newcommand{\comment}[1]{}	

\input xypic
\xyoption{all} 
\numberwithin{equation}{section}

\theoremstyle{plain}
\newtheorem{theorem}{Theorem}[subsection]
\newtheorem{corollary}[theorem]{Corollary}
\newtheorem{lemma}[theorem]{Lemma}
\newtheorem{proposition}[theorem]{Proposition}

\theoremstyle{definition}
\newtheorem{definition}[theorem]{Definition}
\newtheorem{example}[theorem]{Example}
\newtheorem{observation}[theorem]{Observation}

\theoremstyle{remark}
\newtheorem{remark}[theorem]{Remark}

\begin{document}
\bibliographystyle{plain}

\title{Order reduction of $\Lambda$-marked monomial ideals and weak resolutions}

\author{K. Machida}
\thanks{Machida was supported by the David Lachlan Hay Memorial Fund through the University of Melbourne Faculty of Science Postgraduate Writing-Up Award.  }
\address{School of Mathematics and Statistics, University of Melbourne,
VIC 3010, Australia}
\email{machida@student.unimelb.edu.au}

\date{\today}

\begin{abstract}
Borger’s theory of $\Lambda$-spaces imbues algebraic spaces, which include schemes, with an additional structure defined by an extension of the Witt vector functor. Motivated by $\mathbb{F}_1$-geometry, we prove the existence of a weak resolution of singularities in the category of $\Lambda$-schemes. Our arguments are based on standard arguments in characteristic $0$ using the order reduction of an ideal marked with $\Lambda$-equivariant data. This paper is based on work from the author's PhD thesis.
\end{abstract}

\maketitle

\section*{}

Given an algebraic variety $X$, a \textit{resolution of singularities} of $X$ is a morphism $f: X' \rightarrow X$, such that $X'$ is smooth and $f$ is proper and birational. Often $f$ is also required to induce an isomorphism $f^{-1}(X^{reg}) \cong X^{reg}$ over the locus of regular points. 
This paper is concerned with a not often studied variant of resolution of singularities, the question of resolutions of singularities in a non-standard geometry; $\mathbb{F}_1$-geometry. The origins of this theory are attributed to Jacques Tits, who in the 1950s postulated the existence of a hypothetical geometry over a base deeper than $\textrm{Spec} \,\mathbb{Z}$ that would explain analogues between projective geometries and finite sets. There have since been many attempts to define an appropriate category for $\mathbb{F}_1$-geometry and although they differ in their approach a common feature is the presence of (pointed) monoid schemes as defined by Corti\~nas, Haesemeyer, Walker and Weibel in \cite{monoid}. However, although we take inspiration from monoid schemes the category of $\mathbb{F}_1$-geometry in which we consider resolutions of singularities is Borger's theory of $\Lambda$-spaces (see \cite{borger2009lambdarings}, \cite{borger2015basicI} and \cite{borger2015basicII}). Central to the definition of Borger's $\Lambda$-spaces is the monad endofunctor $W^{\ast}$ on the category of spaces $\textbf{Sp}$, the category of sheaves of sets on the category of affine schemes with the \'etale topology. The category of $\Lambda$-spaces has objects consisting of sheafs $X$ with a $\Lambda$-structure, where this is an action of the monad $W^{\ast}$. Morphisms are given by morphims in the category of spaces that are $W^{\ast}$-equivariant. We consider the subcategory of $\Lambda$-schemes; the category of schemes with a $\Lambda$-structure and $W^{\ast}$-equivariant scheme morphisms (see \cite[Section 1]{borger2009lambdarings}). For an affine scheme $X$ flat over $\textrm{Spec} \,\mathbb{Z}$ an action of $W^{\ast}$ reduces to a family of commuting ring endomorphisms $\psi_p:X \rightarrow X$ for each prime $p$, which agree with the $p$-th power Frobenius map. Although we use $W^{\ast}$ more generally in \cref{C:2} most of the $\Lambda$-structures we consider arise from these ring endomorphisms and one can reasonably read the paper keeping only this characterisation in mind. $W^{\ast}$-equivariant morphisms are then just morphisms that commute with each $\psi_p$. The main $\Lambda$-structure we consider is the \textit{toric} $\Lambda$-structure. Let $A$ be a monoid and $\mathbb{Z}[A]$ the monoid algebra over $\mathbb{Z}$. The morphism induced by $a \mapsto a^p, \space a \in A$ gives $\mathbb{Z}[A]$ a $\Lambda$-ring structure and hence $\textrm{Spec} \, \mathbb{Z}[A]$ a $\Lambda$-scheme structure. A ring morphism $f: \mathbb{Z}[A] \rightarrow \mathbb{Z}[B]$ induced by a monoid morphism $A \rightarrow B$ commutes with $a \mapsto a^p$ and so induces a $\Lambda$-equivariant morphism $\textrm{Spec} \,\mathbb{Z}[B] \rightarrow \textrm{Spec} \,\mathbb{Z}[A]$. Combining with properties of colimits in the category of $\Lambda$-schemes (see \cite[Section 1.2]{borger2009lambdarings}) implies that $\mathbb{Z}$-realization (see \cite[Section 5]{monoid}) induces a functor
$$
(-)_{\mathbb{Z}}: \textrm{MSchemes} \rightarrow \Lambda-\textrm{schemes}.
$$
$\Lambda$-schemes given by $\mathbb{Z}$-realizations of monoid schemes are the model for most $\Lambda$-schemes in this paper. In particular, we only consider $\Lambda$-schemes flat over $\textrm{Spec} \,\mathbb{Z}$. This is in line with Borger's own thinking that $\Lambda$-schemes of finite type over $\mathbb{Z}$ should come in some way from toric $\Lambda$-structure (see \cite[Introduction]{borger2009lambdarings}). 
\\The advantage of choosing $\Lambda$-schemes for a $\mathbb{F}_1$-geometry that admits resolutions of singularities is clear. We have access to all usual scheme theory and only need to consider this additional structure. A $\mathbb{F}_1$-resolution of singularities in $\Lambda$-schemes, a \textit{$\Lambda$-resolution}, will then naturally be some variation of a usual resolution compatible with $\Lambda$-structure. More precisely, let $X$ be a $\Lambda$-scheme flat over $\textrm{Spec} \,\mathbb{Z}$. A $\Lambda$-resolution of singularities of $X$ will be a $\Lambda$-equivariant morphism $f:X' \rightarrow X$ such that $X'$, motivatived by monoid schemes, is smooth over $\mathbb{Z}$ and $f$ is proper and birational. 

Motivated by the existence of $\Lambda$-resolutions of singularities arising from the $\mathbb{Z}$-realizations of monoid schemes (see \cite[Theorem 11.1, Theorem 14.1]{monoid}), this paper proves the existence of more general $\Lambda$-resolutions and can be viewed as an exploration of what properties enable a Hironaka-style resolution argument to be adapted to the $\Lambda$-equivariant context.
\\We draw from two main sources for methods of resolutions of singularities. Desingularization of embedded toric varieties over perfect fields by Bierstone and Milman in \cite{BM1} and for general varieties in characteristic $0$ by Kollar in \cite{10.2307/j.ctt7rptq}. Choosing analogous definitions carefully and ensuring compatibility with $\Lambda$-structure form the bulk of the work of the paper. The arguments of resolutions of singularities are often intricate, and so achieving resolutions of singularities in much generality has been a difficult task. Moreover, as resolutions of singularities in the context of $\mathbb{F}_1$-geometry, in particular $\Lambda$-schemes is relatively unexplored, much of the analogous theory for varieties has had to be developed. However, because of this, the results of the paper can be considered novel.  

The standard method to produce resolutions of singularities is by successive blow-ups. In \cref{C:2} we develop a condition on a closed subscheme $Z \subset X$ of a $\mathbb{Z}$-flat $\Lambda$-scheme $X$ such that $Z,X, B_Z X, E_Z X$ are $\Lambda$-schemes, the morphisms are $\Lambda$-equivariant and so the blow-up diagram
occurs in the category of $\Lambda$-schemes. In the affine case this condition is $\mathbb{Z}$-flatness of the normal cone and the existence of generators of ideals that are "toric". The $\Lambda$-structure on $B_Z X$ can then be described easily using Frobenius lifts, which we later show equivalently come from the universal property of the blow-up. Our arguments use the concrete description of the monad structure following from commuting diagrams in \cite{borger2015basicII} that follow from local diagrams in \cite{borger2015basicI}. To extend to the general case, we assume the $\Lambda$-structure of $X$ comes from gluing affine $\Lambda$-structures. This restricts the $\Lambda$-schemes we can consider, but means that we can analyse using affine blow-up algebras and easily describe a $\Lambda$-equivariant strict transform. 

In \cref{C:3} we analyse local blow-ups for smooth schemes over $\mathbb{Z}$. We show that as for varieties in characteristic $0$ we can use local blow-up coordinates, which let us easily compare the local behaviour of blow-ups at closed points. We also introduce simple normal crossings divisors that are flat over $\mathbb{Z}$ and are compatible with $\Lambda$-structure. 

Our resolutions of singularities follow from order reduction of ideals and in \cref{C:4} we develop the necessary properties of order for this purpose. We use ideals marked by nonnegative integers as in standard resolution of singularities and describe their finite sums and relation to higher order differential operators. Using higher order differential operators and Hasse derivatives commuting with fibers we prove upper-semicontinuity of order for certain ideals. This allows us to consider only closed points and compare the support of a marked ideal after blowing up using local blow-up coordinates. Ideals given monomially by locally toric simple normal crossings divisors satisfy upper-semicontinuity of order and define $\Lambda$-marked monomial ideals in \cref{C:5}. We describe properties of their order, pullbacks to fibers and restrictions to normal crossings embeddings in preparation for proving order reduction for $\Lambda$-marked monomial ideals.

In \cref{C:5} we define a $\Lambda$-marked monomial ideal $\underline{\mathcal{I}} = (X, N, P, \mathcal{I}, E, m)$ as a tuple consisting of a marked ideal $\mathcal{I} \subset \mathcal{O}_X$ locally generated by monomials coming from a simple normal crossings divisor $E$ and an embedding of closed $\Lambda$-schemes $P \hookrightarrow N \hookrightarrow X$ smooth over $\mathbb{Z}$. $\Lambda$-marked monomial ideals bring together the objects defined in \cref{C:3} and \cref{C:4} and combines the definitions of marked ideals by Kollar (see \cite[Chapter 3]{10.2307/j.ctt7rptq}) and Bierstone and Milman (see \cite[Section 8]{BM1}). $\Lambda$-marked monomial ideals allow us to use arguments of both Bierstone and Milman and Kollar, where appropriate, and we make use of this to prove order reduction. An order reduction of $\underline{\mathcal{I}}$ is a sequence of blow-ups that makes the order of $\mathcal{I}$ in the blow-ups of $X$ less than $m$. Order reductions of marked ideals are at the heart of many proofs for resolutions of singularities and in \cref{C:5} we prove their existence for $\Lambda$-marked monomial ideals. Our proof follows the standard two-step induction (see \cite[3.70]{10.2307/j.ctt7rptq}). The first step is the most involved and uses the method of Bierstone and Milman described in \cite[Section 8]{BM1}. This method uses explicit local equations instead of derivatives of ideals, which are employed in characteristic $0$. We remark on reasons for this distinction and how monomial derivatives can be a reasonable substitute. Our arguments for this implication are quite detailed as the proof in \cite[Section 8]{BM1} can sometimes be unclear. The second implication is more straightforward and is proved similarly in both \cite{BM1} and \cite{10.2307/j.ctt7rptq}.

Finally, in \cref{C:6} using order reduction for $\Lambda$-marked monomial ideals we apply the arguments of Kollar in \cite{10.2307/j.ctt7rptq} to prove an ideal principalization theorem, which then implies the existence of a weak resolution for certain embedded $\Lambda$-schemes. This weak resolution is not implied by \cite{monoid} and so gives a new way of resolving $\Lambda$-schemes given by the $\mathbb{Z}$-realizations of monoid schemes. We make comment on the limitations of this weak resolution and some desirable improvements, motivating upcoming work.


\section{Blow-ups}
\label{C:2} 
In this section we describe conditions that guarantee the existence of a compatible $\Lambda$-structure of a blow-up of a $\mathbb{Z}$-flat $\Lambda$-scheme. The first section describes affine blow-ups and properties of ideals that induce a $\Lambda$-structure on the Proj of the Rees algebra. We construct this $\Lambda$-structure in two ways. First, using explicit Frobenius lifts; and second, showing that the blow-up is an algebra for $W^{\ast}$ using universal properties. The first perspective is algebraically easier to work with and agreement with the second perspective shows that the $\Lambda$-structure is independent of certain choices. In the second section, we describe a class of $\mathbb{Z}$-flat $\Lambda$-schemes on which our results on affine blow-ups generalise. Finally, we briefly discuss strict transforms of certain closed $\Lambda$-subschemes that are used in \cref{C:5}.


\subsection{Affine case}
In order to understand blow-ups of affine $\Lambda$-schemes, we need to consider $\mathbb{Z}$-flat $\Lambda$-rings and their Rees algebras. Let $A$ be a ring with ideal $I$ and recall the normal cone for the closed immersion defined by $I$ is given by $\textrm{Spec} (\frac{A}{I} \oplus \frac{I}{I^2} \oplus \frac{I^2}{I^3} \oplus \dots)$. We use a characterisation of the normal cone being $\mathbb{Z}$-flat to induce a $\Lambda$-structure on the Rees algebra of a $\mathbb{Z}$-flat ring.

\begin{lemma}\label[lemma]{lemma2101} 
    Let $A$ be a ring and $I$ an ideal. $\frac{A}{I^k}$ is $\mathbb{Z}$-flat for all $k$ if and only if the normal cone is $\mathbb{Z}$-flat.
\end{lemma}
\textit{Proof.} Assume $\frac{A}{I^k}$ is $\mathbb{Z}$-flat for all $k$, then $\frac{A}{I}$ is $\mathbb{Z}$-flat and for all $k$, $\frac{I^k}{I^{k+1}} \hookrightarrow \frac{A}{I^{k+1}}$ is a subgroup of a torsionfree group and so is torsionfree itself. The direct sum of $\mathbb{Z}$-torsion free modules is $\mathbb{Z}$-torsion free and so the normal cone is $\mathbb{Z}$-flat, since $\mathbb{Z}$-flatness is equivalent to being torsionfree. Conversely, suppose that the normal cone is $\mathbb{Z}$-flat and $\frac{I^k}{I^{k+1}}$ is $\mathbb{Z}$-flat for all $k$. We argue by induction, noting that $\frac{A}{I}$ is $\mathbb{Z}$-flat by assumption. Suppose $\frac{A}{I^k}$ is $\mathbb{Z}$-flat. Then
$
0 \rightarrow \frac{I^k}{I^{k+1}} \rightarrow \frac{A}{I^{k+1}} \rightarrow \frac{A}{I^{k}} \rightarrow 0
$
is exact and $\frac{I^k}{I^{k+1}}$ and $ \frac{A}{I^{k}}$ are $\mathbb{Z}$-flat, and hence $\frac{A}{I^{k+1}}$ is also. \qedsymbol

\begin{lemma}
    \label[lemma]{lemma2102} 
    Let $A$ be a $\mathbb{Z}$-torsion free $\Lambda$-ring with Frobenius lifts $\psi_p$ and $\Lambda$-ideal $I$ such that $\psi_p(I) \subseteq I^p$ for all primes $p$. Assume that the normal cone of $I$ is $\mathbb{Z}$ flat, which by \cref{lemma2101} implies that for all primes $p$, $pf \in I^k$ implies $f \in I^k$. Denote by $A[It]:=\bigoplus_{i \geq 0} I^i t^i$ the Rees algebra of $A$ with respect to $I$. Then $A[It]$ is a $\mathbb{Z}$-flat $\Lambda$-ring.
\end{lemma}

\textit{Proof.} Consider the family of maps indexed by primes $p$
$$
\psi_{p, A[It]}: A[It] \rightarrow A[It], \quad \sum_{i=0}^n a_i t^i \mapsto \sum_{i=0}^n \psi_p(a_i) t^{pi},
$$
where $a_i \in I^i$ and $t^i$ indicates that we are considering $a_i \in I^i$. Then $\psi_{p, A[It]}$ gives a family of commuting ring endomorphisms as $\psi_p$ are commuting endomorphisms and for all $p$
\begin{align*}
    \psi_{p, A[It]}\left(\sum_{i=0}^n a_i t^i \right) - \left( \sum_{i=0}^n a_i t^i \right)^p &= \sum_{i=0}^n \psi_p(a_i) t^{pi} - \sum_{i=0}^n (a_i)^p t^{pi}-p(\dots)
    \\ &=\sum_{i=0}^n (\psi_p(a_i)-(a_i)^p) t^{pi}-p(\dots)
    \\ & = \sum_{i=0}^n pb_i t^{pi}-p(\dots)
    \\ & = p(\sum_{i=0}^n b_i t^{pi}-\dots),
\end{align*}
where the $pb_i$ come from the Frobenius lift property of $\psi_p(a_i)$ and $b_i \in I^{pi}$ by the $\mathbb{Z}$-flatness of $\frac{A}{I^k}$. Thus the $\psi_{p, A[It]}$ form a commuting family of Frobenius lifts. $A$ $\mathbb{Z}$-flat implies $I^k \subset A$ is $\mathbb{Z}$-flat for all $k$ as sub-modules of flat modules over principal ideal domains are flat. Therefore $A[It]$, being the direct sum of $\mathbb{Z}$-flat modules, is itself $\mathbb{Z}$-flat and so a $\Lambda$-ring. \qedsymbol

\begin{example}[Monoid algebras]
    \label[example]{example2103}
     Let $A:=\mathbb{Z}[M]$ be a pointed monoid algebra with toric $\Lambda$-structure and $I:=\mathbb{Z}[J]$ be the $\Lambda$-ideal generated by a proper monoid ideal $J \subset M$. For any $m \in J$, $\psi_p(m)=m^p$ so $\psi_p(I) \subseteq I^p$ for all $p$ and for any $k$, $\frac{A}{I^k} \cong \mathbb{Z}\left [ \frac{M}{J^k} \right]$ is $\mathbb{Z}$-torsion free so $A[It]$ has an induced $\Lambda$-structure. Moreover, $A[It] = \mathbb{Z}[M[Jt]]$, where $M[Jt]$ is the graded monoid as in \cite[Section 7]{monoid} and for any $m t^k \in J^k t^k \subset I^k t^k$, $\psi_{p, A[It]}(m t^k)=m^p t^{pk} = (m t^k)^p \in J^{pk} t^{pk} \in M[Jt]$ so the induced $\Lambda$-structure on $A[It]$ is precisely the toric $\Lambda$-structure.
\end{example}

\begin{example}[Completions of finitely generated monoid algebras]
    \label[example]{example2104}
    Let notation be as in the previous example, then the completion $\hat{A}_I$ is a $\mathbb{Z}$-flat $\Lambda$-ring with ideal $\hat{I}$ such that $(\hat{I})^k \cong \hat{(I^k)}$ and $\psi_p(I) \subseteq I^p$ implies $\hat{\psi_p}(\hat{I}) \subseteq (\hat{I})^p$. Assume $M$ is a finitely generated monoid, which implies $J$ and hence $I$ are finitely generated. Then $A$ is Noetherian and hence $\frac{A}{(I)^k}$ are Noetherian and for all $k$,
    $
    \frac{\hat{A}_I}{(\hat{I})^k} \cong \frac{\hat{A}_I}{\widehat{(I^k)}} \cong \widehat{\left (\frac{A}{I^k} \right)},
    $
    which is flat over $\frac{A}{(I)^k} \cong \mathbb{Z}\left [ \frac{M}{J^i} \right]$. As $\frac{\hat{A}}{(\hat{I})^k}$ is flat over $\mathbb{Z}$ for all $k$, $\hat{A}_I[\hat{I}t]$ is a $\Lambda$-ring.
\end{example}

\begin{remark}
    \label[remark]{remark2105} 
    If we further assume $M$ is cancellative and torsionfree then $\hat{A}_I$ is not a free abelian group (it contains the Baer-Specker group as a subgroup) and thus $\hat{A}_I[\hat{I}t]$, which contains $\hat{A}_I$ as a subgroup, is not free abelian and so not a monoid algebra over $\mathbb{Z}$.
\end{remark}

The Frobenius lifts on these Rees algebras are compatible with the grading and so we make the following definition.

\begin{definition}[Graded $\Lambda$-rings]
    Let $S = \bigoplus_{d \geq 0} S_d$ be a graded ring finitely generated in degree $1$ by $a_1, \dots, a_n \in S_1$. $S$ will be a \textit{graded $\Lambda$-ring} if $S$ is $\mathbb{Z}$-flat and has a $\Lambda$-ring structure given by $\psi_p$ such that $\psi_p(a_i)=(a_i)^p$ and each $S_m$ has the property that $pf \in S_m$ implies $f \in S_m$ for all $i, p, m$. Note that this implies $\psi_p(S_m) \subseteq S_{pm}$ and hence $\psi_n(S_m) \subseteq S_{nm}$, where $\psi_n$ is the endomorphism given by the composition of $\psi_p$'s corresponding to the prime decomposition of $n$. 
\end{definition}

\begin{example}[Rees algebras of $\Lambda$-rings]
    \label[example]{example2107} 
    Let $A$ be a $\mathbb{Z}$-flat $\Lambda$-ring and $I$ an ideal finitely generated by $a_1, \dots, a_n$ such that $\psi_p(a_i)=(a_i)^p$. We will say that the $a_i$ are \textit{toric} generators for $I$. Assume the normal cone of $I$ is $\mathbb{Z}$-flat so by \cref{lemma2102} $A[It]$ will be a graded $\Lambda$-ring. In particular, Rees algebras of finitely generated monoid algebras and their completions (\cref{example2103}, \cref{example2104}) are graded $\Lambda$-rings. We will show that the Proj of a graded $\Lambda$-ring is a $\Lambda$-scheme and so obtain $\Lambda$-structures on blow-ups using this example.
\end{example}

\textbf{Characterising toric generators.} Let $A$ be a $\mathbb{Z}$-flat $\Lambda$-ring and $I$ an ideal finitely generated by $a_1, \dots, a_n$. Recall that this is equivalent to a ring morphism $f: \mathbb{Z}[x_1, \dots, x_n] \rightarrow A$ such that $I$ is the ideal theoretic image of $(x_1, \dots, x_n).$ If we give $\mathbb{Z}[x_1, \dots, x_n]$ the toric $\Lambda$-structure then $f$ is $\Lambda$-equivariant if and only if $a_1, \dots, a_n$ are toric generators of $I$.


\subsubsection{$\Lambda$-structure using Frobenius lifts} We construct compatible Frobenius lifts on affine covers of $\textrm{Proj} (S)$. This will give a $\Lambda$-structure by the following observation on flatness for graded rings.

\begin{observation}
    Let $S$ be a graded ring flat over $\mathbb{Z}$. Then as flatness is preserved under localization $S_{(f)} \subset S_{f}$ are both $\mathbb{Z}$-flat rings, where $S_{(f)}$ is the degree $0$ part of the localization $S_{f}$ at homogeneous $f \in S_+$ ($S_+$ is the irrelevant ideal). It follows that $\textrm{Proj} (S) = \bigcup_{f} \textrm{Spec} (S_{(f)})$ will be a scheme flat over $\mathbb{Z}$.
\end{observation}

\begin{lemma}
    \label[lemma]{lemma2112} 
    Let $X$ be a scheme flat over $\textrm{Spec} \, \mathbb{Z}$ and $Z \subset X$ a closed subscheme given by an ideal sheaf $\mathcal{I}$. Then the blow-up $B_Z X$ is flat over $\textrm{Spec} \, \mathbb{Z}$. Further, if the normal cone of $\mathcal{I}$ is $\mathbb{Z}$-flat then the exceptional divisor $E_Z X$ is flat also.
\end{lemma}

\textit{Proof.} Let $\textrm{Spec} \, A$ be an affine open of $X$. Then $\mathcal{I}(\textrm{Spec} \, A)=I \subset A$, where $A$ is $\mathbb{Z}$-flat and the normal cone is given by $\textrm{Spec} \, S_{A, I}$ for
$
S_{A, I} = \frac{A}{I} \oplus \frac{I}{I^2} \oplus \frac{I^2}{I^3} \oplus \dots.
$
As in the proof of \cref{lemma2102} $A[It]$ is flat over $\mathbb{Z}$ and so by the observation, $\textrm{Proj} (A[It])$ is $\mathbb{Z}$-flat. The $\textrm{Proj} (A[It])$ cover $B_Z X$ as $\textrm{Spec} \, A$ varies so $B_Z X$ is $\mathbb{Z}$-flat. If the normal cone is $\mathbb{Z}$-flat then $S_{A,I}$ is $\mathbb{Z}$-flat and so $\textrm{Proj} (S_{A,I})$ is also. As the $\textrm{Proj} (S_{A,I})$ cover
$
E_Z X = \mathcal{P}roj \left(\bigoplus_{k\geq 0} \frac{\mathcal{I}^k}{\mathcal{I}^{k+1}} \right),
$
$E_Z X$ is $\mathbb{Z}$-flat. \qedsymbol

\begin{proposition}
    \label[proposition]{proposition2113} 
    Let $S$ be a graded $\mathbb{Z}$-flat $\Lambda$-ring with generators $a_1, \dots, a_n$ and Frobenius lifts $\psi_p$. Then $\textrm{Proj} (S)$ has a $\Lambda$-scheme structure induced from $S$.
\end{proposition}
\textit{Proof.} Fix $i$ and consider the $\mathbb{Z}$-flat affine scheme $\textrm{Spec} (S_{(a_i)})$, where $S_{(a_i)}$ is the degree 0 part of the localization $S_{a_i}$ and $\textrm{Proj} (S)=\bigcup_i \textrm{Spec} (S_{(a_i)})$. Recall elements of $S_{(a_i)}$ are represented by expressions of the form $\frac{x}{(a_i)^n}$, where $x \in S_n$. Consider
$$
\psi_{p,i}: S_{(a_i)} \rightarrow S_{(a_i)}, \quad \frac{x}{(a_i)^n} \mapsto \frac{\psi_p(x)}{(a_i)^{pn}}.
$$
Note that $\psi_p(x) \in S_{pn}$ as $x \in S_n$ and if $\frac{x}{(a_i)^n} \sim \frac{y}{(a_i)^m}$ then $\frac{\psi_p(x)}{(a_i)^{pn}} \sim \frac{\psi_p(y)}{(a_i)^{pm}}$ so $\psi_{p,i}$ is well-defined. For each prime $p$, $\psi_{p,i}$ defines a ring endomorphism as $\psi_p$ is an endomorphism and the family of endomorphisms will commute as the $\psi_p$'s commute. Finally 
$
\psi_{p,i}\left( \frac{x}{(a_i)^n} \right)-\left(\frac{x}{(a_i)^n} \right)^p = \frac{\psi_p(x)-x^p}{(a_i)^{pn}}=\frac{ph_x}{(a_i)^{pn}} = p \left( \frac{h_x}{(a_i)^{pn}}\right),
$
where $h_x$ comes from the Frobenius lift $\psi_p$ on $S$ and $h_x \in S_{pn}$ as $ph_x \in S_{pn}$. Then for all $i$ the $\psi_{p,i}$'s form a commuting family of Frobenius lifts on $S_{(a_i)}$ and so induce a $\Lambda$-scheme structure on $\textrm{Spec} (S_{(a_i)})$. Now for $i,j$ recall that $\textrm{Spec} (S_{(a_i)}) \cap \textrm{Spec} (S_{(a_j)}) \cong \textrm{Spec} (S_{(a_ia_j)})$, which will again be $\mathbb{Z}$-flat. Define
$$
\psi_{p,ij}: S_{(a_ia_j)} \rightarrow S_{(a_ia_j)}, \quad \frac{x}{(a_ia_j)^n} \mapsto \frac{\psi_p(x)}{(a_ia_j)^{pn}}.
$$
The $\psi_{p,ij}$ will form a commuting family of Frobenius lifts by the same arguments above and thus $\textrm{Spec} (S_{(a_ia_j)})$ is also a $\Lambda$-scheme. Recall that $\textrm{Spec} (S_{(a_ia_j)}) \rightarrow \textrm{Spec} (S_{(a_i)})$ is induced by
$$
S_{(a_i)} \rightarrow (S_{(a_i)})_{\frac{a_j}{a_i}} \cong S_{(a_ia_j)}, \quad \frac{x}{(a_i)^n} \mapsto \frac{x(a_j)^n}{(a_ia_j)^n}
$$
so that
$$
\psi_{p,i}\left( \frac{x}{(a_i)^n} \right) \mapsto \left( \frac{\psi_p(x)(a_j)^{pn}}{(a_ia_j)^{pn}} \right) = \left( \frac{\psi_p(x(a_j)^n)}{(a_ia_j)^{pn}} \right) = \psi_{p.ij} \left( \frac{x(a_j)^n)}{(a_ia_j)^{n}} \right),
$$
i.e. the ring morphism is $\Lambda$-equivariant. Thus we may glue together the affine $\Lambda$-schemes $\textrm{Spec} (S_{(a_i)})$ along their intersections and induce a $\Lambda$-structure on $\textrm{Proj} (S)$. \qedsymbol

\begin{remark}
    Although we assumed $S$ was $\mathbb{Z}$-flat we only needed $S_{(a_i)}$ to be $\mathbb{Z}$-flat for all $i$ in order for Frobenius lifts to be equivalent to a $\Lambda$-structure. In particular, this is equivalent to $\textrm{Proj} (S)$ being $\mathbb{Z}$-flat but need not imply $S$ is $\mathbb{Z}$-flat.
\end{remark} 

\begin{example}[$\Lambda$ blow-up squares]
    \label[example]{example2115} 
    Let $A$ be a $\mathbb{Z}$-flat $\Lambda$-ring with $I=(a_1, \dots, a_n)$ an ideal with toric generators as in \cref{example2107}. Write $X=\textrm{Spec} (A)$ and $Z=\textrm{Spec} (A/I)$. By \cref{proposition2113}, $\textrm{Proj} (A[It]) = B_Z X$ is a $\mathbb{Z}$-flat $\Lambda$-scheme. Further, the blow-up morphism $B_Z X \rightarrow X$ is induced by ring morphisms
    $$
    A \rightarrow S_{a_i}, \quad a \mapsto \frac{a}{1},
    $$
    which are $\Lambda$-equivariant as
    $
    \frac{\psi_p(a)}{1} = \psi_{p,i} \left(  \frac{a}{1}  \right).
    $
    Thus the blow-up morphism is a $\Lambda$-equivariant scheme morphism. The exceptional divisor $E_Z X$, being the pullback of two $\Lambda$-equivariant morphisms, is then a $\mathbb{Z}$-flat (by \cref{lemma2112}) $\Lambda$-scheme. The blow-up square
    $$
    \begin{tikzcd}
    E_{Z} X \arrow[r] \arrow[d] 
    & B_{Z} X \arrow[d] \\
    Z  \arrow[r] & X 
    \end{tikzcd}
    $$
    then occurs in the category of $\Lambda$-schemes, by which we mean the induced morphisms in $\textbf{Sp}$ are equivariant with respect to $W^{\ast}$. 
\end{example}


\begin{example}[Blow-ups of monoid algebras]
    \label[example]{example2116} 
    Let $A=\mathbb{Z}[M]$ be the $\mathbb{Z}$-realization of a monoid $M$ and $I=\mathbb{Z}[J]$ with $J=(a_1, \dots, a_n) \subset M$. For $S=A[It]=\mathbb{Z}[M[Jt]]$, $S_{(a_i)} = \mathbb{Z}[M[Jt]_{(a_i)}]$, where $M[Jt]_{(a_i)}$ is the degree $0$ part of the localization of the graded monoid $M[Jt]$. Let $\frac{x}{(a_i)^n} \in M[Jt]_{(a_i)} \subset S_{(a_i)}$. Now $x  \in M$ and $\psi_p$ on $S$ being induced from the toric $\Lambda$-structure of $A$ implies 
    $
    \psi_{p,i}\left( \frac{x}{(a_i)^n} \right) = \frac{x^p}{(a_i)^{pn}}  = \left( \frac{x}{(a_i)^{n}} \right)^p
    $
    so the $\Lambda$-structure on $S_{(a_i)}$ is the toric $\Lambda$-structure and similarly for $S_{(a_ia_j)}$. Thus, the induced $\Lambda$-scheme structure on $\textrm{Proj}(S)=\textrm{B}_{\textrm{Spec} (A/I)} \textrm{Spec} (A)$ is precisely the $\Lambda$-structure coming from the $\mathbb{Z}$-realization of the monoid scheme
    $$
    \textrm{MProj}(M[Jt]) = \textrm{B}_{\textrm{MSpec}(M/J)}\textrm{MSpec}(M)
    $$
    (see \cite[Section 7]{monoid}).
\end{example}

\begin{example}[Blow-ups of completions of monoid algebras]
    Let $A, I, M$ be as in the previous example such that $M$ is a finitely generated monoid. Then $B_{\textrm{Spec} (\hat{A}_I/\hat{I})} \textrm{Spec} (\hat{A}_I)$ is a $\mathbb{Z}$-flat $\Lambda$-scheme. Further, if $M$ is also cancellative and torsionfree then as in \cref{remark2105} it need not come from the $\mathbb{Z}$-realization of a monoid scheme. 
\end{example}



\subsubsection{$\Lambda$-structure using $W^{\ast}$}

Let $A$ be a $\mathbb{Z}$-torsion free $\Lambda$-ring with Frobenius lifts $\psi_p$ and $I=(a_1, \dots, a_n)$ an ideal meeting the conditions of \cref{example2107}, i.e. $\psi_p(a_i)=(a_i)^p$ and $\frac{A}{I^k}$ is $\mathbb{Z}$-flat for all $k$. Write $Z=\textrm{Spec} (A/I), X=\textrm{Spec} (A)$ and consider the blow-up diagram
$$
\begin{tikzcd}
E_{Z} X \arrow[r] \arrow[d] 
& B_{Z} X \arrow[d] \\
 Z  \arrow[r] & X 
\end{tikzcd}.
$$
By \cref{example2115} we know $B_Z X$ is a $\mathbb{Z}$-flat $\Lambda$-scheme and the diagram is $\Lambda$-equivariant. However, the proof of \cref{proposition2113} required constructing Frobenius lifts on a cover of $B_Z X$ and one might wonder whether other Frobenius lifts are possible and to what extent this $\Lambda$-structure on $B_Z X$ is the right one to choose. We will show that using the universal property of the blow-up, $B_Z X$ is naturally an algebra for the monad $W^{\ast}$ without having to write Frobenius lifts and that the blow-up square is $\Lambda$-equivariant. Moreover, the algebra structure will induce the same Frobenius lifts as in \cref{proposition2113} and show that the $\Lambda$-structure is independent of the generators of $I$. Recall the notation for Witt vector functors on rings from \cite[Section 1]{borger2015basicI}. We first make the following observation for convenience. 

\begin{remark}
    \label[remark]{remark1241} 
    We have a monoid isomorphism $\mathbb{N}^E \rightarrow \mathbb{N}_{>0}, \quad n=(a_1, a_2, \dots) \mapsto 2^{a_1} \cdot 3^{a_2} \cdots$, where the right hand side is $\mathbb{N}_{>0}$ under multiplication. Therefore $\Psi_{\mathbb{Z},E}$-structure in \cite[Section 1]{borger2015basicI} is equivalent to an action of $\mathbb{N}_{>0}$ where $n \in \mathbb{N}_{>0}$ will induce the endomorphism consisting of powers of $\psi_p$'s for $p$ in the prime decomposition of $n$.
\end{remark}

\begin{observation}
    For $n \in \mathbb{N}^{(E)}$ write $n' \in \mathbb{N}_{> 0}$ for the image of $n$ under the monoid isomorphism described in \cref{remark1241}. Let $\psi_n$ be the induced endomorphism then 
    \begin{enumerate}
        \item $\psi_n(a_i)=(a_i)^{n'}$.
        \item $\psi_p(\psi_n(a_i))=(a_i)^{pn'}$ for all primes $p$.
        \item $\psi_n(I^k) \subset I^{kn'}$ for all $k \geq 1$.
    \end{enumerate}
\end{observation}

Recall that as $A$ is a $\mathbb{Z}$-torsion free $\Lambda$-ring we obtain $h:A \rightarrow W(A)=W^{fl}(A)$. Consider the blow-up square now in the category $\textbf{Sp}$ and apply the functor $W^{\ast}$ to give the diagram
$$
\begin{tikzcd}
	& {W^{\ast}(E_ZX)} && {W^{\ast}(B_Z X)} \\
	{E_Z X} && {B_Z X} \\
	& {W^{\ast}(Z)} && {W^{\ast}(X)} \\
	Z && X
	\arrow[from=2-1, to=2-3]
	\arrow[from=4-1, to=4-3]
	\arrow[from=2-3, to=4-3]
	\arrow[from=2-1, to=4-1]
	\arrow[from=1-2, to=1-4]
	\arrow[from=1-4, to=3-4]
	\arrow[from=1-2, to=3-2]
	\arrow[from=3-2, to=3-4]
	\arrow[from=3-4, to=4-3]
	\arrow[from=3-2, to=4-1]
	\arrow[dashed, from=1-2, to=2-1]
	\arrow[dashed, from=1-4, to=2-3]
\end{tikzcd},
$$
where the bottom commuting square exists by $Z \rightarrow X$ being $\Lambda$-equivariant. We will complete the above diagram by producing maps 
$$
W^{\ast}(B_Z X) \rightarrow B_Z X, \quad W^{\ast}(E_Z X) \rightarrow E_Z X
$$
in $\textbf{Sp}$ making $B_Z X$ and $E_Z X$ $\Lambda$-schemes and the original blow-up square $\Lambda$-equivariant. 

\begin{proposition}
    \label[proposition]{proposition2122} 
    There exist morphisms
    $$
    W^{\ast}(B_Z X) \rightarrow B_Z X, \quad W^{\ast}(E_Z X) \rightarrow E_Z X
    $$
    making the diagram commute.
\end{proposition}

\textit{Proof.} We show that the family of schemes $W^{\ast}_n(B_Z X) \times_X Z$ are effective Cartier divisors on $W^{\ast}_{n}(B_Z X)$, where $W^{\ast}_n$ are the finite length Witt vector functors. The universal property of the blow-up then induces $W^{\ast}_n(B_Z X) \rightarrow B_Z X$ and this collection of morphisms is compatible with taking colimits and so induces $W^{\ast}(B_Z X) \rightarrow B_Z X$. Note that $B_Z X$ is flat over $\textrm{Spec} \, \mathbb{Z}$ so locally we are in the $W^{fl}=W$ case (\cite[1.12]{borger2015basicI}). For consistency we will use the notation of \cref{proposition2113}. Recall $B_Z X = \textrm{Proj} (S: = A[It])$ is covered by affine opens of the form $\textrm{Spec} (S_{(a_i)})$ (where $a_i \in I$ is considered in degree 1) such that $B_Z X \rightarrow X$ restricted to $\textrm{Spec} (S_{(a_i)})$ is induced by 
$$
\phi_i: A \rightarrow S_{(a_i)}, \quad a \mapsto \frac{a}{1}.
$$
Now $W^{\ast}_n(B_Z X)$ is a scheme (see \cite[Introduction]{borger2015basicII}) and is covered by affine open subschemes $\textrm{Spec} (W_n(S_{(a_i)})$ with intersections $\textrm{Spec} (W_n(S_{(a_{ij})})$ and thus $W^{\ast}(B_Z X) \rightarrow W^{\ast}(X)$ will be induced by morphisms $W^{\ast}_n(B_Z X) \rightarrow W^{\ast}_n(X)$, which  restricted to $\textrm{Spec} (W_n(S_{(a_i)})$ are given by $W_n(\phi_i): W_n(A) \rightarrow W_n(S_{(a_i)})$. These morphisms commute with composition with any similar system of maps that make up the colimit/limit, i.e. for $W^{\ast}_{n+j}$ and hence in the affine case for $W_{n+j}$. $W^{\ast}(X) \rightarrow X$, being the $\Lambda$-structure map, is induced by morphisms $A \rightarrow W_n(A)$, which factor through $h: A \rightarrow W(A)$ and the diagram involving projections of the ghost spaces, i.e.
$$
\begin{tikzcd}
	{W_{n}(A)} && {W_{n+j}(A)} && {W(A)} && A \\
	\\
	{A^{[0, n]}} && {A^{[0, n+j]}} && {A^{\mathbb{N}^E}}
	\arrow[hook, from=1-5, to=3-5]
	\arrow[two heads, from=3-5, to=3-3]
	\arrow[two heads, from=3-3, to=3-1]
	\arrow[from=1-5, to=1-3]
	\arrow[from=1-3, to=1-1]
	\arrow[hook, from=1-1, to=3-1]
	\arrow[hook, from=1-3, to=3-3]
	\arrow["h"', from=1-7, to=1-5]
\end{tikzcd},
$$
where the vertical maps are injective as we are in the $\mathbb{Z}$-flat case. Write $h_n$ for the composition $A \xrightarrow{h} W(A) \rightarrow W_n(A)$. Then the composition $W^{\ast}(B_Z X) \rightarrow W^{\ast}(X) \rightarrow X$ will be induced by morphisms $W^{\ast}_n(B_Z X) \rightarrow W^{\ast}_n(X) \rightarrow X$ compatible with the filtration, which are in turn induced by 
$$
W_n(\phi_i) \circ h_n:A \rightarrow W_n(A) \rightarrow W_n(S_{(a_i)})
$$
that are compatible with filtrations and intersections, i.e. extending $i,j$ to $W_n(S_{a_ia_j})$. In particular, we can consider the map $W(\phi_i): W(A) \rightarrow W(S_{a_i})$, which will fit together with the maps above in the following commuting diagram
$$
\begin{tikzcd}[cramped,sep=scriptsize]
	& {W_{n}(A)} && {W_{n+j}(A)} && {W(A)} & A \\
	{W_{n}(S_{a_i})} && {W_{n+j}(S_{a_i})} && {W(S_{a_i})} \\
	& {A^{[0, n]}} && {A^{[0, n+j]}} && {A^{\mathbb{N}^E}} \\
	{S_{a_i}^{[0, n]}} && {S_{a_i}^{[0, n+j]}} && {S_{a_i}^{\mathbb{N}^E}}
	\arrow[from=1-2, to=2-1]
	\arrow[hook, from=1-2, to=3-2]
	\arrow[from=1-4, to=1-2]
	\arrow[from=1-4, to=2-3]
	\arrow[hook, from=1-4, to=3-4]
	\arrow[from=1-6, to=1-4]
	\arrow[from=1-6, to=2-5]
	\arrow[hook, from=1-6, to=3-6]
	\arrow["h"'{pos=0.4}, from=1-7, to=1-6]
	\arrow[hook, from=2-1, to=4-1]
	\arrow[from=2-3, to=2-1]
	\arrow[hook, from=2-3, to=4-3]
	\arrow[from=2-5, to=2-3]
	\arrow[hook, from=2-5, to=4-5]
	\arrow[from=3-2, to=4-1]
	\arrow[two heads, from=3-4, to=3-2]
	\arrow[from=3-4, to=4-3]
	\arrow[two heads, from=3-6, to=3-4]
	\arrow[from=3-6, to=4-5]
	\arrow[two heads, from=4-3, to=4-1]
	\arrow[two heads, from=4-5, to=4-3]
\end{tikzcd}.
$$
Commutativity implies that for any $n$, the image of the ideal $I$, $W_n(\phi_i) \circ h_n (I) \subset W_n(S_{(a_i)})$, viewed in $S_{a_i}^{[0,n]}$ will be given by the image under the natural projection $S_{a_i}^{\mathbb{N}^E}  \rightarrow S_{a_i}^{[0,n]}$ of $W(\phi_i) \circ h (I) \subset S_{a_i}^{\mathbb{N}^E}$ and this is compatible with the filtration. Then if any element $f$ in $W(\phi_i) \circ h (I) \subset S_{a_i}^{\mathbb{N}^E}$ could be written as a product $g_f \cdot  \gamma_i$ for a fixed $\gamma_i$, where $g_f, \gamma_i \in W(\phi_i) \circ h (I) \subset S_{a_i}^{\mathbb{N}^E}$, the same fact would hold for each $n$ and be compatible with the filtration. Further, if $\gamma_i$ is a non-zero-divisor in $S_{a_i}^{\mathbb{N}^E}$ (hence also in $W(\phi_i) \circ h (I)$) then it is given by components that are non-zero-divisors in $S_{a_i}$. It follows that for any $n$, the image of $\gamma_i$ given by projecting onto the $[0, n]$ components gives a non-zero-divisor in $S_{a_i}^{[0,n]}$ and hence $W_n(\phi_i) \circ h_n (I)$. For each $i$ the ideal generated by the image of $I$, $\langle W_n(\phi_i) \circ h_n (I) \rangle \subset W_n(S_{a_i})$, is then generated by a single non-zero-divisor, the image of $\gamma_i$, and so the pullback $X \times _Z W^{\ast}_n(B_Z X)$ is an effective Cartier divisor. Thus by the universal property of the blow-up we have
$$
\begin{tikzcd}
	{W^{\ast}_n(B_Z X) \times_X Z} && {W^{\ast}_n(B_Z X)} \\
	{E_Z X} && {B_Z X} \\
	Z && X
	\arrow[from=1-1, to=1-3]
	\arrow[from=2-1, to=2-3]
	\arrow[from=1-3, to=2-3]
	\arrow[from=2-3, to=3-3]
	\arrow[from=2-1, to=3-1]
	\arrow[from=3-1, to=3-3]
	\arrow[from=1-1, to=2-1]
\end{tikzcd}
$$
given by unique morphisms
$$
W^{\ast}_n(B_Z X) \rightarrow B_Z X, \quad W^{\ast}_n(B_Z X) \times _X  Z \rightarrow E_Z X
$$
that factor $W^{\ast}_n(B_Z X) \rightarrow W^{\ast}_n(X) \rightarrow X$ and $W^{\ast}_n(B_Z X) \times_X Z \rightarrow Z$. But by uniqueness these morphisms must be compatible with filtration, i.e. for all $n, j$ we have factorizations 
\begin{align*}
   &\dots \rightarrow W^{\ast}_n(B_Z X) \rightarrow W^{\ast}_{n+j}(B_Z X) \rightarrow \dots \rightarrow B_Z X \rightarrow X 
   \\& \dots \rightarrow W^{\ast}_n(B_Z X) \times_Z X \rightarrow W^{\ast}_{n+j}(B_Z X) \times_Z X \rightarrow \dots \rightarrow E_Z X \rightarrow Z
\end{align*}
and we get an induced morphism $W^{\ast}(B_Z X) \rightarrow B_Z X$ such that the composition with the blow-up morphism $W^{\ast}(B_Z X) \rightarrow B_Z X \rightarrow X$ agrees with $W^{\ast}(B_Z X) \rightarrow W^{\ast}(X) \rightarrow X$. Similarly, we have $W^{\ast}_n(E_Z X) \rightarrow W^{\ast}_n(Z) \rightarrow Z \rightarrow X$, which agrees with $W^{\ast}_n(E_Z X) \rightarrow W^{\ast}_n(B_Z X) \rightarrow B_Z X \rightarrow X$ by the commutativity of the partially filled in blow-up square and is compatible with filtration. Thus we get $W^{\ast}_n (E_Z X) \rightarrow W^{\ast}_n(B_Z X) \times_Z X$ compatible with filtration and hence $W^{\ast}_n (E_Z X) \rightarrow E_Z X$ compatible with filtration. Finally, this induces $W^{\ast}(E_Z X) \rightarrow E_Z X$, which commutes with the rest of the diagram. 
\\
We now find the non-zero-divisor $\gamma_i \in W(\phi_i) \circ h (I) \subset S_{a_i}^{\mathbb{N}^E}$. Fix $i$ and let $x \in I, n \in \mathbb{N}^{(E)}$ so that $\psi_n(x) \subset I^{n'}$. $h(x) \in W(A)$ has $n$-th component $\psi_n(x) \in A$, which implies $W(\phi_i)(h(x))$ has $n$-th component 
$
\frac{\psi_n(x)}{1} \sim \frac{(a_i)^{n'}}{1} \cdot \frac{\psi_n(x)}{(a_i)^{n'}} = \frac{\psi_n(a_i)}{1} \cdot \frac{\psi_n(x)}{(a_i)^{n'}} \in S_{(a_i)}.
$
Then for all $x \in I$, $W(\phi_i)(h(x)) = W(\phi_i)(h(a_i)) \cdot (x_n) \in (S_{(a_i)})^{\mathbb{N}^{(E)}}$, where $(x_n) \in (S_{(a_i)})^{\mathbb{N}^{(E)}}$ is the element with $n$-th component $\frac{\psi_n(x)}{(a_i)^{n'}}$. Now if $(x_n) \in W(S_{(a_i)})$ then this would imply $W(\phi_i)(h(x)) \in \langle W(\phi_i)(h(a_i)) \rangle \subset W(S_{(a_i)})$, where $\langle W(\phi_i)(h(a_i)) \rangle$ is the principal ideal generated by $W(\phi_i)(h(a_i))$ and we would have $\langle W(\phi_i)(h(I)) \rangle =\langle W(\phi_i)(h(a_i)) \rangle$. But note that for all $n \in \mathbb{N}^{(E)}$ the $n$-th component of $W(\phi_i)(h(a_i))$ is $\frac{(a_i)^{n'}}{1} \in S_{(a_i)}$, a non-zero-divisor in $S_{(a_i)}$. Therefore $\gamma_i = W(\phi_i)(h(a_i))$ is a non-zero-divisor in $(S_{(a_i)})^{\mathbb{N}^{(E)}}$ and hence a non-zero-divisor in $W(S_{(a_i)})$ satisfying the desired conditions. It remains to show $(x_n) \in (S_{(a_i)})^{\mathbb{N}^{(E)}}$ is actually in $W(S_{(a_i)})$. We prove a slightly more general fact below, which finishes the proof. \qedsymbol

\begin{lemma}
  For all $k=1,2, \dots $ and any $x \in I^k$, the element $(x_n) \in (S_{(a_i)})^{\mathbb{N}^{(E)}}$ with $n$-th component $\frac{\psi_n(x)}{(a_i)^{kn'}}$ is contained in $W(S_{(a_i)})$.  
\end{lemma}

\textit{Proof.} We show by induction that for all $i$, $(x_n) \in U_i(S_{(a_i)})$ for $U_i(S_{(a_i)})$ as defined in \cite[Section 1]{borger2015basicI} and hence $(x_n) \in W(S_{a_i})$. Clearly, any $(x_n)$ as above is in $U_0(S_{(a_i)}) = (S_{(a_i)})^{\mathbb{N}^{(E)}}$ so assume all such $(x_n)$ are elements of $U_i(S_{(a_i)})$. Let $(x_n)$ be any such element (for specific $x \in I^k$) and let $\Psi_p$ be the commuting endomorphisms describing the $\Psi_{\mathbb{Z}, E}$-ring structure on $(S_{(a_i)})^{\mathbb{N}^{(E)}}$. Then for all $n \in \mathbb{N}^{(E)}$, the $n$-th component of $\Psi_p((x_n)) - (x_n)^p$ will be given by 
\begin{align*}
  &((n+p)\textrm{-th component of }(x_n))-(n-\textrm{th component of }(x_n))^p  
  \\&=\frac{\psi_p(\psi_n(x))}{(a_i)^{kpn'}}- \left( \frac{\psi_n(x)}{(a_i)^{kn'}}  \right)^p =  \frac{\psi_p(\psi_n(x)) -\psi_n(x)^p}{(a_i)^{pkn'}}   
  \\&= \frac{\psi_n(\psi_p(x)-x^p)}{(a_i)^{pkn'}}=\frac{p\psi_n(h)}{(a_i)^{pkn'}} = p \left( \frac{\psi_n(h)}{(a_i)^{pkn'}} \right),
\end{align*}
where by $\psi_p$ being a Frobenius lift on $A$ we have $\psi_p(x)-x^p=ph \in I^{pk}$ and $ h \in I^{pk}$ by the $\mathbb{Z}$-flatness of $\frac{A}{I^{pk}}$. Thus $\Psi_p((x_n)) - (x_n)^p = p(h_n)$, where $(h_n) \in (S_{(a_i)})^{\mathbb{N}^{(E)}}$ has $n$-th component $\frac{\psi_n(h)}{(a_i)^{pkn'}}$. But by induction $(h_n) \in U_i(S_{(a_i)})$ and thus $(x_n) \in U_{i+1}(S_{(a_i)})$. \qedsymbol

\begin{proposition}
    \label[proposition]{proposition2124}
    The morphisms 
    $$
    W^{\ast}(B_Z X) \rightarrow B_Z X, \quad W^{\ast}(E_Z X) \rightarrow E_Z X
    $$
    of \cref{proposition2122} make $B_Z X, E_Z X$ algebras for $W^{\ast}$ and so the blow-up square is $\Lambda$-equivariant.
\end{proposition}

\textit{Proof.} As the  morphism $W^{\ast}(E_Z X) \rightarrow E_Z X$ is induced from a pullback if the morphisms $W^{\ast}_n(B_Z X) \rightarrow B_Z X$ are compatible with the monad structure of $W^{\ast}$, then the morphisms $W^{\ast}_n(E_Z X) \rightarrow E_Z X$ will be compatible also. Hence as the diagrams of \cite[10.6]{borger2015basicII} induce the monad structure, $W^{\ast}(B_Z X) \rightarrow B_Z X, W^{\ast}(E_Z X) \rightarrow E_Z X$ will induce a $\Lambda$-structure on $B_Z X, E_Z X$. Thus it suffices to show $W^{\ast}_n(B_Z X) \rightarrow B_Z X$ is compatible with the monad structure of $W^{\ast}$, which we will argue locally using \cite[2.4]{borger2015basicI}.  
Locally the scheme morphisms $W^{\ast}_n(B_Z X) \rightarrow B_Z X$ will correspond to morphisms $S_{a_i} \rightarrow W_n(S_{a_i})$ and via the blow-up property we have a morphism $S_{a_i} \rightarrow W(S_{a_i})$. By limit properties for all $n, n+j$ we have
$$
S_{a_i} \rightarrow W(S_{a_i}) \rightarrow W_{n+j}(S_{a_i}) \rightarrow W_{n}(S_{a_i}),
$$
where $S_{a_i} \rightarrow W(S_{a_i})$ and the compositions $S_{a_i} \rightarrow W(S_{a_i}) \rightarrow W_{n+j}(S_{a_i})$ are induced from the blow-up. If the morphism $S_{a_i} \rightarrow W(S_{a_i})$ makes $S_{a_i}$ a co-algebra for $W$ then this implies the morphisms $S_{a_i} \rightarrow W_{n}(S_{a_i})$ are compatible with the diagrams defined by \cite[2.4]{borger2015basicI}. This will imply that the scheme morphisms $W^{\ast}_n(B_Z X) \rightarrow B_Z X$ are compatible with the corresponding diagrams of \cite[10.6]{borger2015basicII} as desired. It remains to show the induced map $S_{a_i} \rightarrow W(S_{a_i})$ makes each $S_{a_i}$ a coalgebra for the comonad $W$. 
The induced ring morphism is
$$
h_i: S_{(a_i)} \rightarrow W(S_{(a_i)}), \quad \frac{x}{(a_i)^d} \mapsto \frac{W(\phi_i)(h(x))}{(W(\phi_i)(h(a_i)))^d},
$$
where $x \in I^d$ implies $W(\phi_i)(h(x)) \in \langle (W(\phi_i)(h(a_i)))^d \rangle$ and we divide off by 
$
(W(\phi_i)(h(a_i)))^d.
$
For each $n \in \mathbb{N}^{(E)}$, $\psi_n(x) \in I^{dn'}$, which implies the $n$-th component of $W(\phi_i)(h(x))$ is
$
\frac{\psi_n(x)}{1} \sim \frac{(a_i)^{dn'}}{1} \cdot \frac{\psi_n(x)}{(a_i)^{dn'}} = \frac{(\psi_n(a_i))^d}{1} \cdot \frac{\psi_n(x)}{(a_i)^{dn'}} \in S_{(a_i)}.
$
Dividing off by $(W(\phi_i)(h(a_i)))^d$ gives $\frac{\psi_n(x)}{(a_i)^{dn'}}$, i.e. the $n$-th component of $h_i\left (  \frac{x}{(a_i)^d}  \right)$ is $\frac{\psi_n(x)}{(a_i)^{dn'}}$. Note this gives a well-defined morphism as argued in the $n=p$ case in the construction of Frobenius lifts of $S_{(a_i)}$ in the proof of \cref{proposition2113}.
\\We now show $h_i$ commutes with the diagrams defining the comonad structure of $W^{fl}$ and so induces a coalgebra structure on $S_{(a_i)}$, completing the proof. 
\\Let $\frac{x}{(a_i)^d} \in S_{(a_i)}$. The $0$-th component of $h_i \left( \frac{x}{(a_i)^d} \right)$ is $\frac{\psi_0(x)}{(a_i)^{d}} = \frac{x}{(a_i)^{d}}$, i.e. 
\\$
\eta_{S_{(a_i)}}(h_i \left( \frac{x}{(a_i)^d} \right)) = \left( \frac{x}{(a_i)^d} \right),
$
showing the first required diagram commutes. For the second to commute we want to show 
\begin{align*}
    W(h_i)(h_i \left( \frac{x}{(a_i)^d} \right)) = \mu_{S_{(a_i)}}(h_i \left( \frac{x}{(a_i)^d} \right)) \in W^2(S_{(a_i)}) &\subset (W(S_{(a_i)}))^{\mathbb{N}^{(E)}} 
    \\&\subset ((S_{(a_i)})^{\mathbb{N}^{(E)}})^{\mathbb{N}^{(E)}}.
\end{align*}
It suffices to show this on each $n$-th component, $n \in \mathbb{N}^{(E)}$, i.e. 
$$
\left( W(h_i)\left(h_i \left( \frac{x}{(a_i)^d} \right) \right) \right)_n = \left(\mu_{S_{(a_i)}}\left(h_i \left( \frac{x}{(a_i)^d} \right) \right) \right)_n \in W(S_{(a_i)}) \subset (S_{(a_i)})^{\mathbb{N}^{(E)}},
$$
for which it suffices to show on each $m$-th component, $m \in \mathbb{N}^{(E)}$, i.e.
$$
\left(\left( W(h_i)\left(h_i \left( \frac{x}{(a_i)^d} \right) \right) \right)_n\right)_m = \left (\left(\mu_{S_{(a_i)}}\left(h_i \left( \frac{x}{(a_i)^d} \right) \right) \right)_n \right)_m \in S_{(a_i)}.
$$
Now the $n$-th component of $\mu_{S_{(a_i)}}\left(h_i \left( \frac{x}{(a_i)^d} \right) \right)$ is $\Psi_n \left( h_i \left( \frac{x}{(a_i)^d} \right)  \right)$, which has $m$-th component the $(n+m)$-th component of $h_i \left( \frac{x}{(a_i)^d} \right)$; $\frac{\psi_{n+m}(x)}{(a_i)^{dn'm'}} = \frac{\psi_{n}(\psi_m(x))}{(a_i)^{dn'm'}} = \frac{\psi_{m}(\psi_n(x))}{(a_i)^{dn'm'}}$. The $n$-th component of $W(h_i)\left(h_i \left( \frac{x}{(a_i)^d} \right) \right)$ is 
\\$h_i(n\textrm{-th component of }h_i \left( \frac{x}{(a_i)^d} \right)) = h_i \left(\frac{\psi_{n}(x)}{(a_i)^{dn'}} \right)$, which has $m$-th component $\frac{\psi_m(\psi_{n}(x))}{(a_i)^{dmn'}}$. \qedsymbol

\begin{proposition}
    \label[proposition]{proposition2126} 
    The two $\Lambda$-scheme structures on $B_Z X$, the first given by explicitly writing out Frobenius lifts (\cref{example2115}) and the second constructed via blow-up universal properties (\cref{proposition2122}), are the same.
\end{proposition}

\textit{Proof.} We show that the morphisms $W^{\ast}_n(B_Z X) \rightarrow B_Z X$ induced by $\Lambda$-structure from explicit Frobenius lifts in the proof of \cref{proposition2113} agree with those constructed from the universal property. This will show that the algebra structures for $W^{\ast}$ and hence $\Lambda$-structures are the same. For each $i,j$ and $S_{(a_i)}, S_{(a_i)}, S_{(a_ia_j)}$ the Frobenius lifts induce coalgebra structures given by 
$
S_{(a_i)} \rightarrow W(S_{(a_i)}), \quad S_{(a_ia_j)} \rightarrow W(S_{(a_ia_j)}),
$
which will be compatible with the filtration and so we have the following commutative diagram
$$
\begin{tikzcd}[cramped]
	{W_{n}(S_{a_i})} && {W_{n+j}(S_{a_i})} && {W(S_{a_i})} && {S_{a_i}} \\
	\\
	{W_{n}(S_{a_ia_j})} && {W_{n+j}(S_{a_ia_j})} && {W(S_{a_ia_j})} && {S_{a_ia_j}}
	\arrow[from=1-1, to=3-1]
	\arrow[two heads, from=1-3, to=1-1]
	\arrow[from=1-3, to=3-3]
	\arrow[two heads, from=1-5, to=1-3]
	\arrow[from=1-5, to=3-5]
	\arrow[from=1-7, to=1-5]
	\arrow[from=1-7, to=3-7]
	\arrow[two heads, from=3-3, to=3-1]
	\arrow[two heads, from=3-5, to=3-3]
	\arrow[from=3-7, to=3-5]
\end{tikzcd}.
$$
This induces the morphisms 
$
W^{\ast}_n(B_Z X) \rightarrow B_Z X, \quad W^{\ast}_{n+j}(B_Z X) \rightarrow B_Z X,
$
which have the necessary compatibility conditions to induce a map $W^{\ast}(B_Z X) \rightarrow X$ making $B_Z X$ an algebra for $W^{\ast}$. Composing with the blow-up morphism $B_Z X \rightarrow X$ gives $W^{\ast}_n(B_Z X) \rightarrow B_Z X \rightarrow X$, which over $\textrm{Spec} (S_{(a_i)})$ is 
$$
\textrm{Spec} (W_n(S_{(a_i})) \rightarrow \textrm{Spec} (S_{(a_i)}) \rightarrow \textrm{Spec} (A)
$$
corresponding to 
$$
A \rightarrow S_{(a_i)} \rightarrow W_n(S_{(a_i)}) \subset S_{(a_i)}^{[0,n]}, \quad a \mapsto \frac{a}{1} \mapsto \left(\frac{\psi_n(a)}{1}  \right)_{n'} \space \forall \, n' \in [0,n],
$$
which agrees with 
$$
W_n(\phi_i)(h_n): A \rightarrow W_n(A) \rightarrow W_n(S_{(a_i)}), \quad a \mapsto (\psi_n(a))_{n'} \mapsto \left(\frac{\psi_{n}(a)}{1}  \right)_{n'} \space \forall \, n' \in [0,n].
$$
Thus each morphism $W^{\ast}_n(B_Z X) \rightarrow B_Z X \rightarrow X$ glued from Frobenius lifts factors the $W^{\ast}_n(B_Z X) \rightarrow X$ coming from the universal property. By the uniqueness of the universal property the glued Frobenius morphisms agree with those previously constructed. Finally, by the uniqueness of maps to the colimit we have the same map $W^{\ast}(B_Z X) \rightarrow X$ in $\textbf{Sp}$ and so the same $\Lambda$-structure. \qedsymbol

\begin{corollary}
    The two $\Lambda$-blow-up squares we constructed are the same.
\end{corollary}
\textit{Proof.} As the $\Lambda$-structure on the blow-up is the same and we induced the $\Lambda$-structure on $E_Z X$ by pullbacks we get the same $\Lambda$-structure on $E_Z X$ also. \qedsymbol

\begin{corollary}
    \label[corollary]{corollary2128} 
    Let $A$ be a $\mathbb{Z}$-flat $\Lambda$-ring and suppose $I=(x_1, \dots, x_n)=(y_1, \dots, y_m)$. Assume both sets of generators are toric, i.e. $\psi_p(x_i)=(x_i)^p$ and $\psi_p(y_i)=(y_i)^p$ and $I$ has $\mathbb{Z}$-flat normal cone as in \cref{example2107}. Then the two $\Lambda$-structures on the blow-up $B_Z X$ as in \cref{example2115} using $y_i$ or $x_i$ are the same.
\end{corollary}

\textit{Proof.} By \cref{proposition2126} the $\Lambda$-structures on $B_Z X$ coming from $y_i$ or $x_i$ are given as in \cref{proposition2122} by $W^{\ast}_n(B_Z X) \times_X Z$ being an effective Cartier divisor. This induces via the universal property a unique morphism $W^{\ast}_n(B_Z X) \rightarrow B_Z X$ that does not depend on the toric generators of $I$. \qedsymbol




\subsection{General case} 
We would like to extend to blow-ups of general $\mathbb{Z}$-flat $\Lambda$-schemes that over affine opens recover the $\Lambda$-structure of our affine blow-ups. However, to generalise the constructions of the previous section we require $\Lambda$-schemes with $\Lambda$-equivariant affine open covers, which generally need not exist (see \cite[2.6]{borger2009lambdarings}). 
\\
\\\textbf{$\Lambda$-equivariant localizations and blow-up algebras.} Let $A$ be a $\mathbb{Z}$-flat $\Lambda$-ring and let $\psi_{p,A}$ be the Frobenius lifts defining the $\Lambda$-structure. Let $A_f$ be the localization of $A$ at $f$ underlying a distinguished open of $\textrm{Spec} \, A$. Assume that $\textrm{Spec} (A_f)=D(f) \subset \textrm{Spec} \, A$ is also a $\Lambda$-scheme, writing $\psi_{p, A_f}$ for its Frobenius lifts. Write $i:A \rightarrow A_f$ for the canonical localization and suppose it is $\Lambda$-equivariant so
$
i(\psi_{p, A}(a)) = \frac{\psi_{p, A}(a)}{1} = \psi_{p, A_f}\left(\frac{a}{1}\right)
$
and $D(f) \subset \textrm{Spec} \, A$ is a $\Lambda$-equivariant open immersion.

\begin{example}
    Let $f \in A$ be \textit{toric} so $\psi_p(f)=f^p$ for all $p$. Then the natural induced $\Lambda$-structure on $A_f$ is given by
    $
    \psi_{p, A_f}\left( \frac{a}{f^n}  \right) = \frac{\psi_{p, A}}{f^{pn}},
    $
    which are well-defined by the arguments in the proof of \cref{proposition2113}. Then $D(f) \subset \textrm{Spec} \, A$ is a $\Lambda$-equivariant open immersion.
\end{example}

Let $I=(x_1, \dots, x_n) \subset A$ be an ideal as in \cref{example2115} for which $\psi_p(x_i)=(x_i)^p$ and a $\Lambda$-structure on $\textrm{Proj} (A[It])$ is given by Frobenius lifts on a $x_i$-chart $\psi_{p, A, i}$. Let $f \in A$ such that $D(f) \subset \textrm{Spec} \, A$ is a $\Lambda$-equivariant open immersion. Now $I_f$ is generated by $\frac{x_i}{1}$, $\psi_{p, A_f}\left( \frac{x_i}{1}  \right) = \left( \frac{x_i}{1}  \right)^p$ and so we can induce a $\Lambda$-structure on $\textrm{Proj} (A_f[I_f t])$ with Frobenius lifts on a $x_i$-chart denoted by $\psi_{p, A_f, i}$.

\begin{lemma}
     The open immersion $\textrm{Proj} (A_f[I_f t]) \hookrightarrow \textrm{Proj} (A[It])$ is $\Lambda$-equivariant.
\end{lemma}
\textit{Proof.} Recall $\textrm{Proj} (A_f[I_f t])$ will be an open subscheme of $\textrm{Proj} (A[It])$, on a $x_i$-chart induced by the ring morphism
$$
A\left[ \frac{I}{x_i}\right] \rightarrow A_f\left[ \frac{I_f}{\frac{x_i}{1}}\right], \quad \frac{a}{x_i^n} \mapsto \frac{\frac{a}{1}}{\left(\frac{x_i}{1}\right)^n} ,
$$
which is equivalent to localization at $\frac{f}{1} \in A[\frac{I}{x_i}]$ in degree $0$.
Then 
$$
\psi_{p, A, i} \left(  \frac{a}{x_i^n} \right) = \left(  \frac{\psi_{p, A}(a)}{x_i^{pn}} \right) \mapsto \frac{\frac{\psi_p(a)}{1}}{\left(\frac{x_i}{1}\right)^{pn}} = \frac{\psi_{p, A_f}\left(\frac{a}{1}\right)}{\left(\frac{x_i}{1}\right)^{pn}} = \psi_{p, A_f, i} \left(\frac{\frac{a}{1}}{\left(\frac{x_i}{1}\right)^n} \right)
$$
and so the open immersion will be $\Lambda$-equivariant. \qedsymbol
\\
\\We combine the above results to define the schemes on which we can induce a $\Lambda$-structure on blow-ups.

\begin{definition}[Affine $\Lambda$-covers]
    \label[definition]{AffineLambdaCover} 
    Let $X$ be a $\mathbb{Z}$-flat $\Lambda$-scheme such that there exists an affine open cover $X=\bigcup_i \textrm{Spec} \, A_i$ given by affine $\Lambda$-schemes such that intersections $\textrm{Spec} \, A_i \cap \textrm{Spec} \, A_j$ can be covered by distinguished affine open $\Lambda$-schemes and the canonical localizations are given by $\Lambda$-equivariant ring morphisms as in the construction above. If $X$ is separated then we can just take the intersection. Moreover, we assume that the $\Lambda$-structure on $X$ comes from the colimit in $\Lambda$-schemes taken over this affine system. In particular, the open embeddings $\textrm{Spec} (A_i) \rightarrow X$ will be $\Lambda$-equivariant. We call such an affine open cover an \textit{affine $\Lambda$-cover} for $X$
\end{definition}

\begin{example}[Monoid schemes]
    \label[example]{example2204} 
   $\mathbb{Z}$-realizations of monoid schemes have affine $\Lambda$-covers given by the $\mathbb{Z}$-realizations of affine open monoid subschemes. 
\end{example}


\begin{definition}
     Let $X$ be a $\mathbb{Z}$-flat $\Lambda$-scheme with an affine $\Lambda$-cover and $\mathcal{I} \subset \mathcal{O}_X$ a quasicoherent ideal sheaf. $\mathcal{I}$ is \textit{locally toric} for the cover of $X$ if restricted to each $\textrm{Spec} (A_i)$, $\mathcal{I}$ is given by an ideal $I_i=(x_1, \dots, x_{n_i})$ such that $\psi_{p, i}(x_k)=(x_k)^p$ for all $p$ and $k$ where $\psi_{p, i}$ are the Frobenius lifts of $A_i$. We call the $x_k$ \textit{toric} generators of $I_{i}$.
\end{definition}

\begin{remark}
    Let $\mathcal{I}$ be locally toric for an affine $\Lambda$-cover of $X$. Then on a distinguished open $\textrm{Spec} (A_{ij})$ in the intersection $\textrm{Spec} (A_i) \cap \textrm{Spec} (A_j)$, $\mathcal{I}(\textrm{Spec} (A_{ij}))$ will be locally toric. However, it may have different sets of toric generators coming from the localizations of $I_{i}$ and $I_{j}$.
\end{remark}

\begin{example}
    Let $X_{\mathbb{Z}}$ be the $\mathbb{Z}$-realization of a monoid scheme $X$ and $Z \subset X$ an equivariant closed monoid subscheme of $X$. Then the ideal sheaf defining the closed subscheme $Z_{\mathbb{Z}} \subset X_{\mathbb{Z}}$ will be locally toric for any affine $\Lambda$-cover of $X_{\mathbb{Z}}$ coming from an affine monoid scheme cover of $X$.
\end{example}

\begin{proposition}[Gluing blow-up $\Lambda$-structures]
    \label[proposition]{proposition2208} 
    Let $X$ be a $\mathbb{Z}$-flat $\Lambda$-scheme with an affine $\Lambda$-cover. Let $Z \subset X$ be a closed subscheme defined by $\mathcal{I}\subset \mathcal{O}_X$ such that the normal cone of $\mathcal{I}$ is $\mathbb{Z}$-flat and $\mathcal{I}$ is locally toric for the cover. Then the blow-up $B_Z X$ is a $\mathbb{Z}$-flat $\Lambda$-scheme with an affine $\Lambda$-cover and $\pi: B_Z X \rightarrow X$ is $\Lambda$-equivariant.
\end{proposition}
\textit{Proof.} Let $\textrm{Spec} A_i$ be an affine open in the $\Lambda$-cover and $I_i=\mathcal{I}(\textrm{Spec} \, A_i)$. As $\mathcal{I}$ has $\mathbb{Z}$-flat normal cone and is locally toric for the cover, $I_i$ is as in \cref{example2115} and $\textrm{Proj} (A_i[I_i t])$ is a $\Lambda$-scheme. Recall $B_Z X$ can be covered by $\textrm{Proj} (A_i[I_i t])$, which can be covered by affine open $\textrm{Spec} (A_i\left[\frac{I_i}{x_{j}}  \right])$ for toric generators $x_{j}$ of $I_{i}$. Let $\textrm{Spec} (A_{ii'})$ be a distinguished open of the $\Lambda$-cover in the intersection $\textrm{Spec} (A_{i'}) \cap \textrm{Spec} (A_i)$. There are toric generators $x_1, \dots, x_{n(i)}$ of $I_{ii'}$ coming from $I_i$, which induce affine open $\textrm{Spec} (A_{ii'}\left[\frac{I_{ii'}}{x_j}  \right])$ of $\textrm{Proj} (A_{ii'}[I_{ii'} t])$ that are distinguished open $\Lambda$-subschemes of $\textrm{Spec} (A_i\left[\frac{I_i}{x_j}  \right])$ as in the construction above. Thus we can induce a $\Lambda$-equivariant open embedding 
$
\textrm{Proj} (A_{ii'}[I_{ii'} t]) \rightarrow \textrm{Proj} (A_{i}[I_{i} t]).
$
Similarly, we have toric generators $y_1, \dots, y_{n(i')}$ of $I_{ii'}$ coming from $I_{i'}$ inducing $\Lambda$-equivariant 
$
\textrm{Proj} (A_{ii'}[I_{ii'} t]) \rightarrow \textrm{Proj} (A_{i'}[I_{i'} t]).
$
But by \cref{corollary2128} the $\Lambda$-structure on $\textrm{Proj} (A_{ii'}[I_{ii'} t])$ is independent of the toric generators and the $\textrm{Proj} (A_{ii'}[I_{ii'} t])$ cover the intersection $\textrm{Proj} (A_i[I_i t]) \cap \textrm{Proj} (A_{i'}[I_{i'} t])$. Thus it follows that we can glue along the $\textrm{Proj} $'s a well-defined $\Lambda$-structure on $B_Z X$. Moreover, as the blow-up morphisms 
\\$\textrm{Proj} (A_{i}[I_{i} t]) \rightarrow \textrm{Spec} (A_i)$ are $\Lambda$-equivariant, $\pi: B_Z X \rightarrow X$ will be $\Lambda$-equivariant also. For the affine $\Lambda$-cover of $B_Z X$ we consider the affine open $\textrm{Spec} (A_i\left[\frac{I_i}{x_{j}}  \right])$ covering each $\textrm{Proj} (A_i[I_i t])$. Recall we have shown in the affine case that the intersection of a $\textrm{Spec} (A_i\left[\frac{I_i}{x_j}  \right])$ and $\textrm{Spec} (A_i\left[\frac{I_i}{x_{j'}}  \right])$ within a $\textrm{Proj} (A_i[I_i t])$ is an affine open $\Lambda$-scheme given by a localization. Consider the intersection $\textrm{Spec} (A_i\left[\frac{I_i}{x_j}  \right]) \cap \textrm{Spec} (A_{i'}\left[\frac{I_{i'}}{y_{j'}}  \right])$. Write $x_j$ and $y_{j}'$ for the toric generators of $I_{ii'}$ in $A_{ii'}$ induced by the respective localizations. Then $x_jy_{j'}$ will be toric, i.e. $\psi_p(x_jy_{j'})=(x_jy_{j'})^p$ and so $A_{ii'}\left[\frac{I_{ii'}}{(x_jy_{j'})}  \right]$ will have a $\Lambda$-structure given by Frobenius lifts as in the affine construction. This implies we have
$
A_i\left[\frac{I_i}{x_j}  \right] \rightarrow A_{ii'}\left[\frac{I_{ii'}}{(x_j)}  \right] \rightarrow A_{ii'}\left[\frac{I_{ii'}}{(x_jy_{j'})}  \right],
$
where each morphism is $\Lambda$-equivariant and given by the localization at a section. It follows that 
$
\textrm{Spec} (A_{ii'}\left[\frac{I_{ii'}}{(x_jy_{j'})}  \right]) \rightarrow \textrm{Spec} (A_i\left[\frac{I_i}{x_j}  \right])
$
induced by the composition will be a $\Lambda$-equivariant open immersion given by a distinguished open and similarly for the embedding into $\textrm{Spec} (A_{i'}\left[\frac{I_i}{y_{j'}}  \right])$. Then as 
$
\textrm{Spec} (A_{ii'}\left[\frac{I_{ii'}}{(x_jy_{j'})}  \right]) \textrm{ corresponds to } D(x_{j}) \cap D(y_{j'}) \subset \textrm{Proj} (A_{ii'}[I_{ii'} t])
$
and $\bigcup D(x_{j}), \bigcup D(y_{j'})$ cover $\textrm{Proj} (A_{ii'}[I_{ii'} t])$ the $D(x_{j}) \cap D(y_{j'})$'s will cover 
$
\textrm{Spec} (A_i\left[\frac{I_i}{x_j}  \right]) \cap \textrm{Spec} (A_{i'}\left[\frac{I_{i'}}{y_{j'}}  \right])
$
as $A_{ii'}$ and $x_j, y_{j'}$ vary. \qedsymbol

\begin{example}
    Let $X_{\mathbb{Z}}$ be the $\mathbb{Z}$-realization of a monoid scheme $X$ and $Z \subset X$ an equivariant closed monoid subscheme of $X$. By \cref{proposition2208} $B_{X_{\mathbb{Z}}} Z_{\mathbb{Z}} \cong (B_Z X)_{\mathbb{Z}}$ is a $\mathbb{Z}$-flat $\Lambda$-scheme. Moreover, it follows from \cref{example2116} that the $\Lambda$-structure is precisely the toric $\Lambda$-structure.
\end{example}

\begin{lemma}
   Let $X$ be a $\mathbb{Z}$-flat $\Lambda$-scheme and suppose we have $\{U_i \}, \{V_j \}$ two affine $\Lambda$-covers on $X$. Let $Z$ be a closed subscheme with corresponding ideal sheaf $\mathcal{I}$ with $\mathbb{Z}$-flat normal cone and assume that $\mathcal{I}$ is locally toric for $\{U_i\}$. Then $B_{Z_j} V_j$ has a $\Lambda$-structure and $\pi_j: B_{Z_j} V_j \rightarrow V_j$ is $\Lambda$-equivariant, where $Z_j = Z \cap V_j$. 
\end{lemma}
\textit{Proof.} As $\mathcal{I}$ is locally toric for $\{ U_i \}$ we have a $\Lambda$-structure on $B_Z X$ with $\pi: B_Z X \rightarrow X$ $\Lambda$-equivariant and so we can construct the pullback
$$
\begin{tikzcd}[cramped]
	{B_Z X} && X \\
	\\
	{B_{Z_j} V_j} && {V_j}
	\arrow["\pi", from=1-1, to=1-3]
	\arrow[from=3-1, to=1-1]
	\arrow["{\pi_j}", from=3-1, to=3-3]
	\arrow[from=3-3, to=1-3]
\end{tikzcd}.
$$
As $V_j \rightarrow X$ and $\pi$ are $\Lambda$-equivariant we can take the pullback in $\Lambda$-schemes, which agrees with the pullback in schemes. Thus $B_Z V_j$ has a unique $\Lambda$-structure compatible with the pullback and $\pi_j$ is $\Lambda$-equivariant. \qedsymbol

\begin{corollary}
    If $\mathcal{I}$ is also locally toric for $\{V_j \}$ then the $\Lambda$-structure on $B_Z X$ induced from $\{V_j \}$ agrees with that induced from $\{U_i\}$.
\end{corollary}

\textit{Proof.} $\mathcal{I}$ being locally toric for $\{V_j\}$ implies we can induce on $B_{Z_j} V_j$ a $\Lambda$-structure as in the proof of \cref{proposition2208}. But this structure will be compatible with the pullback diagram in the lemma and so must be the same $\Lambda$-structure. It follows that the $\Lambda$-structure on $B_Z X$ induced from the $\{V_j \}$, which comes from gluing the $B_{Z_j} V_j$ must be the same as the $\Lambda$-structure from the $\{U_i\}$. \qedsymbol



\subsection{Strict transforms}

We briefly describe the strict transforms of closed $\Lambda$-subschemes with affine $\Lambda$-covers.

\subsubsection{Affine case}
Let $A, B$ be $\mathbb{Z}$-flat $\Lambda$-rings with Frobenius lifts $\psi_{p, A}, \psi_{p, B}$ and let $f:A \rightarrow B$ be a $\Lambda$-equivariant morphism. Suppose $x_1, \dots, x_n \in A$ are toric so $\psi_{p, A}(x_i)=(x_i)^p$ and $f(x_i)$ are toric on $B$. Let $I=(x_1, \dots, x_n) \subset A$ and $J:=BI=(f(x_1), \dots, f(x_n)) \subset B$ be the ideal in $B$ generated by the image of $I$. Suppose both $I, J$ have $\mathbb{Z}$-flat normal cones so that assuming $f(x_i)$ is non-zero we have an induced $\Lambda$-structure on $A[\frac{I}{x_i}], B[\frac{J}{f(x_i)}]$.

\begin{lemma}
    \label[lemma]{lemma2311} 
    The induced morphism
$$
f_i: A\left[\frac{I}{x_i}\right] \rightarrow B\left[\frac{J}{f(x_i)}\right], \quad \frac{a}{(x_i)^n} \mapsto \frac{f(a)}{(f(x_i))^n}
$$
is $\Lambda$-equivariant.
\end{lemma}
\textit{Proof.} 
\begin{align*}
    f_i(\psi_{p, A, i}\left(\frac{a}{(x_i)^n} \right))&=f_i(\left(\frac{\psi_{p, A}(a)}{(x_i)^{pn}} \right))=\frac{f(\psi_{p, A}(a))}{(f(x_i))^{pn}} 
    \\&=\frac{\psi_{p, B}(f(a))}{(f(x_i))^{pn}}=\psi_{p, B, i}\left(\frac{f(a)}{(f(x_i))^n} \right)=
    \psi_{p, B, i}(f_i\left(\frac{a}{(x_i)^n} \right)). \quad \qed
\end{align*}

\begin{remark}
    In general, even if $I$ has $\mathbb{Z}$-flat normal cone and induces a $\Lambda$-structure on the blow-up, $J$ need not. We will give quite specific conditions on when a closed embedding of a smooth $\Lambda$-scheme will satisfy this condition, which will be used to define $\Lambda$-marked monomial ideals in \cref{C:5}. 
\end{remark}


\subsubsection{General case}
Let $X, Y$ be $\mathbb{Z}$-flat $\Lambda$-schemes with affine open $\Lambda$-covers. Let $f: Y \rightarrow X$ be a $\Lambda$-equivariant affine scheme morphism, which is given by a $\Lambda$-ring morphism on affine opens in the cover of $X$. Suppose that $\mathcal{I}$ is a quasicoherent ideal sheaf on $X$ such that $\mathcal{I}$ is locally toric for the cover. Then $\mathcal{J}:=\mathcal{I}\mathcal{O}_Y$ is a quasicoherent ideal sheaf on $Y$ that is locally toric for the cover of $Y$.

\begin{lemma}
    Suppose $\mathcal{I}, \mathcal{J}$ have $\mathbb{Z}$-flat normal cones. Then writing $Z, W$ for the closed subschemes defined by $\mathcal{I}, \mathcal{J}$ we have a $\Lambda$-equivariant morphism
    $$
    B_W Y \rightarrow B_Z X.
    $$
\end{lemma}
\textit{Proof}. The morphism $ B_W Y \rightarrow B_Z X$ is the natural one induced from the morphism of the affine blow-up algebras as in \cref{lemma2311}. As the $\Lambda$-structures on $B_ZX, B_WY$ are induced from $\Lambda$-structures on the affine blow-up algebras the lemma follows. \qedsymbol

\begin{corollary}[Strict transforms]
    \label[corollary]{corollary2322} 
    Let $Y \hookrightarrow X$ as above be a closed embedding and $\mathcal{I}, \mathcal{J}$ satisfy conditions of the lemma. Then $B_W Y$ is the strict transform of $Y$ via the blow-up $\pi: B_Z X \rightarrow X$ and is a closed $\Lambda$-subscheme with affine $\Lambda$-covers given by $\Lambda$-ring morphisms.
\end{corollary}
\textit{Proof.} Follows from the previous lemma and that the affine morphism on blow-up algebras locally defines the strict transform. \qedsymbol

\section{Smooth blow-ups and divisors over $\mathbb{Z}$}
\label{C:3} 

In this section we describe some properties of blow-ups of smooth schemes over $\textrm{Spec}  \, \mathbb{Z}$. These will be schemes $X$ such that the morphism $X \rightarrow \textrm{Spec}  \, \mathbb{Z}$ is smooth in the sense of \cite[Theorem 25.2.2]{FOAG} (relatively smooth of some dimension in every neighbourhood). We will also call these smooth $\mathbb{Z}$-schemes or smooth schemes over $\mathbb{Z}$. We first show that, as in characteristic $0$, we have local blow-up coordinates, which enable us to locally analyse behaviour under blowing up. Blow-up coordinates will be used throughout the remainder of the paper and are critical in proving order reduction. We then describe simple normal crossings divisors in the $\mathbb{Z}$-flat and $\Lambda$-structure context. These divisors will be used to construct ideals and blow-up centres necessary for our order reduction theorem. For any prime ideal $(p) \in \textrm{Spec}  \, \mathbb{Z}$ we will denote by $X_p$ the fiber over $(p)$. Similarly, for $x \in X$ if $x$ maps to $(p) \in \textrm{Spec}  \, \mathbb{Z}$ then denote by $x_p \in X_p$ the point lying over $x$ such that $\mathcal{O}_{X_p, x_p} \cong \mathcal{O}_{X,x}/(p)$. If $X$ is smooth over $\mathbb{Z}$ then $X_p$ is a smooth $\mathbb{F}_p$-variety if $p \not = 0$ and a smooth $\mathbb{Q}$-variety if $p = 0$. For convenience, we will use $\mathbb{F}_0$ to refer to $\mathbb{Q}$ so that $X_p$ is a smooth $\mathbb{F}_p$-variety for all $(p)$.


\subsection{Local blow-ups}
\label{C:3_Sec:Local blow-ups} 
We first recall how to compare the local blow-ups of a stalk with the stalk of the global blow-up. Let $X=\textrm{Spec}  \, A$ with $x$ corresponding to the prime ideal $\mathfrak{p}$ and $\mathcal{I}=\tilde{I}=(x_1, \dots, x_n) \subset A$. Then $\mathcal{I}_x = I_{\mathfrak{p}} \subset A_{\mathfrak{p}}$. Recall the Rees algebras have the relation $A_{\mathfrak{p}}[I_{\mathfrak{p}}t] \cong S^{-1}A[It]$, where $S = A-\mathfrak{p}$ in degree $0$. Let $A[\frac{I}{x_i}]$ be the affine blow-up algebra corresponding to $x_i$. Then for any prime ideal $\mathfrak{q} \in \textrm{Spec}  (A[\frac{I}{x_i}])$ with $\pi(\mathfrak{q})=\mathfrak{p}$ there exists a corresponding prime ideal $\mathfrak{q}' \in \textrm{Spec}  (A_{\mathfrak{p}}[\frac{I_{\mathfrak{p}}}{x_i}])$ containing the unique maximal ideal of $A_{\mathfrak{p}}$ and moreover, $A[\frac{I}{x_i}]_{\mathfrak{q}} \cong A_{\mathfrak{p}}[\frac{I_{\mathfrak{p}}}{x_i}]_{\mathfrak{q}'}$. As blow-ups and stalks can be computed affine locally we have the following lemma.

\begin{lemma}
    \label[lemma]{lemma3102} 
    Let $X$ be a (Noetherian) scheme and $Z \subset X$ a closed subscheme defined by an ideal sheaf $\mathcal{I} \subset \mathcal{O}_X$. Let $\pi:B_Z X \rightarrow X$ be the blow-up of $X$ along $Z$ and consider a point $x' \in B_ZX$ in the fiber over $x$. Write $R = \mathcal{O}_{X,x}$ and $\mathfrak{m} \subset R$ for the unique maximal ideal and $I :=\mathcal{I}_x \subset R$. Then $x'$ corresponds to a point $y$ in $B_{\textrm{Spec}  (R/I)} R$ and $\mathcal{O}_{B_{\textrm{Spec}  (R/I)} R, y} \cong \mathcal{O}_{B_Z X, x'}$ and these can be computed locally as above.
\end{lemma}

\begin{observation}[Applying Hironaka's results]
    \label{observation3103} 
    It follows from the lemma that assuming regularity conditions the stalk at a point $x \in B_Z X$ in the fiber over a point $y \in X$ is a "monoidal transform" of the stalk at $y$ as in \cite[Section 2]{2a777e2f-398f-34da-9fbc-ecc37bb2c7e1} and we can apply the results of that section.
\end{observation}


\subsection{Local blow-up coordinates} 
We will use \cref{lemma3102} and observation \ref{observation3103} in the context of smooth varieties over perfect fields to produce local blow-up coordinates. Combined with some observations on fibers these will give us local blow-up coordinates over $\mathbb{Z}$.

\subsubsection{Local blow-up coordinates for varieties over perfect fields}
We first describe local blow-up coordinates for smooth varieties over perfect fields using results of \cite{cutkoskyresolution} and observation \ref{observation3103}.  

\begin{proposition}
    \label[proposition]{lemma3211} 
    Let $k$ be a perfect field and $X$ a smooth variety over $k$. Let $Z \subset X$ be a smooth closed subvariety and $B_Z X$ be the smooth variety given by the blow-up. Let $x' \in B_Z X$ be a closed point with $x=\pi(x')$ a closed point in $X$. Let $(x_1, \dots, x_n) $ be regular parameters of $\mathcal{O}_{X, x}$ with $Z$ locally defined by the vanishing of $(x_1, \dots, x_r)$. Then there exist regular parameters $(\bar{x}_1, \dots, \bar{x}_n)$ of $\mathcal{O}_{B_Z X, x'}$ such that
    $$
    \bar{x}_1=x_1, \bar{x}_2=\frac{x_2}{x_1}, \dots, \bar{x}_r=\frac{x_r}{x_1}, \bar{x}_{r+1}=x_{r+1}, \dots, \bar{x}_n=x_n,
    $$
    where $x_1$ corresponds to the exceptional divisor and we assume $x'$ lies in that chart.
    
\end{proposition}
\textit{Proof (cf.\cite[3.60.2]{10.2307/j.ctt7rptq}, \cite[Lemma 6.4]{cutkoskyresolution}).} We can argue as in \cite[Lemma A.16]{cutkoskyresolution} by extending to the algebraic closure, producing a system of regular parameters by \cite[Lemma 4]{2a777e2f-398f-34da-9fbc-ecc37bb2c7e1} then showing that these give regular parameters for $\mathcal{O}_{B_Z X, x'}$ using \cite[Corollary A.2]{cutkoskyresolution}. \qedsymbol

\begin{corollary}[Substitution rule]
    \label[corollary]{corollary3212}
    Let $x, x'$ be as in the previous proposition. Then under the local ring morphism $\mathcal{O}_{X, x} \rightarrow \mathcal{O}_{B_Z X, x'}$,
    $$
    x_i \mapsto 
    \begin{cases} 
       \bar{x}_1\bar{x}_i & 2 \leq i \leq r \\
       \bar{x}_i & \textrm{otherwise}
      
    \end{cases}
    $$
    and we can use this to compute local pullbacks of ideals under $\pi$.
\end{corollary}


\subsubsection{Local blow-up coordinates over $\mathbb{Z}$}
We show that blowing up smooth schemes over $\textrm{Spec}  \, \mathbb{Z}$ by smooth centres over $\textrm{Spec}  \, \mathbb{Z}$ commutes with taking fibers. We utilize this fact together with \cref{lemma3211} to produce our desired local blow-up coordinates.

\begin{definition}
    \label[definition]{definition3221} 
    Let $R$ be a Noetherian ring. A sequence of elements $x_1, \dots, x_n \in R$ is a \textit{permutable regular sequence} if any permutation of the sequence is a regular sequence in $R$. For example, if $R$ is local and $x_1, \dots, x_n$ is a regular sequence, then $x_1, \dots, x_n$ is permutable.
\end{definition}

\begin{lemma}
    \label[lemma]{lemma3222} 
     Let $R$ be a Noetherian ring and $x_1, \dots, x_n \in R$ a permutable regular sequence. Let $I=(x_1, \dots, x_r),  r \leq n$ and $f \in R$. Assume $x_if \in I^n$ for $1 \leq i \leq n$ then if $x_i \in I$, $f \in I^{n-1}$ or if $x_i \not \in I$, $f \in I^{n}.$
\end{lemma}

\textit{Proof.} Follows from a theorem of Rees (see \cite[Theorem 1.1.7]{Bruns_Herzog_1998}), which we will make use of frequently. \qedsymbol
\\
\\We apply this lemma to compute strict transforms of certain blow-ups. Let $I=(x_1, \dots, x_r)$ as above. For $1 \leq i \leq n$ define $\bar{R}:=R/(x_i)$ and $\bar{I}:=I+(x_i) \subset \bar{R}$. We can form the commuting squares of blow-up algebras corresponding to $x_j, j \not = i$
$$
\begin{tikzcd}[cramped]
	{\bar{R}[\frac{\bar{I}}{x_j}]} && {R[\frac{I}{x_j}]} \\
	\\
	{\bar{R}} && R
	\arrow[two heads, from=3-3, to=3-1]
	\arrow[from=3-3, to=1-3]
	\arrow[from=3-1, to=1-1]
	\arrow[two heads, from=1-3, to=1-1]
\end{tikzcd}
$$
and $\bar{R}[\frac{\bar{I}}{x_j}] \cong R[\frac{I}{x_j}]/(\frac{x_i}{1})/J$, where $J$ consists of $x_j$-power torsion elements (see \cite[Lemma 10.7.3]{stacks-project}).

\begin{lemma}
    \label[lemma]{lemma3223} 
    $\bar{R}[\frac{\bar{I}}{x_j}] \cong R[\frac{I}{x_j}]/(\frac{x_i}{1})$ if $x_i \not \in I$ and $\bar{R}[\frac{\bar{I}}{x_j}] \cong R[\frac{I}{x_j}]/(\frac{x_i}{x_j})$ if $x_i \in I$.
\end{lemma}
\textit{Proof.} Suppose $x_i \not \in I$ and we have $\frac{a}{x_j^n} \in R[\frac{I}{x_j}]$ with $x_j$-power torsion. Then there exists $s, k, m$ and $b \in R$ such that $x_j^{s+k+m}a=x_j^{s+m}x_ib \in R$ and it follows that $x_j^{s+k+m}a \in (x_i)$. But $x_j$ will be a non-zero divisor in $R/(x_i)$ so $a \in (x_i)$ and hence $a=x_i g \in I^n$. Then by \cref{lemma3222} we have $g \in I^n$ and $\frac{a}{x_j^n} = \frac{x_i}{1}\frac{g}{x_j^n}$, which implies $J=0$. Now suppose $x_i \in I$, then as 
$\frac{x_i}{1}=\frac{x_j}{1}\frac{x_i}{x_j}$ we have $(\frac{x_i}{x_j}) \subset J$. Conversely, by the same arguments as above, for any element of $J$ we have $a=x_i g \in I^n$. But $x_i \in I$ implies $g \in I^{n-1}$ and so $\frac{a}{x_j^n} = \frac{x_i}{x_j}\frac{g}{x_j^{n-1}}$. \qedsymbol

\begin{corollary}
    \label[corollary]{corollary3224} 
    Let $R$ be a regular local ring and $x_1, \dots, x_n$ a regular sequence of parameters in $R$ and $I$ as in the lemma. Then every prime ideal $M'$ of $\bar{R}[\frac{\bar{I}}{x_j}]$ containing the unique maximal ideal of $\bar{R}$ is a prime ideal $M$ of $R[\frac{I}{x_j}]$ containing the unique maximal ideal of $R$ and $R_{M} \rightarrow \bar{R}_{M'}$ is given by taking the quotient by $(\frac{x_i}{1})$ or $\frac{x_i}{x_j}$. Moreover, if $x_i \not \in I$ then any prime ideal $M$ of $R[\frac{I}{x_j}]$ containing the unique maximal ideal of $R$ induces a prime ideal $M'$ of $\bar{R}[\frac{\bar{I}}{x_j}]$ containing the unique maximal ideal of $\bar{R}$.
\end{corollary}

\begin{observation}
    \label{observation3225} 
   Let $X$ be a smooth scheme over $\textrm{Spec}  \, \mathbb{Z}$ and $Z \subset X$ a closed subscheme smooth over $\textrm{Spec}  \, \mathbb{Z}$. Let $x \in X$ be contained in the fiber $X_p$.
    \begin{enumerate}
        \item  If $p \not = 0$ then $X_p$ is an effective Cartier divisor of $X$ and by smoothness over $\mathbb{Z}$ we can lift a regular system of parameters for $x_p \in X_p$ to regular parameters for $x \in X$. If $p = 0$ then $\mathcal{O}_{X,x} \cong \mathcal{O}_{X_p, x_p}$ and a regular system of parameters for $x_p \in X_p$ is a regular system of parameters for $x \in X$.
        \item Suppose $x \in Z \subset X$. As $Z$ is smooth over $\mathbb{Z}$, $\mathcal{O}_{Z_p, x_p} = \mathcal{O}_{X_p,x_p}/I_{x_p}$, where $I_{x_p} = (x_1, \dots, x_r)$ for $x_1, \dots, x_r$ regular parameters in the regular local ring $\mathcal{O}_{X_p, x_p}$. By (1) we can find $x_{r+1}, \dots x_n$ such that a regular system of parameters of $\mathcal{O}_{X,x}$ is 
        $$
        \begin{cases}
            x_0=p, x_1, \dots, x_r, x_{r+1}, \dots, x_n & \textrm{ if } p \not = 0
            \\x_1, \dots, x_r, x_{r+1}, \dots, x_n & \textrm{ if } p = 0
        \end{cases}
        $$
        and $\mathcal{O}_{Z, x}=\mathcal{O}_{X,x}/I_{x}$, where $I_{x} = (x_1, \dots, x_r)$. 
    \end{enumerate}   
\end{observation}

\begin{lemma}
    \label[lemma]{lemma3226} 
   With notation as in the above observation, for all prime ideals $(p) \subseteq \mathbb{Z}$ the fiber of the blow-up is the blow-up of the fiber; $(B_Z X)_p \cong B_{Z_p} X_p$. Further, the following diagram commutes
    $$
    \begin{tikzcd}[cramped]
	   {B_{Z_p} X_p \cong (B_Z X)_p} && {B_Z X} \\
	   \\
	   {X_p} && X
	   \arrow[from=1-1, to=1-3]
	   \arrow[from=3-1, to=3-3]
	   \arrow["{\pi_p}"', from=1-1, to=3-1]
	   \arrow["\pi", from=1-3, to=3-3]
    \end{tikzcd},
    $$
    where the horizontal maps are projections from the fibers and the vertical maps are the blow-ups. 
\end{lemma}
\textit{Proof.} By flatness we know this holds for $p=0$. Now suppose $p \not = 0$. In general we have $B_{Z_p} X_p \rightarrow (B_Z X)_p \rightarrow B_Z X$, where by the blow-up closure lemma $B_{Z_p} X_p$ is the strict transform of the closed subscheme $X_p \rightarrow X$ and $(B_Z X)_p$ is the total transform. Let $x' \in B_Z X$ with $x=\pi(x) \in X$ and assume $x'$ and $x$ are contained in $p$-fibers with corresponding $x_p' \in (B_Z X)_p$ and $x_p \in X_p$. Consider $R:=\mathcal{O}_{X, x} \rightarrow \bar{R}: = \mathcal{O}_{X_p, x_p}$ given by mod $p$ (a regular parameter) and apply \cref{corollary3224}. Now $x' \in B_Z X$ will correspond to a prime ideal $M \subset R[\frac{I}{x_j}]$ (assuming $x'$ is contained in an $x_j$-chart) containing the maximal ideal of $R$ and so induces a prime ideal $M' \subset \bar{R}[\frac{\bar{I}}{x_j}]$ containing the maximal ideal of $\bar{R}$. In particular, this corresponds to a point $y_p' \in B_{Z_p} X_p$ and moreover, the induced morphism $R[\frac{I}{x_j}]_M \rightarrow \bar{R}[\frac{\bar{I}}{x_j}]_{M'}$ is given by mod $p$. It follows that $
\mathcal{O}_{(B_Z X)_p, x_p'} \cong \mathcal{O}_{B_{Z_p} X_p, y_p'} \cong \mathcal{O}_{ B_Z X, x'} /(p)$ and so $y_p'$ maps to $x_p'$ under $B_{Z_p} X_p \rightarrow (B_Z X)_p$.  \qedsymbol

\begin{lemma}[local blow-up coordinates] 
    \label[lemma]{lemma3227} 
    Let $x' \in B_Z X$ be a point such that $x_p' \in (B_Z X)_p \cong B_{Z_p} X_p$ is closed. Let $x=\pi(x')$ with $x_p$ closed in $X_p$. Let a system of regular parameters for $\mathcal{O}_{X,x}$ be
    $$
    \begin{cases}
        x_0=p, x_1, \dots, x_r, x_{r+1}, \dots, x_n & \textrm{ if } p \not = 0
        \\x_1, \dots, x_r, x_{r+1}, \dots, x_n & \textrm{ if } p = 0
    \end{cases}
    $$
    with $Z$ defined by $x_1, \dots, x_r$. Then there exists a system of regular parameters of $\mathcal{O}_{B_Z X, x'}$
    $$
    \begin{cases}
        \bar{x}_0=p, \bar{x}_1=x_1, \bar{x}_2=\frac{x_2}{x_1}, \dots, \bar{x}_r=\frac{x_r}{x_1}, \bar{x}_{r+1}=x_{r+1}, \dots, \bar{x}_n=x_n & \textrm{ if } p \not = 0
        \\\bar{x}_1=x_1, \bar{x}_2=\frac{x_2}{x_1}, \dots, \bar{x}_r=\frac{x_r}{x_1}, \bar{x}_{r+1}=x_{r+1}, \dots, \bar{x}_n=x_n & \textrm{ if } p = 0
    \end{cases}
    ,
    $$
    where $x_1$ corresponds to the exceptional divisor and we assume $x'$ lies in such a chart.
\end{lemma}
\textit{Proof.} As $x_1, \dots, x_n$ are regular parameters on $\mathcal{O}_{X_p, x_p}$ it follows from the local blow-up coordinates for varieties over perfect fields and commutativity of the diagram in \cref{lemma3226} that there exist elements $\bar{x}_1=x_1, \bar{x}_2=\frac{x_2}{x_1}, \dots, \bar{x}_r=\frac{x_r}{x_1}, \bar{x}_{r+1}=x_{r+1}, \dots, \bar{x}_n=x_n \in \mathcal{O}_{B_Z X, x'}$ given by the local blow-up such that mod $p$ they form a regular sequence of parameters in $\mathcal{O}_{(B_Z X)_p, x'_p}$. If $p \not = 0$ then together with $p$ they form regular parameters in $ \mathcal{O}_{B_Z X, x'}$. \qedsymbol

\begin{remark}
    \label[remark]{remark3229} 
    As in \cref{corollary3212} under the local ring morphism $\mathcal{O}_{X, x} \rightarrow \mathcal{O}_{B_Z X, x'}$ we have the same substitution rule with the inclusion of $p \mapsto p$ when $p \not = 0$ and can similarly compute pullbacks under $\pi$ locally. 
\end{remark}


\subsection{Simple normal crossings divisors}
\label{C3_S4} 
We describe $\mathbb{Z}$-flat simple normal crossings divisors for schemes $X$ smooth over $\textrm{Spec}  \, \mathbb{Z}$. When $X$ has a $\Lambda$-structure we will describe additional properties on divisors that are compatible with the $\Lambda$-structure and produce blow-up centres as in \cref{C:2}.

\begin{definition}
    (cf. \cite[Definition 3.24]{10.2307/j.ctt7rptq}).
    Let $X$ be a Noetherian regular scheme. Then $E=\sum E^i$ is a simple normal crossings divisor on $X$ if each $E^i$ defines a divisor that is regular and for each point $x \in X$ there exist regular parameters $x_1, \dots, x_n$ of $\mathcal{O}_{X, x}$ such that for each $i$
    \begin{enumerate}
     \item $x \not \in E^i$ or
     \item $E^i$ is locally in a neighbourhood of $x$ defined by the vanishing of $x_{c(i)}$.
     \item $c(i) \not = c(i')$ if $i \not = i'$.
    \end{enumerate}
    We make clear that each $E^i$ is distinct and so locally $E$ is given by the vanishing of $x_{c(1)} \cdot x_{c(2)} \cdots x_{c(k)}$ for the $k$ components containing $x$. We will also sometimes write $E=\bigcup E^i$ to emphasize that $E^i$ are the components of $E$. This is also known as a strict normal crossings divisor (see \cite[Definition 41.21.1]{stacks-project}). We will often use the acronym "snc" for simple normal crossings.
\end{definition}

A closed subscheme $Z$ of $X$ has \textit{normal crossings} with $E$ if regular parameters can be chosen such that in addition $Z$ is locally defined in a neighbourhood of $x$ by the vanishing of $x_{j_1}, \dots, x_{j_k}$. In particular, $Z \subseteq X$ is a regular embedding. A closed subscheme $Z$ with normal crossings to $E$ will be \textit{transverse} to $E$ if locally at $x$ none of the $x_{j_r}$ defining $Z$ correspond to an $E^i$. It follows that for $Z$ transverse any finite intersection $\bigcap_{j \in J} E^j$ will meet $Z$ transversally in the usual sense.

\begin{lemma}
    \label[lemma]{lemma3302} 
    Let $Z$ have normal crossings with $E=\sum_i E^i$ and suppose $E$ does not contain any component of $Z$ and $E \cap Z \not = \emptyset$. Then $E|_Z$ is a simple normal crossings divisor on $Z$.
    
\end{lemma} 

\textit{Proof.} If $E$ is transverse to $Z$ then clearly $E|_Z = \sum_i E^i|_Z$ is a simple normal crossings divisor on $Z$. Suppose first that $Z$ is irreducible. Each $E^i$ does not contain the generic point of $Z$ and so $E^i$ meets $Z$ transversally and hence $E$ meets $Z$ transversally. Then $E|_Z = (\sum_{i: Z \cap E^i \not = \emptyset}  E^i)|_Z = \sum_{i: Z \cap E^i \not = \emptyset} E^i|_Z$ is a simple normal crossings divisor on $Z$. Now as $Z$ is regular Noetherian it is a disjoint union of its components $Z=\bigsqcup_j V_j$. But each $V_j$ will have normal crossings with $E$ and so $E|_Z=\bigsqcup_j E|_{V_j} = \bigsqcup_j \sum_{i: V \cap E^i \not = \emptyset} E^i|_{V_j} $ is a disjoint union of simple normal crossings divisors on $Z$ and hence is simple normal crossings. \qedsymbol

\begin{example}
    If $E=\sum_j E^j$ is a simple normal crossings divisor and we take a finite intersection $Z=\bigcap_j E^j$ then $Z$ will have normal crossings with $E$.
\end{example}

\textbf{$\mathbb{Z}$-flat simple normal crossings divisors.} We are most interested in the $\mathbb{Z}$-flat setting, where $\Lambda$-structures can be easily considered. Let $X$ be a smooth scheme over $\textrm{Spec}  \, \mathbb{Z}$ and $E=\sum E^i$ a simple normal crossings divisor. 

\begin{observation}
    \label[observation]{observation3305}  
    If $Z$ has normal crossings with $E$ and is flat over $\mathbb{Z}$ then $Z$ is smooth over $\mathbb{Z}$. This follows as locally at any $x$, flatness implies $Z$ is not cut out by $p$. The fibers of $Z$ (given locally by mod $p$) will then have the same local dimension, given by the difference of the local dimension of $X$ and the number of parameters defining $Z$. Further, by the proof of \cref{lemma3302} if $E, Z$ are $\mathbb{Z}$-flat and $E$ contains no components of $Z$ then $(\sum_{i: Z  \cap E^i \not = \emptyset} E^i)|_Z = \sum_{i: Z \cap E^i \not = \emptyset} E^i|_Z$ is a $\mathbb{Z}$-flat simple normal crossings divisor on $Z$. 
\end{observation}

\begin{example}
    \label[example]{example3306} 
    Let $X$ be smooth over $\mathbb{Z}$, $E$ be a $\mathbb{Z}$-flat simple normal crossings divisor and $Z=\cap_{j \in J} E^j$ be a finite intersection. Then $Z$ has normal crossings to $E$ and is $\mathbb{Z}$-flat so is smooth over $\mathbb{Z}$. In particular, $Z$ will be a disjoint union of its components.
\end{example}

\begin{lemma}
    \label[lemma]{lemma3308} 
  Let $X$ be smooth over $\mathbb{Z}$, $E$ a $\mathbb{Z}$-flat simple normal crossings divisor, and $(p)$ a prime ideal of $\mathbb{Z}$. Then the pullback $E_p$ is a simple normal crossings divisor on $X_p$.  
\end{lemma}
\textit{Proof.} If $p = 0$ then it is clear by flatness. Now suppose $p \not = 0$. Let $x_p \in X_p$ correspond to $x \in X$ then $\mathcal{O}_{X_p, x_p} \cong \mathcal{O}_{X,x}/(p)$ and so if $x_0=p, x_1, \dots, x_n$ is a regular system of parameters for $\mathcal{O}_{X,x}$ then $x_1, \dots, x_n$ is a regular system of parameters for $\mathcal{O}_{X_p, x_p}$. It follows that the pullback $(E^i)_p$ of each $E^i$ is a regular embedding of codimension 1 on $X_p$ and so an effective Cartier divisor. The local condition on parameters is clear. \qedsymbol

\begin{definition}[Total transforms]
    \label[definition]{definition3309} 
    Let $\pi:B_Z X \rightarrow X$ be the blow-up of $X$ along $Z$ smooth over $\mathbb{Z}$ with normal crossings to a $\mathbb{Z}$-flat simple normal crossings divisor $E=\sum_{j \in J} E^j$. Note that $Z$ is flat over $\mathbb{Z}$ and so locally can be written in local parameters such that none are $p$. Define the total transform of $E$ as 
    $$
    \pi_{\ast}^{-1}(E)=\sum_{j \in J} \pi_{\ast}^{-1}(E^j)+F,
    $$ 
    where $\pi_{\ast}^{-1}(E^j)$ is the strict transform of $E^j$ and $F$ is the exceptional divisor. If $J$ is totally ordered then $\pi_{\ast}^{-1}(E)=\sum_{j \in J'} E^j$, where $j'$ corresponds to $F$ and $J'=J \cup \{j'\}$ is totally ordered with $j <  j'$ for all $j \in J.$
\end{definition}

\begin{lemma}
    \label[lemma]{lemma33011}  
    The total transform of $E$ is a $\mathbb{Z}$-flat simple normal crossings divisor.
\end{lemma}
\textit{Proof.} Let $x' \in B_Z X$ such that $x_p'$ is closed in $(B_Z X)_p$ and $x=\pi(x')$. Choose a system of regular parameters $x_0=p, x_1, \dots, x_n$ for $\mathcal{O}_{X, x}$ (no $x_0$ if $p = 0$ see \cref{lemma3227}), where $E^i$ containing $x$ corresponds to $x_i$ and $Z$ is defined by the vanishing of $x_1, \dots, x_r$. We use the local blow-up coordinates of $\mathcal{O}_{B_Z X, x'}$ corresponding to a chart given by an $x_j$. Suppose a strict transform ${E^i}'$ contains $x'$ then $E^i$ contains $x$ and so corresponds to a parameter $x_i \not = p$. By \cref{lemma3223}, ${E^i}'$ corresponds to $\bar{x}_i \not = p$ and the exceptional divisor corresponds to the parameter defining the chart. This will hold in a neighbourhood of such points and so it follows that each ${E^i}'$ is a $\mathbb{Z}$-flat divisor on $B_Z X$, which proves the lemma. \qedsymbol

\begin{corollary}
    \label[corollary]{corollary33012} 
    Let $(p) \in \textrm{Spec}  \, \mathbb{Z}$ and write $X_p$ for the fiber, a smooth variety over $\mathbb{F}_p$. Let $Z$ be a centre as above and write $\pi_p:B_{Z_p} X_p \rightarrow X_p$ for the blow-up induced by $Z_p$. Then the total transform of $E_p$ by $\pi_p$ is
    $$
    (\pi_p)_{\ast}^{-1}(E_p):=\sum (\pi_p)_{\ast}^{-1}((E^i)_p)+F_p \cong (\pi_{\ast}^{-1}(E))_p,
    $$
    which is a simple normal crossings divsor on $B_{Z_p} X_p \cong (B_Z X)_p$.
\end{corollary}

\textit{Proof.} Follows using the same blow-up coordinate arguments as in \cref{lemma33011}. \qedsymbol

\textbf{Affine open neighbourhoods.} In order for our divisors to be compatible with affine $\Lambda$-covers we need to consider what happens to divisors in affine neighbourhoods. For this we recall permutable regular sequences (see \cref{definition3221}).

\begin{lemma}[Permutable regular sequence in a neighbourhood]
    \label[lemma]{lemma33013} 
    Let $A$ be a Noetherian ring and $x_1, \dots, x_n$ a sequence of elements in $A$ with the following property: For each $\mathfrak{p} \in \textrm{Spec}  \, A$, writing $x_{i_1}, \dots x_{i_k}, k \leq n$ for the $x_i \in \mathfrak{p}$, $x_{i_1}, \dots, x_{i_k}$ form a permutable regular sequence in $A_{\mathfrak{p}}$ and $(x_{i_1}, \dots, x_{i_k}) \not = A_{\mathfrak{p}}$. Then $x_1, \dots, x_n$ is a permutable regular sequence in $A$.
\end{lemma}

\textit{Proof.} Consider $A/(x_1, \dots, x_{m-1})$ and the endomorphism given by multiplying by $x_m, m \leq n$. Let $J_m$ be the kernel of this map. Then $x_m$ is a non-zero divisor on $A/(x_1, \dots, x_{m-1})$ if and only if $J_m$ is $0$ (see \cite[Exercise 8.4.G]{FOAG}). But for all prime ideals $\mathfrak{p} \in \textrm{Spec}  \, A/(x_1, \dots, x_{m-1})$,
$
(J_m)_{\mathfrak{p}} \subset A_{\mathfrak{p}}/(x_1, \dots, x_{m-1})_{\mathfrak{p}}
$
is the kernel of the multiplication map by $x_m$ in the stalk at $\mathfrak{p}$ and hence is $0$ by the assumption that $x_1, \dots, x_{m-1}$ is a regular sequence. As this is true at all $\mathfrak{p}$, $J_m=0$ and this holds for any $m$ and any permutation of the sequence. \qedsymbol

\begin{corollary}
    \label[corollary]{corollary33014}
    Let $X$ be smooth over $\mathbb{Z}$ and $E=\sum_i E^i$ a $\mathbb{Z}$-flat simple normal crossings divisor. Let $\textrm{Spec}  \, A$ be an affine open where $E$ is defined by a non-zero divisor $x$. In particular, $x=x_1 \cdots x_n$, where each $x_i \in A$ corresponds to the $E^i$ that intersect $\textrm{Spec}  \, A$. Then $x_1, \dots, x_n$ is a permutable regular sequence on $A$.
\end{corollary}

\textit{Proof.} Follows from \cref{lemma33013} as locally $x_1, \dots, x_n$ form part of a system of parameters in regular Noetherian local rings. \qedsymbol

\begin{lemma}[Total transforms in an open neighbourhood]
    \label[lemma]{lemma33015} 
    Let $X$ be a smooth scheme over $\mathbb{Z}$, $E$ a $\mathbb{Z}$-flat simple normal crossings divisor, and $Z$ smooth over $\textrm{Spec}  \, \mathbb{Z}$ with normal crossings to $E$. Let $x, x'$ be as in \cref{lemma33011}. Then there exists an affine open neighbourhood $\textrm{Spec}  \, A$ of $x$ where $Z$ is given by $I=(x_1, \dots, x_r)$ and each $E^i$ intersecting $\textrm{Spec}  \, A$ is defined by the vanishing of a $x_i \in A$ such that the affine blow-up chart $\textrm{Spec}  \, A[I/x_j]$ is an affine open neighbourhood of $x'$ and each ${E^i}'$ is defined on $\textrm{Spec}  \, A[I/x_j]$ by $\bar{x}_i=\frac{x_i}{1}$ or $\frac{x_i}{x_j}$. Moreover, these affine open neighbourhoods commute with the local blow-up description when localizing at $x, x'$.
\end{lemma}
\textit{Proof.} Let $\textrm{Spec}  \, A$ containing $x$ be the intersection of the neighbourhood described in \cref{corollary33014}, where $E$ is defined by a non-zero divisor, and a neighbourhood where $Z$ is defined by the vanishing of $x_1, \dots, x_r$, which become parameters locally at $x$ and at all points in the neighbourhood. It follows from \cref{lemma33013} that on $A$ the sequence given by combining $x_1, \dots, x_r$ with the $x_i's$ corresponding to the $E^i$ that intersect $\textrm{Spec}  \, A$ form a permutable regular sequence. Then $\textrm{Spec}  \, A[I/x_j]$ will be an affine open neighbourhood of $x'$ and by \cref{lemma3223} the ${E^i}'$ will be defined in $\textrm{Spec} \, A[I/x_j]$ by $\bar{x}_i$. By the universal property of the local blow-up, $\bar{x}_i$ in the stalk at $x'$ will precisely be the regular parameter given by the local blow-up. \qedsymbol 

\begin{lemma}
    Let ${E^i}'$ be the strict transforms of the $E^i$ intersecting $\textrm{Spec}  \,A$ as in the previous lemma, defined on $A[I/x_j]$ by $\bar{x}_i$. Then the $\bar{x}_i$'s and $\bar{x}_j$ defining the exceptional divisor form a permutable regular sequence on $A[I/x_j]$.
\end{lemma}
\textit{Proof.} Let $x' \in \textrm{Spec}  (A[I/x_j])$ be a point closed in its fiber, which implies $x=\pi(x) \in \textrm{Spec}  \,A$ is as well. It follows that some of the $x_i \in A$ will be part of a system of regular parameters locally at $x$ and so by the local blow-up the $\bar{x}_i$ on which $x'$ vanishes will locally be part of a system of regular parameters at $x'$. These such points are dense and so it follows that we have local permutable regular sequences at all points in $\textrm{Spec}  (A[I/x_j])$ and we can apply \cref{lemma33013}. \qedsymbol

\label{3.3.1}

\subsubsection{Locally toric simple normal crossings divisors}

We define an additional property that makes the divisors compatible with affine $\Lambda$-covers. These will be used to define the centres of blowing up in the order reduction of \cref{C:5}.

\begin{definition}[Locally toric simple normal crossings divisors]
    \label[definition]{definition3311} 
    Let $X$ be a $\mathbb{Z}$-flat $\Lambda$-scheme with an affine $\Lambda$-cover. Let $E$ be a $\mathbb{Z}$-flat simple normal crossings divisor on $X$. We will say $E$ is \textit{locally toric} for the affine $\Lambda$-cover on $X$ if for each affine open $\textrm{Spec}  \, A$ in the cover the $E^i$ intersecting $\textrm{Spec}  \, A$ are given by $x_i \in A$ that form permutable regular sequences (see \cref{corollary33014}) and $\psi_p(x_i)=(x_i)^p$ for all $i$ and $p$, where $\psi_p$ are the Frobenius lifts of $A$.
\end{definition}

\begin{lemma}
    \label[lemma]{lemma3313} 
    Let $X$ be a $\mathbb{Z}$-flat $\Lambda$-scheme with an affine $\Lambda$-cover. Let $E$ be a $\mathbb{Z}$-flat simple normal crossings divisor on $X$ that is locally toric for the cover and let $Z :=\bigcap_{j \in J} E^j$ be a finite intersection with corresponding ideal sheaf $\mathcal{I}$. Then 
    \begin{enumerate}
        \item $\mathcal{I}$ has $\mathbb{Z}$-flat normal cone and is locally toric for the $\Lambda$-cover of $X$.
        \item $Z$ is a $\mathbb{Z}$-flat $\Lambda$-scheme with an affine $\Lambda$-cover and $Z \hookrightarrow X$ is $\Lambda$-equivariant.
    \end{enumerate}
    
\end{lemma}

\textit{Proof.} By the definition of $E$, $\mathcal{I}$ will be locally toric for the cover. $Z$ will have normal crossings with $E$ by \cref{example3306}. Let $\textrm{Spec}  \, A$ be an affine open in the $\Lambda$-cover of $X$ and $I=(x_1, \dots, x_r)=\mathcal{I}(\textrm{Spec}  \, A)$. By \cref{lemma3222}, $\frac{A}{I^k}$ is $\mathbb{Z}$-flat for all $k$, proving (1) by \cref{lemma2101}. By \cref{example3306}, $Z$ is $\mathbb{Z}$-flat and so $A/I$ is a $\mathbb{Z}$-flat ring. As $I$ is $\Lambda$-equivariant the Frobenius lifts of $A$ descend to $A/I$ so that $A/I$ is a $\Lambda$-ring and $A \rightarrow A/I$ is $\Lambda$-equivariant. The $\textrm{Spec}  (A/I)$ clearly form an affine $\Lambda$-cover of $Z$ and induce a $\Lambda$-equivariant $Z \hookrightarrow X$, which proves (2). \qedsymbol

\begin{corollary} 
    \label[corollary]{corollary3314} $B_Z X$ is a $\Lambda$-scheme with an affine $\Lambda$-cover and $\pi:B_Z X \rightarrow X$ is $\Lambda$-equivariant.
\end{corollary}

\textit{Proof.} Follows from \cref{proposition2208}. \qedsymbol

\begin{lemma}  
    \label[lemma]{lemma3315}  
    Let $X$ be a $\mathbb{Z}$-flat $\Lambda$-scheme with an affine $\Lambda$-cover. Let $E$ be a $\mathbb{Z}$-flat simple normal crossings divisor on $X$ that is locally toric for the cover. Let $Z=\bigcap_{j \in J} E^j$ and $\pi: B_Z X \rightarrow X$ the blow-up. Let $E'$ be the total transform of $E$, which is a $\mathbb{Z}$-flat normal crossings divisor by \cref{lemma33011}. Then $E'$ is locally toric for the affine $\Lambda$-cover of $B_Z X$.
\end{lemma}

\textit{Proof.} Recall the affine $\Lambda$-cover of $B_Z X$ is given by $\textrm{Spec}  (A_i\left[\frac{I_i}{x_j} \right])$, where $\textrm{Spec}  (A_i)$ is part of the cover of $X$. It follows from \cref{lemma33015} that ${E^i}'$ intersecting $\textrm{Spec}  (A_i\left[\frac{I_i}{x_j} \right])$ is given by $\bar{x}_i \in A_i\left[\frac{I_i}{x_j} \right]$, which is a toric element under the Frobenius lifts of $A_i\left[\frac{I_i}{x_j} \right]$. \qedsymbol

\begin{observation}[$\mathbb{Z}$-smoothness]
    \label{observation3316} 
    Let notation be as in \cref{lemma3313} and suppose $X$ is smooth over $\mathbb{Z}$. By \cref{example3306}, $Z$ will be smooth over $\mathbb{Z}$ and so a disjoint union of its irreducible components, which are smooth over $\mathbb{Z}$. From the proof of \cref{lemma3313} each component $Z_i$ is a $\mathbb{Z}$-smooth $\Lambda$-scheme with $Z_i \hookrightarrow X$ $\Lambda$-equivariant. 
\end{observation}

\subsubsection{Intersection distinguished divisors}

\label{332} 
To do order reduction in maximal order in \cref{C:5} we require the existence of distinguished points (cf. \cite[Section 6.2]{BM1}) for affine opens in the $\Lambda$-cover.

\begin{definition}
    \label[definition]{definition3321} 
    Let $X$ be a smooth $\Lambda$-scheme over $\textrm{Spec}  \, \mathbb{Z}$ with an affine $\Lambda$-cover and let $E=\sum_i E^i$ be a $\mathbb{Z}$-flat snc divisor that is locally toric on the cover of $X$. We will say $E$ is \textit{intersection distinguished} for the cover if for any $\textrm{Spec}  \, A$ in the $\Lambda$-cover 
    $$
    \bigcap_{k \in K} E^k \cap \textrm{Spec}  \, A \not = \emptyset, \textrm{ where } k \in K \iff E^k \cap \textrm{Spec}  \, A \not = \emptyset,
    $$
    i.e there exists a distinguished point in $\textrm{Spec}  \, A$ contained in each $E^k$ that intersects $\textrm{Spec}  \, A$. Note for any $x \in \bigcap_{k \in K} E^k \cap \textrm{Spec}  \, A$, $y \in \bar{\{x\}} \subset \textrm{Spec}  \, A$ will have $y \in \bigcap_{k \in K} E^k \cap \textrm{Spec}  \, A$. In particular, there will be a closed distinguished point in $\bigcap_{k \in K} E^k \cap \textrm{Spec}  \, A$ and we may use local blow-up coordinates.
\end{definition}

\begin{proposition}
    \label[proposition]{lemma3323} 
    Let $X$ be a smooth $\Lambda$-scheme over $\textrm{Spec}  \, \mathbb{Z}$ with an affine $\Lambda$-cover and let $E=\sum_i E^i$ be a $\mathbb{Z}$-flat snc divisor that is locally toric on the cover of $X$. Then either $E$ is intersection distinguished for this cover or we can refine the affine $\Lambda$-cover to one on which $E$ is intersection distinguished.
\end{proposition}

\textit{Proof.} Let $\textrm{Spec}  \, A$ be an affine open in the $\Lambda$-cover. We first show that we can cover $\textrm{Spec}  \, A$ by distinguished affine opens on which $E$ is intersection distinguished. As above let $K$ be the index set for the components intersecting $\textrm{Spec}  \, A$. Then either $\textrm{Spec}  \, A \cap \bigcap_{k \in K} E^k \not = \emptyset$ and $\textrm{Spec}  \, A$ is intersection distinguished or we would like to cover $\textrm{Spec}  \, A$ with distinguished opens $\textrm{Spec}  \, A_{K'}$ for $K' \subset K$ with
$$
A_{K'}:= A_{x^{\alpha}}, \quad x^{\alpha} :=\Pi_{k \in K'} x_k,
$$
where $x_k$ corresponds to $E^k$. The divisors on $\textrm{Spec}  \, A_{K'}$ will then be $E^k: k \in K-K'$ and if
$$
\bigcap_{k \in K-K'} E^k \cap \textrm{Spec}  \, A \bigcap_{k \in K'} (E^k)^C \not = \emptyset 
$$
then $\textrm{Spec}  \, A_{K'}$ will be intersection distinguished. Let $x \in \textrm{Spec}  \, A$. If $x \not \in E$ then $x \in \textrm{Spec}  \, A_{K'}$ with $K'=K$ and $ \textrm{Spec}  \, A_{K'}$ intersects no components of $E$. Suppose $x$ is contained in at least one $E^k$, then there exists $K'' \subset K$ such that $x \in \bigcap_{k \in K''} E^k$ and $x \not \in E^k$ if $k \not \in K''$. Then there exists $K'$ with $K''=K-K'$, which implies $x \not \in E^k$ if $k \in K'$ and so $x \in \bigcap_{k \in K-K'} E^k \cap \textrm{Spec}  \, A \bigcap_{k \in K'} (E^k)^C$.
\\Now since each $x_i$ is toric for the $\Lambda$-structure of $A$, each $A_{K'}$ will have a natural induced $\Lambda$-structure such that $A \rightarrow A_{K'}$ is $\Lambda$-equivariant. Moreover, the divisors of $E$ intersecting $\textrm{Spec}  \, A_{K'}$ will be given by toric elements. It remains to show that the $\textrm{Spec}  \, A_{K'}$ meet the intersection conditions for an affine $\Lambda$-cover. The intersection of $\textrm{Spec}  \, A_{K_1'}$ and $\textrm{Spec}  \, A_{K_2'}$ in $\textrm{Spec}  \, A$ will be $\textrm{Spec}  \, A_{K_3'}$ where $K_3' = K_1' \cup K_2'$, which is a distinguished $\Lambda$-equivariant open of both. Let $\textrm{Spec}  \, A, \textrm{Spec}  \, B$ be two affine opens and $\textrm{Spec}  \, A_{K_1'},\textrm{Spec}  \, B_{K_2'}$ two opens in the refined covers of $\textrm{Spec}  \, A, \textrm{Spec}  \, B$. If $\textrm{Spec}  \, A_{K_1'},\textrm{Spec}  \, B_{K_2'}$ intersect then there exists $f \in A, g \in B$ with $A_f \cong B_g$ such that $A \rightarrow A_f, B \rightarrow B_g$ are $\Lambda$-equivariant and we write $J$ to represent the subset of divisors of $E$ that intersect $\textrm{Spec}  (A_f) \cong \textrm{Spec}  (B_g)$, which will be a subset of the divisors intersecting $\textrm{Spec}  \, A$ and those intersecting $\textrm{Spec}  \, B$. Let $J_1', J_2' \subset J$ correspond to the divisors coming from $K_1', K_2'$ respectively. The $\Lambda$-equivariance of the localization morphisms imply the elements given by $J$ are toric in $A_f \cong B_g$. Then $\textrm{Spec}  \, A_f \cap \textrm{Spec}  \, A_{K_1'}=\textrm{Spec}  \, (A_f)_{J_1'}, \textrm{Spec}  \, B_g \cap \textrm{Spec}  \, (B_{K_2'})=\textrm{Spec}  (B_g)_{J_2'}$ and
$
\textrm{Spec}  (B_g)_{J_2'} \cap \textrm{Spec}  (A_f)_{J_1'} = \textrm{Spec}  (A_f)_{J_3'} \cong \textrm{Spec}  (B_g)_{J_3'}, 
$
where $J_3'=J_1' \cup J_2'$ indexes the union of the divisors from $J_1', J_2'$ in $A_f \cong B_g$. Clearly, being localizations by toric elements $(A_f)_{J_1'} \rightarrow (A_f)_{J_3'}, (B_g)_{J_2'} \rightarrow (B_g)_{J_3'}$ are $\Lambda$-equivariant. Letting $x^{\alpha}$ correspond to $K_1'$ we have
\begin{align*}
    &A_{K_1'} \rightarrow (A_{K_1'})_f \cong (A_f)_{K_1'}
    \\& \frac{a}{(x^{\alpha})^n} \mapsto \frac{\frac{a}{(x^{\alpha})^n}}{1} \mapsto \frac{\frac{a}{1}}{(\frac{x^{\alpha}}{1})^n},
\end{align*}
which implies for all $p$
$$
\frac{\psi_p(a)}{(x^{\alpha})^{pn}} \mapsto \frac{\frac{a}{(x^{\alpha})^{pn}}}{1} \mapsto \frac{\frac{\psi_p(a)}{1}}{(\frac{x^{\alpha}}{1})^{pn}} = \frac{\psi_{p, A_f}(\frac{a}{1})}{(\frac{x^{\alpha}}{1})^{pn}}
$$
and so $A_{K_1'} \rightarrow (A_{K_1'})_f \cong (A_f)_{J_1'}$ is $\Lambda$-equivariant and hence $A_{K_1'} \rightarrow (A_{K_1'})_f \cong (A_f)_{J_1'} \rightarrow (A_f)_{J_3'} $ is also. The same holds for $B$ and the $\textrm{Spec}  (A_f)_{J_3'}$ will cover $\textrm{Spec}  \, A_{K_1'} \cap \textrm{Spec}  \, B_{K_2'}$ as $f, g$ vary. Thus the $\textrm{Spec}  \, A_{K'}$ form an affine $\Lambda$-cover for $X$. \qedsymbol

\begin{lemma}
    \label[lemma]{lemma3324} 
    Let $X$ be a smooth $\Lambda$-scheme over $\textrm{Spec}  \, \mathbb{Z}$ with an affine $\Lambda$-cover and a $\mathbb{Z}$-flat snc divisor $E$ that is locally toric and intersection distinguished on the cover. Let $Z = \bigcap_{j \in J} E^j$ be a finite intersection, $X'$ the blow-up by $Z$ and $E'$ the total transform of $E$. Then $E'$ is intersection distinguished for the induced affine $\Lambda$-cover of $X'$.
\end{lemma}

\textit{Proof.} Let $\textrm{Spec}  \, A$ be an affine open in the $\Lambda$-cover and $I=(x_1, \dots, x_r)$ be the ideal corresponding to $Z$ (assume $Z \cap \textrm{Spec}  \, A \not = \emptyset$) and $\textrm{Spec}  \, A[I/x_i]$ an affine open in the $\Lambda$-cover of $X'=B_Z X$. Let $\bigcap_{k \in K} E^k \cap \textrm{Spec}  \, A \not = \emptyset$ so that $J \subseteq K$ and there exists $\mathfrak{p} \in \bigcap_{k \in K} E^k \cap \textrm{Spec}  \, A$ (closed in its fiber) and also $\mathfrak{p} \in Z$. Recall $\bigcap_{k \in K} {E^k}'$ are all the components of $E'$ that intersect $\textrm{Spec}  \, A[I/x_i]$. Let $x_1, \dots, x_r, \dots, x_n \in A$ correspond to the $E^k$ so that $x_1, \dots, x_n \in \mathfrak{p}$. We show there exists a $\mathfrak{q} \in \textrm{Spec}  \, A[I/x_i]$ with $\mathfrak{q} \in \bigcap_{k \in K} {E^k}'$, i.e. $\frac{x_1}{x_i}, \dots, \frac{x_r}{x_i}, x_{r+1}, \dots x_n \in \mathfrak{q}$. Consider the corresponding local blow-up at $A_{\mathfrak{p}}$ given by $I_{\mathfrak{p}}=(x_1, \dots, x_r) \in A_{\mathfrak{p}}$. Then by \cite[Theorem 1.1.7]{Bruns_Herzog_1998} there exists a prime ideal of $\textrm{Spec}  (A_{\mathfrak{p}}[I_{\mathfrak{p}}/x_i])$ containing $\frac{x_1}{x_i}, \dots, \frac{x_r}{x_i}, x_{r+1}, \dots x_n \in A_p[I_{\mathfrak{p}}/x_i]$. It follows from our discussion on local blow-ups (see \ref{C:3_Sec:Local blow-ups}) that this corresponds to a prime ideal $\mathfrak{q} \in \textrm{Spec}  \, A[I/x_i]$ with $\frac{x_1}{x_i}, \dots, \frac{x_r}{x_i}, x_{r+1}, \dots x_n \in \mathfrak{q}$ and $\mathfrak{q}$ lies over $\mathfrak{p}$. \qedsymbol

\section{Order and marked ideals}
\label{C:4}

In this section we describe properties of order on ideal sheaves on smooth schemes over $\textrm{Spec}  \, \mathbb{Z}$. We define marked ideals as in standard resolution of singularities and describe their finite sums and relation to higher order differential operators. In the final sections we prove facts about ideals and embeddings defined by locally toric simple normal crossings divisors, which lay the foundation for our study of $\Lambda$-marked monomial ideals in \cref{C:5}.




\subsection{Order of an ideal sheaf and multiplicity of the exceptional divisor}
\label{C4_1} 
As in standard resolutions of singularities we use order as an invariant to measure the complexity of an ideal sheaf or function at a point in a scheme smooth over $\mathbb{Z}$. We show that a number of facts for order of ideals in smooth varieties (in characteristic $0$ or over perfect fields) hold in our setting of smooth schemes over $\mathbb{Z}$. We prove some results more generally for regular Noetherian schemes but when we compare under blowing up we make use of blow-up coordinates and so assume smoothness over $\mathbb{Z}$. 
\\
\\Let $X$ be a regular Noetherian scheme, e.g. smooth over $\textrm{Spec}  \, \mathbb{Z}$.

\begin{definition}[Order]
    Let $f$ be a section of $\mathcal{O}_X$, the \textit{multiplicty} or \textit{order} of $f$ at a point $x \in X$ in an open set where $f$ is defined is
    $$
    ord_x f: = max\{n: f_x \in \mathfrak{m}_x^n \},
    $$
    where $\mathfrak{m}_x \subset \mathcal{O}_{X, x}$ is the unique maximal ideal of the stalk (a Noetherian regular local ring) and $f_x \in \mathcal{O}_{X, x}$ is the image of $f$. 
\end{definition}

Let $\mathcal{I} \subset \mathcal{O}_X$ be a quasicoherent ideal sheaf.

\begin{definition}
    (cf. \cite[Section 3.4]{10.2307/j.ctt7rptq}).
    Let $x \in X$ and $\mathcal{I}_x \subset \mathcal{O}_{X,x}$ be the ideal in the regular local ring $\mathcal{O}_{X,x}$ with unique maximal ideal $\mathfrak{m}_x$. Then by Krull's intersection theorem we can define the \textit{order of vanishing} or \textit{order} of $\mathcal{I}$ at $x$ to be
    $$
    ord_x \mathcal{I} := max\{n: \mathcal{I}_x \subset \mathfrak{m}_x^n \}.
    $$
\end{definition}

\begin{definition}
    Define the \textit{maximal order} of $\mathcal{I}$ to be
    $$
    \textrm{max-ord }\mathcal{I}: = \textrm{sup}\{ord_x \mathcal{I}:x \in X   \}.
    $$
\end{definition}

\begin{lemma}
    \label[lemma]{lemma4106} 
    If $\textrm{max-ord }\mathcal{I} < \infty$ then the maximal order will occur at a closed point.
\end{lemma}
\textit{Proof.} Let $x_1, \dots, x_t$ be regular parameters for $\mathcal{O}_{X, x}$ then there exists an affine open neighbourhood $\textrm{Spec}  \, R$ of $x$ such that $(x_1, \dots, x_t) = \mathfrak{p} \in \textrm{Spec}  \, R$ corresponds to $x$ (using $X$ Noetherian). Now the closed subscheme defined by $\mathfrak{p}$ is a regular embedding locally at $\mathfrak{p}$ so it follows that after shrinking $\textrm{Spec}  \, R$ further if necessary we can assume that for every $\mathfrak{q} \in V(\mathfrak{p}) \subset \textrm{Spec}  (R)$, $x_1, \dots, x_t$ forms a regular sequence locally at $\mathfrak{q}$. Take such a $\mathfrak{q}$ and consider regular parameters $x_1, \dots, x_t, x_{t+1} \dots,x_r$ in $R_{\mathfrak{q}}$. Then we have
$$
R_{\mathfrak{q}} \xrightarrow[]{} R_{\mathfrak{p}}, \quad x_i \mapsto x_i, i \leq t
$$
given by localizing at $S:=R_{\mathfrak{q}}-\mathfrak{p}R_{\mathfrak{q}}$. Let $I=\mathcal{I}(\textrm{Spec}  \, R)$ and suppose $ord_x \mathcal{I} = ord_{\mathfrak{p}}I \geq m$. Let $g \in I_{\mathfrak{q}}$, then $\frac{g}{1} \in I_{\mathfrak{p}}$ so $\frac{g}{1} \in (\mathfrak{m}_{\mathfrak{p}})^m = (x_1, \dots, x_t)^m \subset R_{\mathfrak{p}}$. Thus there exists $s \in S$ such that $sg \in (x_1, \dots, x_t)^m \subset R_{\mathfrak{q}}$. Then $g \in (x_1, \dots, x_t)^m \subset (\mathfrak{m}_{\mathfrak{q}})^m$ by \cite[Theorem 1.1.7]{Bruns_Herzog_1998} and so $ord_x \mathcal{I} \leq ord_{\mathfrak{q}} \mathcal{I}$. We can take $\mathfrak{q}$ to be a maximal ideal in $R$ and hence a closed point in $X$, finishing the proof. \qedsymbol
\\
\\\textbf{Order at a closed subscheme.} Let $Z \hookrightarrow X$ be an irreducible regular closed subscheme, e.g. $Z, X$ smooth over $\textrm{Spec}  \, \mathbb{Z}$. The \textit{order of }$\mathcal{I}$ along $Z$ is
$$
ord_Z \mathcal{I}:=ord_{\eta} \mathcal{I}
$$
for $\eta$ the generic point of $Z$. If $Z$ is not irreducible, then $ord_Z \mathcal{I} = m$ (respectively, $ord_Z \mathcal{I} \geq m$) will mean the order of $\mathcal{I}$ at every generic point is $m$ (respectively, $\geq m$).

\begin{lemma}
    \label[lemma]{lemma4107} 
    Let $X, Z$ and $\mathcal{I}$ be as above. Let $\mathcal{J}$ be the ideal defining $Z$ and suppose $ord_Z \mathcal{I} \geq m$. Then for any point $x \in Z$, $\mathcal{I}_x \subset  (\mathcal{J}_{x})^m \subset \mathcal{O}_{X,x}.$
\end{lemma}
\textit{Proof.} Let $x \in Z$. As $Z \hookrightarrow X$ is a regular embedding, $\mathcal{J}_x = (x_1, \dots, x_t) \subset \mathcal{O}_{X, x}$ for some regular parameters $x_1, \dots, x_t, x_{t+1}, \dots, x_r$ in $\mathcal{O}_{X, x}$. Then as in the proof of the previous lemma we can choose an affine open neighbourhood $\textrm{Spec}  \, R$ of $x$ such that the generic point $\eta$ of any component of $Z$ containing $x$ is given by $(x_1, \dots, x_t) \subset R$. Then we have the localization map $\mathcal{O}_{X,x} \rightarrow \mathcal{O}_{X, \eta}$ as in the previous lemma and $\mathcal{I}_{\eta} \subset (x_1, \dots, x_t)^m = (\mathfrak{m}_{\eta})^m \subset \mathcal{O}_{X, \eta}$. By the same arguments it follows that for any $g \in \mathcal{I}_x, g \in (x_1, \dots, x_t)^m = (\mathcal{J}_x)^m \subset \mathcal{O}_{X, x}$. \qedsymbol
\\
\\\textbf{Multiplicity of Exceptional divisors.} As noted in \cite[3.4]{10.2307/j.ctt7rptq} the most useful property of order as an invariant is to describe the multiplicity of the exceptional divisor. In the following arguments we need to make use of blow-up coordinates so assume $X$ is smooth over $\textrm{Spec}  \, \mathbb{Z}$ and $Z$ is a closed subscheme smooth over $\textrm{Spec}  \, \mathbb{Z}$ with $\mathcal{J} \subset \mathcal{O}_X$ the corresponding ideal sheaf.  Let  $\mathcal{I} \subset \mathcal{O}_X$ be an arbitrary ideal.

\begin{proposition}
    \label[proposition]{prop4108}
   Suppose $ord_Z \mathcal{I} \geq n > 0$ then
    $$
    \pi^{\ast}\mathcal{I} \subset \mathcal{O}(-E_Z X)^n,
    $$
    where $\pi:B_Z X \rightarrow X$ is the blow-up of $X$ along $Z$ and $\mathcal{O}(-E_Z X)=\pi^{\ast}\mathcal{J} \subset \mathcal{O}_{B_Z X}$ is the ideal sheaf corresponding to the exceptional divisor $E_Z X$ of $B_Z X$. 
\end{proposition}
\textit{Proof.} Let $x$ be a point in $X$ closed in a fiber over $(p)$ and $x' \in B_Z X$ a closed point lying over $x$ and assume that $x \in Z$. Then we may compute $\pi^{\ast}\mathcal{I}$ locally at $x'$ using blow-up coordinates and the substitution rule of \cref{remark3229}. But by the above lemma any $g \in \mathcal{I}_x$ is contained in $(x_1, \dots, x_t)^n=(\mathcal{J}_x)^n$, where $x_0=p, x_1, \dots, x_t, x_{t+1}, \dots, x_r$ (no $x_0$ if $p = 0$ see \cref{lemma3227}) are regular parameters in $\mathcal{O}_{X, x}$ and by $\mathbb{Z}$-flatness of $Z$, $x_i \not = p, 1 \leq i \leq t$. Suppose $x'$ is contained in the chart corresponding to $x_i$. Then by the substitution rule the local pullback of $g$, which generates $(\pi^{\ast}\mathcal{I})_{x'}$, is contained in $(\bar{x}_i)^n$ and the proposition follows. \qedsymbol

\begin{corollary}[Fibers]
    \label[corollary]{corollary4109} 
    Let $Z,X, \mathcal{I}$ be as above and let $(p) \in \textrm{Spec}  \, \mathbb{Z}$ be a prime ideal. Write $X_p, Z_p, \mathcal{I}_p$ for the fibers and pullbacks over $(p)$ and $\pi_p: (B_Z X)_p \cong B_{Z_p} X_p \rightarrow X_p$. Then
    $$
    (\pi_p)^{\ast}\mathcal{I}_p \subset \mathcal{O}(-E_{Z_p} X_p)^n.
    $$
\end{corollary}

\textit{Proof.} Follows by the proof of the proposition and applying the arguments of \cref{lemma3227} and \cref{lemma3226}. \qedsymbol


\subsection{Marked ideals and birational transforms}
\label{C4:Markedidealsbirationaltransforms} 
To prove order reduction in \cref{C:5} we will make use of an induction on dimension and will need to restrict to smooth closed subschemes of smaller dimension. However, in general taking order does not commute with these restrictions (see \cite[Section 3.58, 3.9]{10.2307/j.ctt1b7x7vc} for this issue in characteristic $0$). To deal with this we use marked ideals and their birational transforms. For notational convenience we define the birational transform of a marked ideal as an ideal without marking. A marking for the transform will be included when we consider $\Lambda$-marked monomial ideals in \cref{C:5}. 

\begin{definition}[Marking data]
    Let $X$ be a scheme. Let $m$ be a non-negative integer, $f$ a section of $\mathcal{O}_X$, and $\mathcal{I} \subset \mathcal{O}_X$ an ideal. A \textit{marked} ideal or section is a pair of the form $(f, m), (\mathcal{I}, m)$. Define the support of a marked section or ideal as 
    $$
    supp(f, m):=\{x \in X: ord_x f \geq m     \}, \quad supp(\mathcal{I}, m):=\{x \in X: ord_x \mathcal{I}\geq m     \}.
    $$
\end{definition}

\begin{definition}
    Let $X$ be a smooth scheme over $\mathbb{Z}$ and $Z \subset X$ a smooth closed subscheme. Let $\pi:B_Z X \rightarrow X$ be the blow-up and $E$ the exceptional divisor. Assume $ord_Z \mathcal{I} \geq m$ and define the \textit{birational transform} of $(\mathcal{I}, m)$ by $\pi$ as the unmarked ideal
    $$
    \pi_{\ast}^{-1}(\mathcal{I}, m) : = \mathcal{O}(mE) \cdot \pi^{\ast}\mathcal{I} \subset \mathcal{O}_{B_Z X},
    $$
    which is an ideal sheaf by \cref{prop4108}. 
\end{definition}

\begin{definition}[Fibers]
    \label[definition]{definition4203}  
     Let notation be as in the previous definition. By \cref{corollary4109}, for any prime ideal $(p) \subset \mathbb{Z}$, writing $\pi_p: B_{Z_p} X_p \rightarrow X_p$ for the blow-up of $X_p$ along $Z_p$ and $E_p$ the exceptional divisor we define the \textit{birational transform} of $(\mathcal{I}_p, m)$ by $\pi_p$ as the unmarked ideal
     $$
    (\pi_p)_{\ast}^{-1}(\mathcal{I}_p, m) : = \mathcal{O}(mE_p) \cdot (\pi_p)^{\ast}\mathcal{I}_p \subset \mathcal{O}_{B_{Z_p} X_p}.
    $$
\end{definition}

\begin{definition}[Permissibility for marked ideals]
    \label[definition]{definition4204} 
    Let $X$ be a scheme smooth over $\textrm{Spec}  \, \mathbb{Z}$ and $(\mathcal{I}, m)$ a marked ideal. A closed subscheme $Z \subset X$ will be a \textit{permissible centre} for the pair $(\mathcal{I}, m)$ if $Z$ is smooth over $\textrm{Spec}  \, \mathbb{Z}$ and $Z \subset supp(\mathcal{I}, m)$ so that $ord_Z \mathcal{I} \geq m$ and we can define the birational transform of $(\mathcal{I}, m)$ as an ideal.
\end{definition}

\textbf{Blow-up coordinates and birational transforms.} Let $X$ be smooth over $\textrm{Spec}  \, \mathbb{Z}$, $\mathcal{I}$ an ideal of $X$ and $Z$ a permissible centre for $(\mathcal{I}, m)$ on $X$. Consider a closed point $x' \in B_Z X$ lying over a closed point $x \in X$. Assume $x \in Z$ and $x',x$ are contained in the fiber over $(p)$. Let $p, x_1, \dots, x_n$ be regular parameters of $\mathcal{O}_{X,x}$ with $Z$ defined by $x_1, \dots, x_r$ and $p, \bar{x}_1, \dots, \bar{x}_n$ be blow-up coordinates corresponding to the $x_1$ chart (no $p$ if $p = 0$ see \cref{lemma3227}).

\begin{definition}
    \label[definition]{definition4205} 
    (See \cite[3.60.3]{10.2307/j.ctt7rptq}).
    Let $f \in \mathcal{I}_{x}$ and write $\pi^{-1}(f)$ for the image of $f$ under the morphism of stalks $\mathcal{O}_{X,x} \rightarrow \mathcal{O}_{B_Z X, x'}$. Then as in \cref{lemma4107}, $f \in (x_1, \dots, x_r)^m$ and $\pi^{-1}(f) \in (\bar{x}_1)^m$. Write $\pi^{-1}(f)=\bar{x}_1^m\bar{f}$ and for the marked section $(f, m)$ define its birational transform by $\pi$ as the unmarked section
    $$
    \pi_{\ast}^{-1}(f, m): = \bar{x}_1^{-m}\pi^{-1}(f)=\bar{f},
    $$
    which is not unique but will generate $(\pi_{\ast}^{-1}(\mathcal{I}, m))_{x'}.$
\end{definition}

\begin{lemma}
    \label[lemma]{lemma4206} 
    (cf. \cite[Lemma 3.61]{10.2307/j.ctt7rptq}).
    Let $\mathcal{I}$ be an ideal sheaf on $X$ and $Z \subset X$ a smooth closed subscheme with $ord_Z I = \textrm{max-ord } \mathcal{I} = m$. Then  
    $$
    \textrm{max-ord }{\pi}^{-1}_{\ast}(\mathcal{I}, m) \leq m.
    $$
\end{lemma}
 
\textit{Proof.} Recall from \cref{lemma4106} the maximal order occurs at a closed point $x$. If $x \not \in Z$ then the stalk at $x'$ is isomorphic and the order stays the same. Thus it suffices to show that for a closed point $x' \in B_Z X$ lying over a closed point $x \in Z$ the order cannot increase. We use the blow-up coordinates and assume $x, x'$ are in the fiber over $(p) \not = 0$. Let $f \in \mathcal{I}_{x}$ have order precisely $m$. As in \cref{lemma4107} $f \in (x_1, \dots, x_r)^m$ and so we can assume (after potentially multiplying by a unit) that $f$ contains a monomial of degree $m$ in the parameters $x_1, \dots, x_r$ corresponding to $\mathcal{J}_x$. If no such monomial existed then the order of $f$ would exceed $m$. By the substitution rule $\pi^{-1}(f)$ contains a monomial of degree $\leq 2m$ in the blow-up coordinates, which implies $\pi_{\ast}^{-1}(f, m)$ contains a monomial of degree $\leq 2m-m=m$ and so
$
ord_{x'} {\pi}^{-1}_{\ast}(\mathcal{I}, m) \leq m. \quad \qed
$

\begin{lemma}[Birational transforms commute with fibers]
    \label[lemma]{lemma4207} 
    For all prime ideals $(p) \subset \mathbb{Z}$, and with $ (\pi_p)_{\ast}^{-1}(\mathcal{I}_p, m)$ as in \cref{definition4203},
    $$
    (\pi_p)_{\ast}^{-1}(\mathcal{I}_p, m) = (\pi^{-1}_{\ast}(\mathcal{I}, m))|_{B_{Z_p} X_p = (B_Z X)_p}.
    $$
\end{lemma}

\textit{Proof.} It suffices to show this locally at points $x' \in B_Z X$ closed in $(B_Z X)_p$ lying over $x \in X_p$. By the proof of \cref{lemma3226} we have for any $f \in \mathcal{I}_x \subset \mathcal{O}_{X,x}$ $(\pi_p)^{-1}([f]) = [\pi^{-1}(f)] \in \mathcal{O}_{B_{Z_p} X_p, x_p}$, where $[f] \in \mathcal{O}_{X_p, x_p}$ represents the class mod $p$. Noting that multiplying by $\bar{x}_1^{-m}$ commutes with mod $p$ finishes the proof. \qedsymbol

\subsubsection{Sums of marked ideals}
\label{431}
 We describe how to take the sum of finitely many marked ideals.

\begin{lemma}
    Let $(R, \mathfrak{m})$ be a regular local ring and $x_1, \dots, x_n$ be regular parameters for $\mathfrak{m}$. Let $I_i, m_i, i = 1, \dots, k$ be a collection of ideals of $R$ and $m_i$ non-negative integers. For each $i$ define $m_{\hat{i}} = \Pi_{j \not = i} m_j$ and $m = \Pi_i m_i$ so that $m_{\hat{i}}m_i = m$. Then
    $$
    \sum_{i=1}^k ((I_i)^{m_{\hat{i}}}) \subset \mathfrak{m}^m \iff I_i \subset \mathfrak{m}^{m_i}\textrm{ for all } i.
    $$
    
\end{lemma}

\textit{Proof.} One direction is clear. The other follows from \cite[Theorem 1.1.7]{Bruns_Herzog_1998}. \qedsymbol

\begin{corollary}
    \label[corollary]{corollary4213} 
    (See \cite[page 79]{10.2307/j.ctt7rptq} for the sum of two ideals).
    Let $\mathcal{I}_i,  i = 1, \dots, k$ be a collection of quasicoherent ideals on a regular Noetherian scheme $X$ and $m_i$ non-negative integers. As before let $m_{\hat{i}} = \Pi_{j \not = i} m_j$ and $m = \Pi_i m_i$. Let $x \in X$ then
    $$
    \textrm{ord}_x(\sum_{i=1}^k(\mathcal{I}_i)^{m_{\hat{i}}}) \geq m \iff \textrm{ord}_x{\mathcal{I}}_i \geq m_i \textrm{ for all } i.
    $$
    It follows that 
    $$
    supp(\sum_{i=1}^k(\mathcal{I}_i)^{m_{\hat{i}}}, m) = \bigcap_{i=1}^k supp(\mathcal{I}_i, m_i) \subset X.
    $$
\end{corollary}

\begin{lemma}[Sums commute with transforms]
    \label[lemma]{lemma4214} 
    Let $X$ be a smooth scheme over $\mathbb{Z}$, $\mathcal{I}_i,  i = 1, \dots, r$ a collection of quasicoherent sheaves on $X$, and $m_i$ non-negative integers. As before let $m_{\hat{i}} = \Pi_{j \not = i} m_j$ and $m = \Pi_i m_i$. Let $Z$ be a blow-up centre permissible for $(\sum_{i=1}^r \mathcal{I}_i^{m_{\hat{i}}}, m)$ and hence each $(\mathcal{I}_i, m_i)$ and let $\pi: X'=B_Z X \rightarrow X$ be the blow-up morphism. Write $(\sum_{i=1}^r \mathcal{I}_i^{m_{\hat{i}}})'$ for $\pi_{\ast}^{-1}(\sum_{i=1}^r \mathcal{I}_i^{m_{\hat{i}}}, m)$ and $(\mathcal{I}_i)'$ for $\pi_{\ast}^{-1}( \mathcal{I}_i, m_i)$. Then
    $$
    (\sum_{i=1}^r \mathcal{I}_i^{m_{\hat{i}}})' = \sum_{i=1}^r ((\mathcal{I}_i')^{m_{\hat{i}}}) \subset \mathcal{O}_{X'}.
    $$
    
\end{lemma}

\textit{Proof.} It suffices to show the statement holds in the local blow-up at a closed point $x' \in X'$ lying over $x \in X$ in the fiber over $(p) \not = 0$. Let $x_0=p, x_1, \dots, x_s$ be a system of regular parameters for $\mathcal{O}_{X,x}$. Let $K_i$ be an index set for a set of generators $f_k \in \mathcal{O}_{X, x}$ of $(\mathcal{I}_i)_x$ so
$$
(f_k, k \in K_i) = (\mathcal{I}_i)_x.
$$
Consider the local blow-up corresponding to a $x_t$-chart. Then  
$$
(\mathcal{I}_i')_{x'} = (\bar{f}_k, k \in K_i),
$$
where $\bar{f}_k$ is as in \cref{definition4205}. Now
$$
(\mathcal{I}_i^{m_{\hat{i}}})_x = (\Pi_{k \in K_i} (f_k)^{l_q }: l \in \mathbb{N}^{|K_i|} \textrm{ with } \sum_q l_q = m_{\hat{i}})
$$
so
$$
(\sum_{i=1}^r \mathcal{I}_i^{m_{\hat{i}}})_x = ((\Pi_{k \in K_i} (f_k)^{l_q }: l \in \mathbb{N}^{|K_i|} \textrm{ with } \sum_q l_q = m_{\hat{i}}), i = 1, \dots, r).
$$
Again, considering a $x_t$-chart
\begin{align*}
    ((\sum_{i=1}^r \mathcal{I}_i^{m_{\hat{i}}})')_{x'} &= ({\bar{x}_t}^{-m}\Pi_{k \in K_i} ((\bar{x}_t)^{m_i}\bar{f}_k)^{l_q}:  l \in \mathbb{N}^{|K_i|} \textrm{ with } \sum_q l_q =  m_{\hat{i}}, i=1, \dots, r)\\
    &= ({\bar{x}_t}^{-m}(\bar{x}_t)^{m_{\hat{i}}m_i}\Pi_{k \in K_i} (\bar{f}_k)^{l_q}:  l \in \mathbb{N}^{|K_i|} \textrm{ with } \sum_q l_q =  m_{\hat{i}}, i=1, \dots, r)\\
    & =( \Pi_{k \in K_i} (\bar{f}_k)^{l_q}: l \in \mathbb{N}^{|K_i|} \textrm{ with } \sum_q l_q =  m_{\hat{i}}, i =1, \dots, r)\\
    & = (\sum_{i=1}^r ((\mathcal{I}_i')^{m_{\hat{i}}}))_{x'}. \quad \qed
\end{align*}
 
We will use the above results in \cref{C:5} when we consider the sum of $\Lambda$-marked monomial ideals and for convenience make the following definition.

\begin{definition}
    \label[definition]{definition4215} 
    Let $\mathcal{I}_i,  i = 1, \dots, k$ be a collection of quasicoherent ideals on a regular Noetherian scheme $X$ and $m_i$ non-negative integers. As before let $m_{\hat{i}} = \Pi_{j \not = i} m_j$ and $m = \Pi_i m_i$. Define the sum of the marked ideals 
    $$
    \sum_{i=1}^k (\mathcal{I}_i, m_i) := (\sum_{i=1}^k(\mathcal{I}_i)^{m_{\hat{i}}}, m).
    $$
\end{definition}


\subsection{Differential operators}
\label{C4_3} For the purposes of order reduction we use higher order differentials for two reasons. First, to show upper-semicontuinity of the order function for marked ideals introduced in \cref{C4:Markedidealsbirationaltransforms}, which follows from higher order derivatives commuting with pullbacks to $(p)$-fibers. Second, classical order reduction in characteristic $0$ uses derivatives of ideals and we compare this with our setting over $\mathbb{Z}$ in the final remarks of \cref{C:5_Sec:Analysis:FM}. Much of the theory of higher order differentials can be found in \cite[Chapter 16]{PMIHES_1967__32__5_0} but we learnt mostly from the excellent honours thesis Differential Operators and Algorithmic Weighted Resolution by JongHyun Lee \cite{DifferentialOperatorsAlgorithmicWeightedResolution}, which introduces the theory and proves facts for smooth varieties over perfect fields. We make use of Hasse derivatives \cite[Section 5.2]{DifferentialOperatorsAlgorithmicWeightedResolution}and that if $X \xrightarrow{f} S$ is of finite presentation, e.g. smooth, then the global constructions restricted to affine opens or stalks can be computed by considering the analogous constructions on rings. Most of our notation is borrowed from this honours thesis. In particular, $\Omega_{A/B}$ will denote the module of K{\"a}hler differentials and ${P^n}_{A/B}$ will denote the module of principal parts of order $n$.

\begin{theorem}[Upper-semicontinuity of order for varieties over perfect fields]
    (See \cite[Theorem 13]{DifferentialOperatorsAlgorithmicWeightedResolution}, \cite[Thereom A.19]{cutkoskyresolution}).
    \label[theorem]{theorem4325} 
    Let $X$ be a smooth variety over a perfect field $k$, $\mathcal{I}$ a coherent ideal and $m$ a non-negative integer. Then
    $$
    supp(\mathcal{I}, m) = V({\mathcal{D}^{\leq m}}_{X/S} \mathcal{I}),
    $$
    hence $supp(\mathcal{I}, m)$ is a closed subset of $X$.
\end{theorem}


\textbf{Hasse derivatives of smooth schemes over $\mathbb{Z}$.} We prove an upper-semicontinuity theorem for smooth schemes over $\mathbb{Z}$ and certain ideals. We do this by showing how Hasse derivatives and regular sequences behave under quotients $mod$ $p$ and utilizing \cref{theorem4325} applied to the fibers $X_p$. 

\begin{proposition}
     Let $X$ be a scheme smooth over $\textrm{Spec}  \, \mathbb{Z}$ and $\mathcal{I}$ a coherent ideal. Let $x \in X$ be a closed point in the fiber over $(p) \not = 0$. Write $R:= \mathcal{O}_{X, x}$ so that $\bar{R} := R/(p) = \mathcal{O}_{X_p, x_p}$. Define
    $$
    \Delta: = ker(R \otimes_{\mathbb{Z}} R \xrightarrow{\textrm{multiplication}} R), \quad \Delta_p:= ker(\bar{R} \otimes_{\mathbb{F}_p} \bar{R} \xrightarrow{\textrm{multiplication}} \bar{R})
    $$
    so that
    $$
    \Delta/\Delta^2 \cong \Omega_{R/\mathbb{Z}} \quad \Delta_p/{\Delta^2}_p \cong \Omega_{\bar{R}/\mathbb{F}_p}. 
    $$
    Let $p, x_1, \dots, x_n$ be a regular system of parameters for $R$ so that $[x_1], \dots, [x_n]$ form a regular system of parameters for $\bar{R}$, where $[x_i]$ is the class of $x_i$ in $\bar{R}$. Then $dx_1, \dots, dx_n \in \Delta$ is a quasi-regular sequence for $\Delta$, $d_p[x_1], \dots, d_p[x_n] \in \Delta_p$ is a quasi-regular sequence for $\Delta_p$, where $dx_i = 1 \otimes x_i-x_i \otimes_1 \in R \otimes_{\mathbb{Z}} R$ and $d_p[x_i] = 1 \otimes [x_i] - [x_i] \otimes 1 \in \bar{R} \otimes_{\mathbb{F}_p} \bar{R}$. Moreover, $d_p[x_i] = [dx_i]$.
\end{proposition}

\textit{Proof.} First note that the last statement makes sense and follows as the following diagram commutes:
$$
\begin{tikzcd}[cramped]
	R && {R \otimes_{\mathbb{Z}} R} \\
	\\
	{\bar{R}} && {\bar{R} \otimes_{\mathbb{F}_p} \bar{R} \cong R \otimes_{\mathbb{Z}}R \otimes_{\mathbb{Z}} \mathbb{Z}/(p)}
	\arrow["d", from=1-1, to=1-3]
	\arrow["{\textrm{mod p}}"', two heads, from=1-1, to=3-1]
	\arrow["{\textrm{mod p}}"', two heads, from=1-3, to=3-3]
	\arrow["{d_p}", from=3-1, to=3-3]
\end{tikzcd}.
$$
Now $[x_1], \dots, [x_n]$ being regular parameters for $\bar{R}$ implies by \cite[Corollary 8]{DifferentialOperatorsAlgorithmicWeightedResolution} that the $d_p[x_i]$ form a quasi-regular sequence for $\Delta_p$ as $\mathbb{F}_p$ is perfect. We will use this to deduce the statement for the $dx_i$. Recall that K{\"a}hler differentials behave well with respect to base changes and we have the commuting diagram
$$
\begin{tikzcd}[cramped]
	R && {\Delta \subset{R \otimes_{\mathbb{Z}} R}} && {\Omega_{R /\mathbb{Z}}} \\
	\\
	{\bar{R}} && {\Delta_p \subset \bar{R} \otimes_{\mathbb{F}_p} \bar{R}} && {\Omega_{\bar{R}/\mathbb{F}_p} \cong \Omega_{R /\mathbb{Z}}\otimes_{\mathbb{Z}} \mathbb{Z}/(p)}
	\arrow["d", from=1-1, to=1-3]
	\arrow["{\textrm{mod p}}"', two heads, from=1-1, to=3-1]
	\arrow["{\textrm {mod }\Delta^2}", from=1-3, to=1-5]
	\arrow["{\textrm{mod p}}"', two heads, from=1-3, to=3-3]
	\arrow["{\textrm{mod p}}"', two heads, from=1-5, to=3-5]
	\arrow["{d_p}", from=3-1, to=3-3]
	\arrow["{\textrm{mod }{\Delta^2}_p}", from=3-3, to=3-5]
\end{tikzcd}
$$
where $[\Delta] \subset \Delta_p$ and the top and bottom maps in the outer rectangle give the respective universal derivations, which we also denote $d, d_p$. Let $\mathfrak{m}=(p, x_1, \dots, x_n) \subset R$ be the unique maximal ideal and recall $\Omega_{R /\mathbb{Z}}$ is a finitely generated $R$-module of rank $n$ by smoothness over $\textrm{Spec}  \, \mathbb{Z}$. Now the $dx_i \in \Omega_{R /\mathbb{Z}}$ have $[dx_i] = d_p[x_i]  \in \Omega_{\bar{R}/\mathbb{F}_p} = \Omega_{R /\mathbb{Z}} / (p)$ generate $\Omega_{\bar{R}/\mathbb{F}_p}$ as a $\bar{R}$-module by \cite[Proposition 11]{DifferentialOperatorsAlgorithmicWeightedResolution}. It follows that $[dx_i] \in \Omega_{R /\mathbb{Z}}/(\mathfrak{m}) = \Omega_{\bar{R}/\mathbb{F}_p} / (x_1, \dots, x_n)$ will generate the module and so by Nakayamas lemma (for local rings) the $dx_i$ will generate $\Omega_{R /\mathbb{Z}}$. Moreover, as $\Omega_{R /\mathbb{Z}}$ is free of rank $n$ the $dx_i$ form a basis. Thus we can similarly apply \cite[Lemma 10]{DifferentialOperatorsAlgorithmicWeightedResolution} and \cite[Lemma 5]{DifferentialOperatorsAlgorithmicWeightedResolution}, which apply generally to smooth morphisms. \qedsymbol

\begin{corollary}
    $$
    {P^m}_{R/\mathbb{Z}} = \bigoplus_{|s| \leq m} R \cdot dx^{(s)}, \quad {P^m}_{\bar{R}/\mathbb{F}_p} = \bigoplus_{|s| \leq m} \bar{R} \cdot d_p[x^{(s)}], 
    $$
    where $d x_i = 1 \otimes x_i - x_i \otimes 1 \in {P^m}_{R/\mathbb{Z}} \textrm{ and } d_p[x_i] = [x_i] \otimes 1 - 1 \otimes [x_i] \in {P^m}_{\bar{R}/\mathbb{F}_p}$. Moreover, $d_p[x_i] = [dx_i]$, where $[dx_i]$ is $dx_i$ mod $p$.
\end{corollary}

\textit{Proof.} ${P^m}_{R/\mathbb{Z}},{P^m}_{\bar{R}/\mathbb{F}_p} $ being free modules with the given bases follows from \cite[Lemma 6]{DifferentialOperatorsAlgorithmicWeightedResolution}. The second statement holds because the following diagram commutes
$$
\begin{tikzcd}[cramped]
	{{R \otimes_{\mathbb{Z}} R}} && {{P^m}_{R/\mathbb{Z}}} \\
	\\
	{ \bar{R} \otimes_{\mathbb{F}_p} \bar{R}} && {{P^m}_{\bar{R}/\mathbb{F}_p}}
	\arrow["{\textrm {mod }\Delta^{m+1}}", from=1-1, to=1-3]
	\arrow["{\textrm{mod p}}"', two heads, from=1-1, to=3-1]
	\arrow["{\textrm{mod p}}"', from=1-3, to=3-3]
	\arrow["{\textrm{mod }{\Delta^{m+1}}_p}", from=3-1, to=3-3]
\end{tikzcd},
$$
where $[\Delta] \subset \Delta_p$ implies $[\Delta^m] \subset {\Delta^m}_p$ so the right vertical arrow is well-defined. \qedsymbol

\begin{corollary}
    \label[corollary]{corollary4328} 
    $$
    {{D}^{\leq m}}_{R/\mathbb{Z}} = \bigoplus_{|s| \leq m} R \cdot D^{(s)}, \quad {{D}^{\leq m}}_{\bar{R}/\mathbb{F}_p} = \bigoplus_{|s| \leq m} R \cdot {D^{(s)}}_p, 
    $$
    where $\{ D^s\}, \{ {D^s}_p\}$ are the respective Hasse derivatives. Moreover, for each $s \in \mathbb{N}^n$, ${D^s}_p=[D^s]$, where $[D^s]$ is the class of $D^s$ mod $p$.
\end{corollary}

\textit{Proof.} ${{D}^{\leq m}}_{R/\mathbb{Z}},{{D}^{\leq m}}_{\bar{R}/\mathbb{F}_p} $ being free modules with the bases given by Hasse derivatives follows from the previous corollary. The last statement follows from the fact that we can extend differential operators along the mod $p$ morphism (see \cite[Lemma 10.133.1]{stacks-project}), the previous corollary and the definition of Hasse derivatives. \qedsymbol

\begin{theorem}
    \label[theorem]{theorem4329} 
    Let $X$ be a scheme smooth over $\textrm{Spec}  \, \mathbb{Z}$ and $\mathcal{I}$ a coherent ideal. Then for each prime ideal $(p) \subset \mathbb{Z}$
    $$
    ({\mathcal{D}^{\leq m}}_{X/\textrm{Spec}  \, \mathbb{Z}} \mathcal{I})_p \cong {\mathcal{D}^{\leq m}}_{X_p/\textrm{Spec} (\mathbb{F}_p)} \mathcal{I}_p,
    $$
    where $\pi_p: X_p \rightarrow X$ is the projection from the fiber over $(p)$, $({\mathcal{D}^{\leq m}}_{X/\textrm{Spec} \, \mathbb{Z}} \mathcal{I})_p =  (\pi_p)^{\ast}(\mathcal{I})$ and $\mathcal{I}_p =(\pi_p)^{\ast} \mathcal{I}$ (both pullbacks as ideals).
\end{theorem}

\textit{Proof.} It suffices to show this holds locally at all closed points of $X_p$. Let $x \in X$ be closed in $X_p$ and $x_p \in X_p$ the corresponding closed point. Let $R:= \mathcal{O}_{X,x}, \bar{R}:=\mathcal{O}_{X_p, x_p}=R/(p)$ and so $\pi_p: R \rightarrow \bar{R}$ is given by mod $p$. If $p = 0$ then $R = \bar{R}$ and the theorem is clear. Assume $p \not = 0$. Let $I:=\mathcal{I}_x $ so $(\mathcal{I}_p)_{x_p}=[I] \subset \bar{R}$ then by the local construction on rings we have 
$$
(({\mathcal{D}^{\leq m}}_{X/\textrm{Spec} \, \mathbb{Z}} \mathcal{I})_p)_{x_p} \cong [{{D}^{\leq m}}_{R/\mathbb{Z}}I], \quad ({\mathcal{D}^{\leq m}}_{X_p/\textrm{Spec} (\mathbb{F}_p)} \mathcal{I}_p)_{x_p} \cong {{D}^{\leq m}}_{\bar{R}/\mathbb{F}_p} [I].
$$
But
$$
\begin{tikzcd}[cramped]
	{{D^{\leq n}}_{R/\mathbb{Z}}  \times I} && R \\
	\\
	{{D^{\leq n}}_{\bar{R}/\mathbb{F}_p}  \times [I]} && {\bar{R}}
	\arrow["eval", from=1-1, to=1-3]
	\arrow["{\textrm{(mod p, mod p)}}"', from=1-1, to=3-1]
	\arrow["{\textrm{mod p}}"', from=1-3, to=3-3]
	\arrow["eval", from=3-1, to=3-3]
\end{tikzcd}
$$
commutes by \cref{corollary4328}. \qedsymbol

\begin{corollary}
    \label[corollary]{corollary43212} 
    With notation as in \cref{theorem4329}
    $$
    supp(\mathcal{I}, m) \subseteq V({\mathcal{D}^{\leq m-1}}_{X/\textrm{Spec} \, \mathbb{Z}} \mathcal{I}).
    $$
    Further, if $\mathcal{I}$ has the property that for all $x$ with corresponding $x_p \in X_p$ 
    $$
    \textrm{ord}_x \mathcal{I} = \textrm{ ord}_{x_p} \mathcal{I}_p,
    $$
    then $x \mapsto \textrm{ord}_x \mathcal{I}$ is an upper-semicontinuous function on $X$.
\end{corollary}

\textit{Proof.} Let $x \in X$ with corresponding $x_p \in X_p$. Suppose $x \in supp(\mathcal{I}, m) $, which implies $x_p \in supp(\mathcal{I}_p, m) \subset X_p$. Then by \cref{theorem4325} and with notation as in \cref{theorem4329} and its proof, $[{{D}^{\leq m-1}}_{R/\mathbb{Z}}I] \subset [\mathfrak{m}] \subset \bar{R}$, where $\mathfrak{m}$ is the unique maximal ideal of $R$. As $(p) \subset \mathfrak{m}$ it follows that 
$
({\mathcal{D}^{\leq m-1}}_{X/\textrm{Spec} \, \mathbb{Z}} \mathcal{I})_x \cong {{D}^{\leq m-1}}_{R/\mathbb{Z}}I \subset \mathfrak{m}.  
$
If $\mathcal{I}$ has the stated property then we show $V({\mathcal{D}^{\leq m-1}}_{X/\textrm{Spec} \, \mathbb{Z}} \mathcal{I}) \subseteq supp(\mathcal{I}, m)$. Let $x \in V({\mathcal{D}^{\leq m-1}}_{X/\textrm{Spec} \, \mathbb{Z}} \mathcal{I})$. By the theorem above, $x_p \in V({\mathcal{D}^{\leq m-1}}_{X_p/\textrm{Spec} (\mathbb{F}_p)} \mathcal{I}_p)$, which implies by \cref{theorem4325} that 
$
\textrm{ord}_x \mathcal{I} = \textrm{ ord}_{x_p} \mathcal{I}_p \geq m. \quad \qed
$

\begin{remark}
    Order being upper-semicontinuous implies the support of a marked ideal is always closed. Thus to compare supports of marked ideals with supports of marked ideals in blow-ups it is enough to check at points closed in their fibers, which we can analyse using local blow-up coordinates. We construct examples of ideals satisfying this property in \cref{C:4_Locally monomial ideals} and compare supports using blow-up coordinates throughout \cref{C:5}.
\end{remark}


\subsection{Ideals locally generated by monomials}
\label{C:4_Locally monomial ideals} 

We now use locally toric simple normal crossings divisors to construct ideals on which order is upper-semicontinuous. These ideals form part of the data for $\Lambda$-marked monomial ideals we introduce in \cref{C:5}.

\begin{definition}
    Let $(X,\mathcal{I}, E)$ be a triple with $X$ smooth over $\textrm{Spec} \, \mathbb{Z}$ with an affine $\Lambda$-cover, $E=\sum_{j \in J} E^j$ a $\mathbb{Z}$-flat snc divisor locally toric on the cover and $\mathcal{I} \subset \mathcal{O}_X$ a coherent ideal. Then $\mathcal{I}$ will be \textit{locally generated by monomials} for $E$ if 
    $$
    \mathcal{I}=\sum_{j=1}^r {\mathcal{I}}_j, \quad {\mathcal{I}}_j = \mathcal{O}_X(-\sum_{i_j \in J_j} a_{i_j} E^{i_j}),
    $$
    where $J_j \subset J$ is a finite index set and $a_{i_j}$ are natural numbers, i.e. $\mathcal{I}$ is generated by finitely many $\mathcal{I}_j$, each of which is a monomial given by finite products in the $E^j$.
\end{definition}

\begin{definition}[Fibers are also generated by monomials]
    Let $(p)$ be a prime ideal of $\mathbb{Z}$, then with the notation as above we have $X_p$ a smooth $\mathbb{F}_p$-variety with $E_p=\sum_i (E^i)_p$ a snc divisor on $X_p$. Then by the local descriptions of $\mathbb{Z}$-flat snc divisors
    $$
    \mathcal{I}_p=\sum_{j=1}^r ({\mathcal{I}}_j)_p, \quad ({\mathcal{I}}_j)_p= \mathcal{O}_{X_p}(-\sum_{i_j \in J_j} a_{i_j} (E^{i_j}_p))
    $$
    and so locally the pullback $\mathcal{I}_p$ will be "generated" by monomials.
    
\end{definition}

We prove some easy facts about ideals locally generated by monomials.

\begin{lemma}
    \label[lemma]{lemma4403} 
    Let $X$ be a smooth $\Lambda$-scheme over $\textrm{Spec} \, \mathbb{Z}$ with an affine $\Lambda$-cover, $E$ a $\mathbb{Z}$-flat simple normal crossings divisor locally toric on the cover and $\mathcal{I}$ an ideal locally generated by monomials for $E$. Let $Z= \cap_{k \in K} E^k$ be the intersection of finitely many $E^i$ and assume that $Z\subset supp(\mathcal{I}, m)$. Then writing $\pi:B_Z X = X' \rightarrow X$ for the blow-up
    $$
    {\pi_{\ast}}^{-1}(\mathcal{I}, m) = \sum_{j=1}^r \mathcal{I}_j \subset \mathcal{O}_{X'}, \quad \mathcal{I}_j = \mathcal{O}_{X'}(-\sum_{i_j \in J_j} a_{i_j} {E^{i_j}}'-a_{i_j'}F),
    $$
    which is an ideal locally generated by monomials for the total transform $E'$. Here ${E^{i_j}}'$ are the strict transforms of $E^{i_j}$, $F$ is the exceptional divisor and $a_{i_j'} := \sum_{i_j: i_j \in K} a_{i_j}-m$. Moreover, the birational transform on an affine of $X'$ on which $E'$ is locally toric can be computed using the substitution rule applied to the monomials on the affine of $X$. 
    
\end{lemma}

\textit{Proof.} We argue locally using blow-up coordinates at closed points. Assume $p \not = 0$. Let $x' \in B_Z X$ closed in $B_{Z_p}X_p$ with $x=\pi(x')$ and consider the local coordinates $x_0, \dots, x_n$ for $x$ and $\bar{x}_0, \dots, \bar{x}_n$ where we are in a $x_1$ chart and $x_1$ corresponds to the exceptional divisor. Then $\mathcal{I}_x = (x^{\alpha^l}, l=1, \dots, s)$ and we can compute ${\pi_{\ast}}^{-1}(\mathcal{I}, m)_{x'}$ using the substitution rule given by the local coordinates on $x'$, i.e. 
$$
{\pi_{\ast}}^{-1}(\mathcal{I}, m)_{x'} = (\bar{x}^{(\alpha^l)'}, l=1, \dots, s), \textrm{ where } (\alpha^l)'_t = 
\begin{cases}
    (\alpha^l)_t & t \not = 1\\
    (\sum_{t: t \in K} (\alpha^l)_t)-m & t = 1
\end{cases}
.
$$
But this is precisely $(\sum_j \mathcal{I}_j)_{x'}$. For the second part we note that the substitution rule can be extended to an affine neighbourhood of $x$ on which $E$ is locally toric. \qedsymbol

\begin{lemma}
    \label[lemma]{lemma4404} 
    Let $\mathcal{I}$ be an ideal locally generated by monomials for a $\mathbb{Z}$-flat snc divisor $E$ and $x \in X$. Then 
    $$
    ord_x (\mathcal{I}) = min_j\{\sum_{i_j \in J_j: \space x \in E^{i_j}} a_{i_j}\}. 
    $$
    Moreover, if  $x$ lies in the fiber over $(p)$ and $x_p \in X_p$ is the point lying over $x$ then $ord_x(\mathcal{I}) = ord_{x_p}(\mathcal{I}_p)$.
\end{lemma}

\textit{Proof.} Locally $\mathcal{I}_x$ will be generated by monomials in $x_{i_j}$ for each $E^{i_j}$ with $x \in E^{i_j}$. The first fact follows from the order for a finitely generated ideal being the minimum of the orders of its generators and the order of a monomial being its degree. The last statement is trivial if $p=0$. If $p \not = 0$ then the statement follows as none of the parameters locally defining $\mathcal{I}$ are $p$ . \qedsymbol

\begin{corollary}
    \label[corollary]{corollary4405} 
    If $\mathcal{I}$ is locally generated by monomials then $x \mapsto \textrm{ord}_x \mathcal{I}$ is an upper-semicontinuous function on $X$ and so for any $m$, $supp(\mathcal{I}, m)$ is a closed subset of $X$.    
\end{corollary}

\textit{Proof.} By the lemma we can apply \cref{corollary43212}. \qedsymbol

\label{4.5}

\subsection{Normal crossing embeddings}
We describe embeddings of $\Lambda$-schemes defined by restrictions of $\mathbb{Z}$-flat locally toric simple normal crossings divisors. These embeddings are the final piece of data needed to define $\Lambda$-marked monomial ideals in the following section.
\label{C4_NormalCrossingsEmbeddings}

Let $X$ be a smooth $\Lambda$-scheme with an affine open $\Lambda$-cover and a $\mathbb{Z}$-flat simple normal crossings divisor $E=\sum_{j \in J} E^j$ that is locally toric on the cover. Let $N \hookrightarrow X$ be a $\mathbb{Z}$-smooth closed $\Lambda$-subscheme with an affine cover coming from the cover of $X$ and induced by $\Lambda$-equivariant quotients. Assume that $N$ has normal crossings to $E$ and is irreducible in order to use observation \ref{observation3305}.

\begin{lemma}
    \label[lemma]{lemma4501} 
    Let $J_N = \{j \in J: N \cap E^j \not = \emptyset \textrm{ or } N \}\subset J$ and $E_N = \sum_{j \in J_N} E^j \subset E$. Then $E_N|_N$ is a $\mathbb{Z}$-flat simple normal crossings divisor on $N$ locally toric on the affine cover of $N$.
\end{lemma}
\textit{Proof.} By observation \ref{observation3305}, $E_N|_N$ will be a $\mathbb{Z}$-flat snc divisor on $N$ and being locally toric follows from the embedding being affine $\Lambda$-equivariant on the respective $\Lambda$-covers. \qedsymbol
\\
\\Let $E_N|_N=(\sum_{j \in J_N} E^j)|_N=\sum_{j \in J_N} (E^j|_N)$ be the $\mathbb{Z}$-flat snc divisor on $N$ as above. Let $J_P \subset J_N$ be finite and define $P=\cap_{j \in J_P} (E^j|_N)$. By the lemma above, \cref{lemma3313} and observation \ref{observation3316}, $P \hookrightarrow N$ will be a $\mathbb{Z}$-smooth $\Lambda$-closed subscheme with an affine $\Lambda$-cover, with $\Lambda$-equivariance induced from the affine covers. For example, if $P=N=X$ then $P$ is an empty intersection of divisors.

\begin{observation}[Local picture]
    \label[observation]{observation4502} 
    Let $x \in P \subset N \subset X$. Then there exists a system of regular parameters $x_0=p,x_1, \dots, x_n$ for $\mathcal{O}_{X, x}$ such that
    $$
    \mathcal{O}_{N,x} \cong \mathcal{O}_{X,x}/(x_1, \dots, x_s), 1 \leq s \leq n, \quad \mathcal{O}_{P,x} \cong \mathcal{O}_{X,x}/(x_1, \dots, x_s, x_{s+1}, \dots, x_{r}), r \leq n.
    $$
\end{observation}

\begin{lemma}[Restrictions] 
    \label[lemma]{lemma4503} 
    Let $Z= \cap_{j \in J_0} E^j \subset X$ for finite $J_0 \subset J$ and suppose $P$ is irreducible. Then there exists a $K_N, K_P \subset J_N$ (not necessarily unique) such that 
    $$
    Z \cap N = \bigcap_{j \in K_N} (E^j|_N),  \quad Z \cap P = \bigcap_{j \in K_P} (E^j|_P)
    $$
    and $\sum_{j \in K_N} (E^j|_N), \sum_{j \in K_P} (E^j|_P)$ are $\mathbb{Z}$-flat snc divisor locally toric on $N$ and $P$ respectively.
\end{lemma}

\textit{Proof.} Let $K_N = J_N-J_0$ and $K_P = K_N-J_P$ then by the local picture it is clear that 
$$
Z \cap N = \bigcap_{j \in J_0} E^j \cap N = \bigcap_{j \in J_0: N \cap E^j \not = \emptyset \textrm{ or N}} (E^j)|_N = \bigcap_{j \in K_N} (E^j)|_N. 
$$
Similarly, 
$$
Z \cap P =  \bigcap_{j \in J_0} E^j \cap N \cap P = \bigcap_{j \in K_N} (E^j)|_N \cap P = \bigcap_{j \in K_N: P \cap (E^j)|_N \not = \emptyset \textrm{ or } P}  (E^j)|_P = \bigcap_{j \in K_P} (E^j)|_P. 
$$ 
As $K_N \subset J_N$ by observation \ref{observation3305}, $\sum_{j \in K_N} (E^j|_N)$ is a $\mathbb{Z}$-flat snc divisor on $N$. Similarly, as $K_P$ is disjoint from $J_P$, $\sum_{j \in K_P} (E^j|_N)$ will meet $P$ transversally in $N$ and hence again by observation \ref{observation3305}, $\sum_{j \in K_P} (E^j|_P)$ is a $\mathbb{Z}$-flat snc divisor on $P$. Being locally toric follows from the fact that $\sum_{j \in J_N} (E^j|_N)$ is locally toric on $N$ and $P \hookrightarrow N$ is induced by affine $\Lambda$-equivariant morphisms. \qedsymbol

\begin{remark}
    \label[remark]{remark4504} 
    If in the above lemma $Z \cap N \subsetneq P$ then $K_P, K_N \not = \emptyset.$
\end{remark}

\begin{corollary} 
    \label[corollary]{corollary4505} 
    With notation as above let $Z= \cap_{j \in J_0} E^j \subset X$ be a blow-up centre such that $J_0 \subset J$ is finite. Then taking the blow-ups $X', N', P'$ by $Z, Z \cap N, Z \cap P$ respectively we have closed affine $\Lambda$-equivariant embeddings $P' \hookrightarrow N' \hookrightarrow X'$ given by the strict transforms. Moreover, these embeddings satisfy the same conditions with respect to the total transform $E'$. 
\end{corollary}
\textit{Proof.} By the above two lemmas $Z \cap N, Z \cap P$ are blow-up centres as in \cref{example3306} for the respective $\mathbb{Z}$-smooth $\Lambda$-schemes. It follows from \cref{corollary3314} that $P', N', X'$ are $\mathbb{Z}$-smooth $\Lambda$-schemes with affine $\Lambda$-covers. The embeddings being $\Lambda$-equivariant and coming from affine $\Lambda$-rings follow from \cref{corollary2322}. Now $E$ being a $\mathbb{Z}$-flat snc divisor locally toric on $X$ implies its total transform $(E)'$ is $\mathbb{Z}$-flat locally toric by \cref{lemma3315} and similarly so is the total transform $(E_N)'$. Moreover, by the observation of the local picture and \cref{lemma3227} we can see that $N'$ has normal crossings with $E'$ and is transverse to the exceptional divisor. It follows that the components of $(E_N)'$ are precisely the components of $E'$ that are transverse to $N'$. Again using \cref{lemma3227} the strict transform $P'$ will be given by the closed subscheme $\cap_{j \in J_P} ((E^j)')|_{N'}$ with $J_P \subset J_{N'}=J_N \cup j$, where $j$ indexes the exceptional divisor. Thus $P' \hookrightarrow N'$ is given by a finite intersection of components of $E_{N'}|_{N'}$. \qedsymbol

\section{Marked monomial ideals}
\label{C:5}

In this section we define $\Lambda$-marked monomial ideals and prove the existence of order reduction using classical resolution methods and the results of the previous sections. Our $\Lambda$-marked monomial ideals will be tuples consisting of closed $\Lambda$-subschemes, an ideal locally generated by monomials for a simple normal crossings divisor $E$ and a non-negative integer $m$. This differs from marked monomial ideals of \cite{BM1} as we require an $E$ to generalise codimension one orbit closures. Similarly, the inclusion of embedding data differs from \cite{10.2307/j.ctt7rptq}, which we discuss in \cref{C:6} and allows us to give stronger resolutions for certain embedded $\Lambda$-schemes. We will use the notation of \cite{BM1} and make use of the results on marked ideals developed in \cref{C:4} by considering restrictions to closed subschemes. 


\subsection{First definitions}
\label{C5_1}

\begin{definition}
    \label[definition]{definition5101}  
    A $\Lambda$-marked monomial ideal is a tuple 
    $$
    \underline{\mathcal{I}} := (X, N, P, \mathcal{I}, E, m),
    $$
    where
    \begin{itemize}
        \item $X$ is a smooth $\Lambda$-scheme over $\mathbb{Z}$ with an affine $\Lambda$-cover (see \cref{AffineLambdaCover}) that is connected and has connected fibers.
        \item $E =\sum_{j \in J} E^j$ is a $\mathbb{Z}$-flat simple normal crossings divisor on $X$ that is locally toric for  an affine $\Lambda$-cover of $X$ (see \cref{C3_S4},  \cref{definition3311}). Moreover, $J$ is totally ordered.
        \item $N$ is a smooth $\Lambda$-scheme over $\mathbb{Z}$ that is a closed $\Lambda$-subscheme of $X$ with normal crossings to $E$ as in \cref{4.5}. We further assume that the affine $\Lambda$-cover of $N$ induced by the embedding $N \hookrightarrow X$ is intersection distinguished (see \cref{definition3321}) for $E_N|_N$ as in \cref{lemma4501}. Finally, $N$ is also connected and has connected fibers.
        \item $P$ is a smooth $\Lambda$-subscheme of $N$ given by intersecting $N$ with a finite collection of $E^j$ as in \cref{4.5}. Again, $P$ is connected and has connected fibers.
        \item 
        $$
        \mathcal{I}=\sum_{j=1}^r {\mathcal{I}}_j, \quad {\mathcal{I}}_j = \mathcal{O}_X(-\sum_{i_j \in J_j} a_{i_j} E^{i_j}),
        $$
        is an ideal locally generated by monomials for $E$ on $X$ such that $N$ is transverse to each $E^{i_j}$ and no $E^{i_j}$ contain $P$ (see \cref{C:4_Locally monomial ideals} and \cref{C4_NormalCrossingsEmbeddings}). 
        \item $m$ \textrm{ is a positive integer.}
    \end{itemize}
\end{definition}

These $\Lambda$-marked monomial ideals generalise marked monomial ideals for toric $\Lambda$-schemes introduced in upcoming work, which themselves are modelled on marked monomial ideals as in Bierstone and Milman (see \cite[Section 8]{BM1}). Here the $E^j$ play the role of codimension one orbit closures, which allow us to construct blow-ups that transform this data in a controllable, $\Lambda$-equivariant manner. Note $X,N,P$ smooth over $\mathbb{Z}$ and connected imply they are irreducible.

\begin{definition}
    (See \cite[Section 8.1]{BM1}).
    Let $\underline{\mathcal{I}}=(X, N, P, \mathcal{I},E, m)$ be a $\Lambda$-marked monomial ideal. The \textit{support of $\underline{\mathcal{I}}$} is defined as 
    $$
    \textrm{supp}\underline{\mathcal{I}} := \{x \in P: ord_x  \mathcal{I}\geq m\}. 
    $$
    If for all $x \in P$, $ord_x (\mathcal{I}) \leq m$ then we will say $\underline{\mathcal{I}}$ is of \textit{maximal order}, i.e. either $\textrm{supp}\underline{\mathcal{I}} = \emptyset$ or $m=\textrm{max-ord }\mathcal{O}_P\mathcal{I}$.
\end{definition}

\begin{lemma}
    \label[lemma]{lemma5104} 
    Let $x \in P$. Then 
    $$
    ord_x \mathcal{O}_P \mathcal{I} = ord_x \mathcal{O}_N \mathcal{I} = ord_x  \mathcal{I},
    $$  
    where $\mathcal{O}_P \mathcal{I} \subset \mathcal{O}_P,\mathcal{O}_N \mathcal{I} \subset \mathcal{O}_N$ are the restrictions (as ideals) of the ideal $\mathcal{I}$ to $P, N$ respectively.
\end{lemma}

\textit{Proof.} There exist a system of regular parameters $x_0=p,x_1, \dots, x_n$ (no $x_0$ if $p = 0$ see observation \ref{observation3225}) for $\mathcal{O}_{X, x}$ such that $\mathcal{O}_{N,x} \cong \mathcal{O}_{X,x}/(x_1, \dots, x_s), 1 \leq s \leq n$, $\mathcal{O}_{P,x} \cong \mathcal{O}_{X,x}/(x_1, \dots, x_s, x_{s+1}, \dots, x_{r}), r \leq n$ and each $E^{i_j}$ in $\mathcal{I}$ containing $x$ is cut out by $x_k \in \mathcal{O}_{X,x}, k\geq r+1$. Then $\mathcal{I}_x$ is generated by monomials $x^{\alpha}$ in $x_{r+1}, \dots, x_{n}$, each of which remain regular parameters in $\mathcal{O}_{N,x}, \mathcal{O}_{P,x}$. \qedsymbol

\begin{corollary}
    \label[corollary]{corollary5106} 
    Let $\underline{\mathcal{I}}=(X, N, P, \mathcal{I},E, m)$ be a $\Lambda$-marked monomial ideal and $m'$ an arbitrary non-negative integer. 
    \begin{enumerate}
        \item $supp(\mathcal{O}_P \mathcal{I}, m') = supp(\mathcal{O}_N \mathcal{I}, m') \cap P = supp(\mathcal{I}, m') \cap P$,
    where the supports of the marked ideals are taken in $P,N,X$ respectively.
        \item $\textrm{supp}\underline{\mathcal{I}} = supp(\mathcal{O}_P \mathcal{I}, m) = supp(\mathcal{O}_N \mathcal{I}, m) \cap P = supp(\mathcal{I}, m) \cap P.$
        \item $\textrm{supp}\underline{\mathcal{I}} \subset X$ is closed and so can be analysed using local coordinates.
    \end{enumerate}
\end{corollary}

\textit{Proof.} (1) is clear from the lemma. (2) follows from (1). (3) follows as $P \subset X$ is closed and applying \cref{corollary4405} to $\mathcal{I} \subset \mathcal{O}_X$. \qedsymbol

\begin{definition}
    Let $\underline{\mathcal{I}}=(X, N, P, \mathcal{I},E, m)$ be a $\Lambda$-marked monomial ideal. A \textit{permissible} blow-up centre of $\underline{\mathcal{I}}$ is a subscheme $Z = \bigcap_{j \in J_0} E^j $ for a finite set $J_0 \subset J$ such that $Z \cap N \subset \textrm{supp}\underline{\mathcal{I}}$. 
\end{definition}

\begin{remark}
    \label[remark]{remark5212} 
    For a permissible centre $Z = \bigcap_{j \in J_0} E^j $ and by the definition of $\mathcal{I}$ there will be at least one $E^{i_j}$ involved in $\mathcal{I}_j$ such that $Z \cap N \subset E^{i_j}$. By \cref{remark4504}, either $Z \cap N = \emptyset$ or the $K_N, K_P$ of \cref{lemma4503} are non-empty. Moreover, as $E^{i_j}$ does not contain $P$, $Z \cap N \not = P$ and it follows that $Z$ does not contain the generic point of $M$ and $Z \cap N$ does not contain the generic point of $N$ or $P$. 
\end{remark}

\begin{lemma}
    \label[lemma]{lemma5213} 
    Let $Z = \bigcap_{j \in J_0} E^j $ be a permissible centre. Then $Z \subset supp(\mathcal{I}, m)$.
\end{lemma}

\textit{Proof.} Recall from \cref{lemma4503} that there exists $K_N \subset J_0$ such that $Z \cap N = \bigcap_{j \in K_N} (E^j|_N)$. Let $\eta \in Z \cap N$ be a generic point of a component. A regular sequence of parameters for $\mathcal{O}_{N, \eta}$ is $x_1, \dots, x_n$ where each $x_j$ corresponds to a  $E^j|_N, j \in K_N$ and $(\mathcal{O}_N\mathcal{I})_{\eta} = (x^{\alpha^j}, j \in J_1)$. But $Z$ permissible implies by \cref{lemma4404} that
$$
\textrm{ord}_{\eta} \mathcal{O}_N\mathcal{I} = min_j\{\sum_{i_j: \eta \in E^{i_j}} a_{i_j}\} \geq m.
$$
It follows that for every $x \in \bigcap_{j \in K_N} E^j$, $\textrm{ord}_x \mathcal{I} \geq m$. In particular,  $Z \subset \bigcap_{j \in K_N} E^j$. \qedsymbol

\begin{definition}
    (See \cite[Definition 8.3]{BM1}).
    Let $\pi_Z: X' \rightarrow X$ be the blow-up of $X$ by a permissible centre $Z$ for $\underline{\mathcal{I}}$. We define the \textit{transform} $\underline{\mathcal{I}}'$ of $\underline{\mathcal{I}}$ by $\pi$ as the tuple $\underline{\mathcal{I}}'=(X', N', P', \mathcal{I}', E', m)$, where
    \begin{itemize}
        \item $X'$ is the blow-up.
        \item $E'$ is the total transform of $E$ by $\pi$.
        \item $N'$ is the strict transform of $N$ by $\pi$. 
        \item $P'$ is the strict transform of $P$ by $\pi$.
        \item $\mathcal{I}'= \pi^{-1}_{\ast}(\mathcal{I}, m)$ is the birational transform of the marked ideal $(\mathcal{I}, m)$ (see \cref{C4:Markedidealsbirationaltransforms}).
    \end{itemize}
\end{definition}

\begin{lemma}
    \label[lemma]{lemma5215} 
    $\underline{\mathcal{I}}'$ is a $\Lambda$-marked monomial ideal.
\end{lemma}

\textit{Proof.} By the definition of $\mathcal{I}$ and \cref{lemma5213}, $Z$ is a blow-up centre satisfying the conditions of \cref{corollary4505} and permissible for the marked ideal $(\mathcal{I}, m)$ (see \cref{definition4204}). It follows from \cref{corollary4505} that $P'\hookrightarrow N' \hookrightarrow X'$ satisfies the embedding conditions. Further, by \cref{lemma4403} and its proof, $\mathcal{I}'$ is an ideal locally generated by monomials for $E'$, a totally ordered $\mathbb{Z}$-flat simple normal crossings divisor by \cref{definition3309} and \cref{lemma3315}. The components of $E'$ generating $\mathcal{I}'$ satisfy the conditions on $N', P'$ by the arguments in the proof of \cref{corollary4505}. The induced affine $\Lambda$-cover of $N'$, which is compatible with the cover of $X'$, is intersection distinguished by \cref{lemma3324}. By \cref{remark5212}, $X', N', P'$ are irreducible and hence connected. Finally, the fibers of $X', N', P'$ are connected by \cref{remark5223} and \cref{lemma5224} below. \qedsymbol

\begin{definition}
    \label[definition]{definition5216} 
    A \textit{permissible} sequence of blowings up for $\underline{\mathcal{I}}$ is a sequence of blowings-up
    $$
    X=X_0 \xleftarrow[]{\pi_1} X_1 \leftarrow \dots \xleftarrow[]{\pi_{t+1}} X_{t+1},
    $$
    where each $\pi_{j+1}, j = 0, \dots, t$ is a permissible blow-up for $\underline{\mathcal{I}}_j = (X_j, N_j, P_j, \mathcal{I}_j, m)$ and $\underline{\mathcal{I}}_{j+1} = (X_{j+1}, N_{j+1}, P_{j+1}, \mathcal{I}_{j+1}, m)$ is the transform of $\underline{\mathcal{I}}_j$ by $\pi_{j+1}$, with $\underline{\mathcal{I}}_0 = \underline{\mathcal{I}}$. An \textit{order reduction} of $\underline{\mathcal{I}}$ is a sequence of permissible blowings-up such that $supp \underline{\mathcal{I}}_{t+1} = \emptyset$. By definition a permissible blow-up sequence is $\Lambda$-equivariant.
\end{definition}

Let $Z = \bigcap_{j \in J_0} E^j $ be a permissible centre for $\underline{\mathcal{I}}$. Then by \cref{corollary5106} (2) $Z \cap N$ is a permissble centre for the marked ideal $(\mathcal{O}_P \mathcal{I}, m)$ and we may take the birational transform and denote this by $(\mathcal{O}_P \mathcal{I})' \subset \mathcal{O}_{P'}.$

\begin{lemma}
    \label[lemma]{lemma5216} 
    Let $\mathcal{I}' \subset \mathcal{O}_{X'}$ be the ideal underlying $\underline{\mathcal{I}}'$. Then 
    $$
    \mathcal{O}_{P'}\mathcal{I}' \cong (\mathcal{O}_P \mathcal{I})'.
    $$
\end{lemma}

\textit{Proof.} Let $x' \in P'$ be a closed point in the fiber over $(p) \not = 0$ lying over $x \in P$ and $x \in Z$. Then by our definitions there exist a regular sequence of parameters 
$$
x_0=p, x_1, \dots, x_r, \dots, x_s, \dots, x_l, \dots, x_n,  \quad r < s \leq l \leq n
$$ 
for $\mathcal{O}_{X,x}$ such that $Z$ is defined by $x_1, \dots, x_s$, $\mathcal{O}_{P,x} \cong \mathcal{O}_{X,x}/(x_1, \dots, x_r)$ and $\mathcal{I}_x$ is defined by monomials in $x_{r+1}, \dots, x_l$. Suppose $x'$ is contained in a blow-up chart corresponding to $x_i$ and note that as $x' \in P'$, $r < i \leq s$. Then there are local blow-up coordinates on $\mathcal{O}_{X', x'}$ as in \cref{lemma3227} such that $\mathcal{O}_{P', x'} \cong \mathcal{O}_{X',x'}/(\bar{x}_1, \dots, \bar{x}_r)$, the remaining parameters are the blow-up coordinates of $\mathcal{O}_{P,x}$ and 

$$
\begin{tikzcd}[cramped]
	{\mathcal{O}_{P',x'} \cong \mathcal{O}_{X',x'}/(\bar{x}_1, \dots, \bar{x}_1))} && {\mathcal{O}_{X,x}} \\
	\\
	{\mathcal{O}_{P,x} \cong \mathcal{O}_{X,x}/(x_1, \dots, x_r))} && {\mathcal{O}_{X,x}}
	\arrow[two heads, from=1-3, to=1-1]
	\arrow["{\pi_P}"', from=3-1, to=1-1]
	\arrow["\pi"', from=3-3, to=1-3]
	\arrow[two heads, from=3-3, to=3-1]
\end{tikzcd}
$$
commutes. But for each monomial $x^{\alpha} \in \mathcal{I}_x \subset \mathcal{O}_{X,x}$ the image in $\mathcal{O}_{P,x}$, which will generate $(\mathcal{O}_P\mathcal{I})_x$, is the same monomial in the parameters of $\mathcal{O}_{P,x}$. It follows from the commutativity of the diagram and the way we compute birational transforms (see \cref{definition4205}) that $(\mathcal{O}_{P'}\mathcal{I}')_{x'} \cong ((\mathcal{O}_P \mathcal{I})')_{x'}$. As this holds in the neighbourhood of every closed point of $P'$ the lemma is proved. \qedsymbol

\begin{corollary}
    \label[corollary]{corollary5217} 
    $supp(\underline{\mathcal{I}}') = supp(\mathcal{O}_{P'}\mathcal{I}', m) = supp((\mathcal{O}_P \mathcal{I})', m)$.
\end{corollary}

\begin{remark}
    \label[remark]{remark5218} 
    This corollary allows us to compare the supports of $\Lambda$-marked monomial ideals after transforms by just considering the birational transform of the restriction to $P$, which can be analysed as in \cref{C4:Markedidealsbirationaltransforms} by using local blow-up coordinates and applying the substitution rule.
\end{remark}

\begin{definition}
    Let $(p)$ be a prime ideal of $\mathbb{Z}$ and $\underline{\mathcal{I}}$ a $\Lambda$-marked monomial ideal. Define the tuple
    $$
    \underline{\mathcal{I}}_p = (X_p, N_p, P_p, \mathcal{I}_p,E_p, m),
    $$
    where $X_p, N_p, P_p$ are the smooth connected (hence irreducible) varieties over $\mathbb{F}_p$ given by the fibers over $(p)$, $\mathcal{I}_p \subset \mathcal{O}_{X_p}$ is the pullback ideal and $E_p$ is the pullback divisor, which by \cref{lemma3308} is a simple normal crossings divisor on $X_p$. Note that as $P_p$ is connected each $E^i$ not containing $P$ will correspond to a $(E^i)_p$ not containing $P_p$ using local descriptions. Define 
    $$
    \textrm{supp}\underline{\mathcal{I}}_p := \{x \in P_p: ord_x  \mathcal{I}_p\geq m\}.
    $$
\end{definition}

\begin{observation}
    \label[observation]{observation5222}  
    Let $(p)$ be a prime ideal of $\mathbb{Z}$ and $x \in X$ contained in the fiber over $(p)$ with corresponding point $x_p$. Then as $\mathcal{I}$ is locally generated by monomials, $x \in \textrm{supp}\underline{\mathcal{I}} \iff x_p \in \textrm{supp}\underline{\mathcal{I}}_p$. It follows that as sets $\textrm{supp}\underline{\mathcal{I}} = \bigsqcup_{(p) \in \textrm{Spec} \, \mathbb{Z}} \textrm{supp}\underline{\mathcal{I}}_p$.
\end{observation} 

Let $Z = \cap E^j$ be a permissible blow-up centre for $\underline{\mathcal{I}}$ and $(p)$ a prime ideal of $\mathbb{Z}$. Then the pullback $Z_p = \cap_j (E^j)_p \subset X_p$ will be a smooth $\mathbb{F}_p$-variety such that $Z_p \subset supp(\mathcal{I}_p, m) \subset X_p$ and we have the blow-up $\pi_p: B_{Z_p}X_p \rightarrow X_p$. Define the transform $(\underline{\mathcal{I}}_p)'$ of $\underline{\mathcal{I}}_p$ by $\pi_p$ as the tuple $(\underline{\mathcal{I}}_p)'=((X_p)', (N_p)', (P_p)', (\mathcal{I}_p)', (E_p)', m)$, where
    \begin{itemize}
        \item $(X_p)'$ is the blow-up.
        \item $(E_p)'$ is the total transform of $E_p$ by $\pi_p$ (see \cref{corollary33012}).
        \item $(N_p)'$ is the strict transform of $N_p$ by $\pi_p$. 
        \item $(P_p)'$ is the strict transform of $P_p$ by $\pi_p$.
        \item $(\mathcal{I}_p)'= (\pi_p)^{-1}_{\ast}(\mathcal{I}_p, m)$ is the birational transform of $(\mathcal{I}_p, m)$ by $\pi_p$ (see \cref{definition4203}).
    \end{itemize}

\begin{remark}
    \label[remark]{remark5223} 
    By the same arguments as in \cref{remark5212}, $Z_p$ will not contain the generic points of $X_p$ and $Z_p \cap N_p$ will not contain the generic points of $N_p$ or $P_p$. It follows that $(X_p)', (N_p)', (P_p)'$ are smooth connected $\mathbb{F}_p$-varieties.   
\end{remark}

\begin{lemma}[Transforms and fibers]
    \label[lemma]{lemma5224} 
    Let $Z = \bigcap_{j \in J_0} E^j$ be a permissible blow-up centre for $\underline{\mathcal{I}}$. Let $\underline{\mathcal{I}}'$ be the transform of $\underline{\mathcal{I}}$ and $(\underline{\mathcal{I}}_p)'$ be the transform of $\underline{\mathcal{I}}_p$. Then 
    $$
    (\underline{\mathcal{I}}')_p \cong (\underline{\mathcal{I}}_p)'
    $$
\end{lemma}

\textit{Proof.} $(X_p)' \cong (X')_p, (N_p)' \cong (N')_p, (P_p)' \cong (P')_p$ by \cref{lemma3226}, $(E_p)' \cong (E')_p$ by \cref{corollary33012} and $(\mathcal{I}')_p \cong (\mathcal{I}_p)'$ by \cref{lemma4207}. \qedsymbol

\begin{observation}[Fibers of order reduction]
    Let $\underline{\mathcal{I}}$ be a $\Lambda$-marked monomial ideal and $(p)$ be any prime ideal of $\mathbb{Z}$. Suppose we have an order reduction of $\underline{\mathcal{I}}$, i.e. a sequence of permissible blow-ups of $X$
    $$
    X=X_0 \xleftarrow[]{\pi_1} X_1 \leftarrow \dots \xleftarrow[]{\pi_{t+1}} X_{t+1},
    $$
    such that $supp \underline{\mathcal{I}}_{t+1} = \emptyset$. Taking the fiber over $(p)$ we get a sequence of blow-ups of $X_p$
    $$
    X_p=(X_0)_p \xleftarrow[]{\pi_1} (X_1)_p \leftarrow \dots \xleftarrow[]{\pi_{t+1}} (X_{t+1})_p,
    $$
    such that by \cref{lemma4404}, $supp (\underline{\mathcal{I}}_p)_{t+1} = \emptyset$ (see \cref{C5_1}).
\end{observation}

We describe how to sum $\Lambda$-marked monomial ideals together and properties arising from \cref{431}.

\begin{definition}
    \label[definition]{definition5231} 
    (See \cite[Definition 8.8]{BM1} for the case of the sum of two ideals).
    Let $\underline{\mathcal{I}}_i = (X, N, P, \mathcal{I}_i, E, m_i), i = 1, \dots, k$ be a collection of $\Lambda$-marked monomial ideals with $X,N,P,E$ all the same. As before let $m_{\hat{i}} = \Pi_{j \not = i} m_j$ and $m = \Pi_i m_i$. Define the sum
    $$
    \sum_i \underline{\mathcal{I}}_i = (X, N, P, \sum_i ((I_i)^{m_{\hat{i}}}), E, m),
    $$
    which is a $\Lambda$-marked monomial ideal as $\mathcal{I}_i$ are locally generated by monomials for $E$ and so the same holds for their sums and products.
\end{definition}

\begin{remark}
    \label[remark]{remark5232} 
    Our definition of the finite sum of $\Lambda$-marked monomials generalises Bierstone and Milman's sum of two marked monomial ideals in \cite[Definition 8.8]{BM1} and gives an explicit definition for the ideal. The sum defined in \cite[Definition 8.8]{BM1} is not associative and so we have options for how we could define the sum of more than two ideals using their definition. This is not an issue for the purposes of order reduction as the supports will be the same. However, for convenience in later arguments we give a precise definition.  
\end{remark}

\begin{lemma}   
    \label[lemma]{lemma5233} 
    $$
    supp\sum_i \underline{\mathcal{I}}_i = \bigcap_i supp\underline{\mathcal{I}}_i \subset P.
    $$
\end{lemma}

\textit{Proof.} By \cref{corollary5106} (2)
$$
supp\sum_i \underline{\mathcal{I}}_i = supp(\mathcal{O}_P(\sum_i(\mathcal{I}_i)^{m_{\hat{i}}}), m), \quad supp\underline{\mathcal{I}}_i = supp(\mathcal{O}_P\mathcal{I}_i, m_i) \subset P.
$$
But by using the local description of the ideals in $\Lambda$-marked monomial ideals and observation \ref{observation4502},
$$
\mathcal{O}_P(\sum_i(\mathcal{I}_i)^{m_{\hat{i}}}) = \sum_i(\mathcal{O}_P\mathcal{I}_i)^{m_{\hat{i}}} \subset \mathcal{O}_P.
$$
As $P$ is a regular Noetherian scheme, by \cref{corollary4213}
$$
supp(\sum_i(\mathcal{O}_P\mathcal{I}_i)^{m_{\hat{i}}}, m) = \bigcap_i supp((\mathcal{O}_P\mathcal{I}_i, m_i) \subset P
$$
and thus
\begin{align*}
     supp\sum_i \underline{\mathcal{I}}_i &= supp(\mathcal{O}_P(\sum_i(\mathcal{I}_i)^{m_{\hat{i}}}), m)
     \\&=supp(\sum_i(\mathcal{O}_P\mathcal{I}_i)^{m_{\hat{i}}}, m) = \bigcap_i supp(\mathcal{O}_P\mathcal{I}_i, m_i) 
     \\&= \bigcap_i supp\underline{\mathcal{I}}_i \subset P. \quad \qed
\end{align*}

\begin{corollary}
    A blow-up centre $Z$ is permissible for $\sum_i \underline{\mathcal{I}}_i$ if and only if it is permissible for each $\underline{\mathcal{I}}_i$.
\end{corollary}

\begin{lemma}[Sums commute with transforms]
    \label[lemma]{lemma5235} 
    Let $Z$ be a permissible centre for $\sum_i \underline{\mathcal{I}}_i$ and hence each $\underline{\mathcal{I}}_i$. Write $(\sum_i \underline{\mathcal{I}}_i)'$ for the transform of the sum and $(\underline{\mathcal{I}}_i)'$ for the transform of each summand. Then
    $$
    (\sum_i \underline{\mathcal{I}}_i)' = \sum_i (\underline{\mathcal{I}}_i')
    $$
\end{lemma}

\textit{Proof.} As summing $\Lambda$-marked monomial ideals only changes the underlying ideal and marked integer the proof follows from \cref{lemma4214}. \qedsymbol

\begin{corollary}
    \label[corollary]{corollary5236} 
    $$
    supp((\sum_i \underline{\mathcal{I}}_i)') = \bigcap_isupp(\underline{\mathcal{I}}_i')
    $$
\end{corollary}

\textit{Proof.} Follows from \cref{lemma5233}. \qedsymbol

\begin{lemma}
    \label[lemma]{lemma52316} 
    (See \cite[Lemma 8.9]{BM1} for the sum of two marked ideals).
    Let $\underline{\mathcal{I}}_i$ be a collection of $\Lambda$-marked monomial ideals. 
    \begin{itemize}
        \item A blow-up sequence for $X$ is permissible for $\sum_i \underline{\mathcal{I}}_i$ if and only if it is permissible for each $\underline{\mathcal{I}}_i$.
        \item  Writing $\underline{\mathcal{I}}_i' ,(\sum_i \underline{\mathcal{I}}_i)'$ for the transforms under such a sequence we have 
        $$
        (\sum_i \underline{\mathcal{I}}_i)' = \sum_i (\underline{\mathcal{I}}_i') \textrm{ and } 
        supp((\sum_i \underline{\mathcal{I}}_i)') = \bigcap_isupp(\underline{\mathcal{I}}_i'). 
        $$
    \end{itemize}
    
\end{lemma}

\textit{Proof.} Combine the above lemmas and corollaries and use induction. \qedsymbol

\subsection{Order reduction}
\label{C:5_Sec:Analysis}

Let $\underline{\mathcal{I}}=(X, N, P, \mathcal{I},E, m)$ be a $\Lambda$-marked monomial ideal. We prove the existence of an order reduction for $\underline{\mathcal{I}}$ by induction on $dim_{\mathbb{Z}}P$, the relative dimension of $P$ over $\mathbb{Z}$. First suppose $dim_{\mathbb{Z}}P=0$. If $\mathcal{I}=0$ then $supp\underline{\mathcal{I}} = P$ and we can blow up $Z_P$ as in \cref{4.5} such that $Z_P \cap N = P$ and give an order reduction. If $\mathcal{I} \not = 0$ then for all prime ideals $(p) \subset \mathbb{Z}$, $dim_{\mathbb{F}_p}P_p = 0$ and $P_p$ consists of a single point. But each $(E^i)_p$ defining $\mathcal{I}_p$ does not contain $P_p$ so $\mathcal{O}_{P_p}\mathcal{I}_p$ is the unit ideal and $supp(\mathcal{I}_p, m) = \emptyset$. As this holds for all $(p)$, $supp\underline{\mathcal{I}} = \emptyset$. Now assume $dim_{\mathbb{Z}}P > 0$. We use the same standard two step induction as in \cite[3.70]{10.2307/j.ctt7rptq}, \cite[Section 8]{BM1} and \cite[Section 5]{bierstone2007functoriality}:
\begin{gather*} 
  \textrm{order reduction for general } \underline{\mathcal{I}} \textrm{ in dimensions $\leq$ n-1} \\
                                 \Downarrow                                        \\
    \textrm{order reduction for maximal order } \underline{\mathcal{I}} \textrm{ in dimension n} \\
                                 \Downarrow                                        \\
    \textrm{order reduction for general } \underline{\mathcal{I}} \textrm{ in dimension n}.
\end{gather*}
We prove the implications in the induction in the following three subsections. In \cref{C:5_Sec:Analysis:FM} we prove order reduction for an $\underline{\mathcal{I}}$ of maximal order assuming the existence of order reduction when $dim_{\mathbb{Z}}P \leq n-1$. The second implication is proved in two parts. In \cref{C5_S3_S2} we show that assuming order reduction for maximal order $\underline{\mathcal{I}}$ of $dim_{\mathbb{Z}}P=n$ we can reduce order reduction to the monomial case, which we prove in \cref{C5_S3_S3}.


\subsubsection{Order reduction in the case of maximal order}
\label{C:5_Sec:Analysis:FM} 
In this section we show that we can do order reduction for $\Lambda$-marked monomial ideals of maximal order. We use the "by hand" method of Bierstone and Milman \cite[step 1]{BM1} rather than the derivative ideal method of Kollar (and Hironaka) in characteristic $0$. The reason for this is that we run into similar problems as resolution in characteristic $p$; derivatives of ideals do not behave as nicely as in characteristic $0$. We discuss this difference in the context of monomial derivatives. 
\\
\\The reason we require $N$ to have an affine $\Lambda$-cover that is intersection distinguished is the following lemma.

\begin{lemma}
    \label[lemma]{lemma5311} 
    Let $\underline{\mathcal{I}}$ be a $\Lambda$-marked monomial ideal of maximal order and $\textrm{Spec} \, A$ an affine open in the $\Lambda$-cover of $X$ such that $supp(\underline{\mathcal{I}}) \cap \textrm{Spec} \, A \not = \emptyset$. Let $I =\mathcal{I}(\textrm{Spec} \, A)$, which will be generated by monomials in $x_1, \dots, x_n \in A$ that correspond to the $E^{i_j} \cap \textrm{Spec} \, A \not = \emptyset$, i.e. $I=(x^{\alpha^j}:j \in J)$. Then $|\alpha^j| \geq m$ for all $j$ and there exists at least one $j$ with $|\alpha^j|=m$.
\end{lemma}

\textit{Proof.} By definition each $x^{\alpha^j}$ remains a monomial in the ring $A/I(N)$ underlying $N \cap \textrm{Spec}  (A)$, i.e. none of the $x_1, \dots, x_n$ vanish and will generate $\mathcal{O}_N\mathcal{I}(N \cap \textrm{Spec} \, A)$. $N \cap \textrm{Spec} \, A$ being intersection complete for $E_N|_N$ implies there exists $a \in N \cap \textrm{Spec} \, A$ such that $a \in \bigcap_{j \in J_N} (E^j)|_N$. It follows that $a \in P \cap \textrm{Spec} \, A$ and $\textrm{max-ord }\mathcal{O}_P\mathcal{I}|_{\textrm{Spec} \, A} = \textrm{ord}_a\mathcal{O}_P\mathcal{I} = \textrm{ord}_a\mathcal{O}_N\mathcal{I}=m$ as otherwise $supp(\underline{\mathcal{I}}) \cap \textrm{Spec} \, A = \emptyset$. Then $|\alpha^j| \geq m$ for all $j$ and there exists at least one $j$ with $|\alpha^j|=m$. \qedsymbol
\\


\paragraph{\textbf{Constructing the associated ideal $\underline{\mathcal{C}}(\underline{\mathcal{I}})$.}}
Let $\underline{\mathcal{I}}$ be a $\Lambda$-marked monomial ideal of maximal order and $U_i=\textrm{Spec} \, A_i$ an affine open in the $\Lambda$-cover. Then as our embeddings and divisor $E$ are $\Lambda$-equivariant with respect to the affine cover of $X$
$$
\underline{\mathcal{I}}_i := (U_i, N_i=N|_{U_i}, P_i=P|_{U_i}, \mathcal{I}_i:= \mathcal{I}|_{U_i},  E_i : =E|_{U_i}, m)
$$
will be a $\Lambda$-marked monomial ideal with $dim_{\mathbb{Z}} P_i \leq dim_{\mathbb{Z}} P$. By \cref{lemma5311}, $I=\mathcal{I}(U_i) = (x^{\alpha^j}:j \in J) \subset A$ with $|\alpha^j| \geq m$ for all $m$ and $|\alpha^j|=m$ for at least one $j$. For all $x$-variables occuring in a $x^{\alpha^j}$ with $|\alpha^j|=m$ label these as $z$-variables and the remaining $x$-variables as $w$-variables so that each monomial of $I$ is of the form $x^{\alpha^j}=z^{\eta^j}w^{\zeta^j}$ and there are some monomials with $\zeta^j=\bar{0}, \eta^j=\alpha^j$. By definition each $z_i$ corresponds to an $E^i$ not containing $P$. Then as $P$ is connected $E^i|_{U_i}$ does not contain $P_i$ and so $(z_i)|_{P_i}$ is a $\mathbb{Z}$-flat snc divisor on $P_i$. We construct a $\Lambda$-marked monomial ideal using this data and the following lemma, which will allow us to use the inductive assumption and give an order reduction for $\underline{\mathcal{I}}_i$. 

\begin{lemma} 
    \label[lemma]{lemma5312} 
    Consider the closed subscheme $Y_i=\textrm{Spec}  (B_i:= C_i/L)$, where $C_i$ is the ring underlying $P_i$ and  $L$ is generated by all $z$-variables.
    \begin{enumerate}
        \item $Y_i$ is a smooth $\Lambda$-scheme with an affine $\Lambda$-cover.
        \item For $\zeta^j \not = \bar{0}$, $(w^{\zeta^j})$ is an ideal locally generated by monomials for $E_i$ satisfying the conditions on $N_i, Y_i$.
        \item If $m-|\eta^j| > 0$ the tuple $(Y_i,N_i, U_i (w^{\zeta^j}),  m-|\eta^j| ,E_i)$ is a $\Lambda$-marked monomial ideal.
        \item For any $x \in Y_i\subset P_i$, $\mathcal{O}_{Y_i, x} \cong \mathcal{O}_{P_i, x}/(x_1, \dots, x_s)$, where $x_1, \dots , x_s$ are regular parameters corresponding to the $z$-variables.
    \end{enumerate}
\end{lemma}

\textit{Proof.} $Y_i$ is a closed subscheme of $N_i$ as in \cref{4.5}, where with the notation of that section $J_{P_i} \subset J_{Y_i} \subset J_{N_i}$. It follows that $Y_i$ is a smooth $\Lambda$-scheme over $\textrm{Spec} \, \mathbb{Z}$ with the desired embedding into $N_i$. $(w^{\zeta^j})$ is clearly locally generated by monomials for $E_i$ and already satisfies transversality for $N_i$ and non-containment for $P_i$. As the $w$-variables and $z$-variables are disjoint it satisfies non-containment for $Y_i$. (3) follows from (1) and (2). (4) is clear from the local description of snc divisors. \qedsymbol
\\
\\Write $m_j = m-|\eta^j|$ and define 
\begin{align*}
    \underline{\mathcal{C}}(\underline{\mathcal{I}}_i): &= \sum_{j: m-|\eta^j| > 0} (U_i, N_i, Y_i, (w^{\zeta^j}), E_i, m_j) 
    \\&= (U_i, N_i, Y_i, \mathcal{C}(\mathcal{I}_i), E_i, m_i),
\end{align*}
where 
$$
m_i: = \Pi_{j: m-|\eta^j| > 0} (m-|\eta^j|) \textrm{ and } \mathcal{C}(\mathcal{I}_i) : = (w^{\zeta^jm_i \cdot (m-|\eta^j|)^{-1}}, j :m-|\eta^j| > 0)
$$
as in \cref{definition5231}. Hence $\underline{\mathcal{C}}(\underline{\mathcal{I}}_i)$ is a $\Lambda$-marked monomial ideal by the above lemma. Moreover, $dim_{\mathbb{Z}}Y_i < dim_{\mathbb{Z}}P_i$ and by \cref{lemma5233}
$$
supp(\underline{\mathcal{C}}(\underline{\mathcal{I}}_i))  = \bigcap_{j: m-|\eta^j| > 0} supp(U_i, N_i, Y_i, (w^{\zeta^j}), E_i, m_j) \subset Y_i.
$$

\begin{observation}
    \label[observation]{observation5313} 
    Let $x \in Y_i \subset P_i$ and consider any $j$ and ideal $(w^{\zeta^j})$. By \cref{lemma5312} (4)
    $$
    \textrm{ord}_x \mathcal{O}_{P_i}  (w^{\zeta^j}) = \textrm{ord}_x\mathcal{O}_{Y_i}  (w^{\zeta^j}).
    $$
\end{observation}

\begin{remark}
    Our $\underline{\mathcal{C}}(\underline{\mathcal{I}}_i)$ construction differs slightly from \cite[Section 8 Step 1]{BM1} in that the underlying ideal is explicitly defined by our sum of marked ideals. However, as in \cref{remark5232} for the purposes of order reduction this does not matter (see \cite[Lemma 3.8]{bierstone2007functoriality}). We also only consider summing over $m-|\eta^j| > 0$ as we have defined $\Lambda$-marked ideals only with non-negative marking. We could have summed over all $j$ as in \cite[Section 8 Step 1]{BM1} but as the support of the sum of marked ideals is the intersection only the $j :m-|\eta^j| > 0$ contribute a non-trivial support. Further, it is not clear in \cite[Section 8 Step 1]{BM1} what to do when no such $m-|\eta^j| > 0$ exist, which we account for in the following lemma.
\end{remark}

\begin{lemma}
    If there exist no $j$ with $m-|\eta^j| > 0$ then 
    $$
    supp(\underline{\mathcal{I}}_i) = Y_i.
    $$ 
\end{lemma}

\textit{Proof.} Let $x \in supp(\underline{\mathcal{I}}_i) \subset P_i$. Clearly, for each $j$ with $|\alpha^j|=m$ the prime ideal of $A_i$ corresponding to $x$ must contain the variables in $x^{\alpha^j}$ and $x \in Y_i$. Let $x \in Y_i$ then for each $j$, $|\eta^j| \geq m$ and
$$
\textrm{ord}_x \mathcal{I}_i = \textrm{min}_j \{ \textrm{ord}_x z^{\eta^j}w^{\zeta^j} = |\eta^j|+\textrm{ord}_x w^{\zeta^j}\} \geq m. \quad \qed
$$
In the above situation order reduction for $\underline{\mathcal{I}}_i$ is simple. Let $Z = \cap (E_i)^j \subset U_i$ be an (it need not be unique) intersection of $E_i$ such that $Z \cap N_i = Y_i$.

\begin{lemma} 
    \label[lemma]{lemma5316} 
    Blowing up $Z$ is an order reduction for $\underline{\mathcal{I}}_i$.
\end{lemma}

\textit{Proof.} We argue as in \cref{remark5218}. Let $x' \in (P_i)'$ closed in the fiber over $(p) \not = 0$ lying over closed $x \in P_i \subset U_i$ and $x_0=p, x_1, \dots, x_n$ be a regular system of parameters of $\mathcal{O}_{P, x}$ as in  \cref{lemma5216} at $x$ with $Y_i$ corresponding to $x_1, \dots, x_r$ and consider an $x_i$ chart in $P_i'$. By definition $x_i$ will occur with nonzero power in a monomial $x^{\alpha^j}$ of $(\mathcal{O}_{P_i}\mathcal{I}_i)_{x}$ with $|\alpha^j|=m$ and hence after applying the substitution rule will correspond to a monomial $\bar{x}^{{\alpha^j}'}$ with $|{\alpha^j}'| < |\alpha^j| < m$. Hence $((\mathcal{O}_{P_i}\mathcal{I}_i)')_{x'}$ contains a monomial of degree $< m$ and so $\textrm{ord}_{x'}(\mathcal{O}_{P_i}\mathcal{I}_i)' < m$ and $supp(\underline{\mathcal{I}}_i') = \emptyset$. \qedsymbol 

\begin{remark}
    Later on we incorporate this situation into the general order reduction over each $U_i$. Note that \cite[Section 8 Step 1]{BM1} does not take this into account when $U_i \cap supp(\underline{\mathcal{I}}_i) \not = \emptyset$ and the construction of the ideal $\mathcal{C}(\mathcal{I}_i)$ does not make sense. It is simple to consider, but we include the argument for detail. 
\end{remark}

We now suppose there exists at least one $j$ with $m-|\eta^j| > 0$ and we can construct $\underline{\mathcal{C}}(\underline{\mathcal{I}}_i)$.

\begin{lemma} 
    \label[lemma]{lemma5318} 
    $$
    supp(\underline{\mathcal{C}}(\underline{\mathcal{I}}_i)) = supp(\underline{\mathcal{I}}_i) \subset Y_i.
    $$
\end{lemma}

\textit{Proof.} Let $x \in supp(\underline{\mathcal{I}}_i)$ so that $x \in P_i$. Consider the ideal $( \mathcal{O}_{P_i}\mathcal{I}_i)_x$. By the definition of the $z$-variables and the proof of \cref{lemma5104} the prime corresponding to $x$ in $B_i$ must contain all the $z$-variables and hence $x \in Y_i$. Then recalling the local description of ideals locally generated by monomials
$$
\textrm{ord}_x \mathcal{O}_{P_i}\mathcal{I}_i = \textrm{min}_j \{ \textrm{ord}_x \mathcal{O}_{P_i}(z^{\eta^j}w^{\zeta^j}) = |\eta^j|+\textrm{ord}_x \mathcal{O}_{P_i}(w^{\zeta^j})\} \geq m.
$$
By observation \ref{observation5313}, for all $j :m-|\eta^j| > 0$, $|\eta^j|+\textrm{ord}_x \mathcal{O}_{Y_i}(w^{\zeta^j}) \geq m$ and hence by \cref{corollary5106} (2),
$$
x \in supp(U_i, N_i, Y_i, (w^{\zeta^j}), m_j, E_i),
$$
which proves $supp(\underline{\mathcal{I}}_i) \subseteq supp(\underline{\mathcal{C}}(\underline{\mathcal{I}}_i))$. Conversely, let $x \in supp(\underline{\mathcal{C}}(\underline{\mathcal{I}}_i))$. Then $x \in Y_i \subset P_i$ and by similar arguments 
$\textrm{ord}_x\mathcal{O}_{Y_i}(w^{\zeta^j}) \geq m-|\eta^j| \implies |\eta^j|+\textrm{ord}_x \mathcal{O}_{P_i}(w^{\zeta^j}) \geq m$ for each $j :m-|\eta^j| > 0$. For $j: m-|\eta^j| \leq 0$, $|\eta^j| \geq m$ and so $|\eta^j|+\textrm{ord}_x \mathcal{O}_{P_i}(w^{\zeta^j}) \geq m$. Thus
$$
\textrm{ord}_x \mathcal{O}_{P_i}\mathcal{I}_i = \textrm{min}_j \{ \textrm{ord}_x \mathcal{O}_{P_i}(z^{\eta^j}w^{\zeta^j}) = |\eta^j|+\textrm{ord}_x \mathcal{O}_{P_i}(w^{\zeta^j})\}  \geq m. \quad \qed
$$

\begin{corollary}
    A blow-up centre $Z= \cap_k (E_i)^k = (\cap_k E^k) \cap U_i \subset U_i$ is permissible for $\underline{\mathcal{C}}(\underline{\mathcal{I}}_i)$ if and only if it is permissible for $\underline{\mathcal{I}}_i$.
\end{corollary}

\begin{proposition}
    Let notation be as above and let $Z=\cap_k (E_i)^k$ be a permissible centre for $\underline{\mathcal{I}}_i$ so that $Z \cap N_i$ is contained in $Y_i$ and is permissible for $\underline{\mathcal{C}}(\underline{\mathcal{I}}_i)$. Then writing $\underline{\mathcal{C}}(\underline{\mathcal{I}}_i)',\underline{\mathcal{I}}_i'$ for the transforms
    $$
    supp(\underline{\mathcal{C}}(\underline{\mathcal{I}}_i)') = supp(\underline{\mathcal{I}}_i') \subset Y_i' \subset P_i' \subset U_i'.
    $$
\end{proposition}

\textit{Proof.} $Z \subset U_i$ is defined by a collection of elements 
$$
x_1, \dots, x_l, x_{l+1}, \dots, x_s, x_{s+1}, \dots,  x_r \in A_i
$$ such that by permissibility we can assume that $P_i \cong \textrm{Spec} (A_i/(x_{1}, \dots, x_l)) \cap N_i$ and $Y_i \cong \textrm{Spec} (A_i/(x_1, \dots, x_s)) \cap N_i$ so that the $z$-variables are $x_{l+1}, \dots, x_s$. As in the $\underline{\mathcal{C}}(\underline{\mathcal{I}}_i)$ construction write the $x_{l+1}, \dots, x_s$ as $z$-variables and write the $x_{s+1}, \dots, x_r$ as $v$-variables, i.e. those coming from the blow-up but not containing $Y_i$. Write the remaining variables in the mononomials generating $I_i$ as $u$-variables so each monomial can be written as
$$
x^{\alpha^j}=z^{\eta^j}w^{\zeta^j}=z^{\eta^j}v^{\beta^j}u^{\gamma^j} \in A_i, \quad w^{\zeta^j}=v^{\beta^j}u^{\gamma^j}
$$
and note that none of the monomials have non-zero degree in $x_{1}, \dots, x_l$. We will use an affine blow-up chart of $U_i'$ given by one of the variables of $Z$ to compare the supports. Only those charts corresponding to $v$-variables intersect $Y_i'$ and both supports are contained in $Y_i'$ so it is sufficient to consider a $v_k$-chart $(U_i')_k$. We will compute the transforms of the monomials generating $\underline{\mathcal{I}}_i$ and $\underline{\mathcal{C}}(\underline{\mathcal{I}}_i)$ in $(U_i')_k$ using the substitution rule (see \cref{lemma4403} and the proof of \cref{lemma3315}) and compare.
\\The transforms of monomials of $\underline{\mathcal{I}}_i$ in $(U_i')_k$: 
\begin{align*} 
    z^{\eta^j} \textrm{ with } \space|\eta^j|=m \mapsto \bar{z}^{\eta^j}, \textrm{ with } |\eta^j|=m.
    \\z^{\eta^j}v^{\beta^j}u^{{\gamma^j}} \mapsto
    \begin{cases} 
      \bar{v_k}^{|\eta^j|+|\beta^j|-m}\bar{z}^{\eta^j}\bar{v}^{\beta^j}\bar{u}^{{\gamma^j}} & \textrm{if } (\beta^j)_k = 0  \\
      \bar{v_k}^{|\eta^j|+|\beta^j|-m}\bar{z}^{\eta^j}\bar{v}^{{\beta^j}'}\bar{u}^{{\gamma^j}} &  \textrm{if } (\beta^j)_k \not = 0 
    \end{cases}.
\end{align*}
By \cref{lemma52316} the local transform of $\underline{\mathcal{C}}(\underline{\mathcal{I}}_i)$ in $(U_i')_k$ will be the sum of the transforms of $w^{\zeta^j}, \space j :m-|\eta^j| > 0$, which will be as follows: 
\begin{align*}
    w^{\zeta^j}=v^{\beta^j}u^{{\gamma^j}} \mapsto
    \begin{cases} 
      \bar{v_k}^{|\beta^j|-(m-|\eta^j|)}\bar{v}^{\beta^j}\bar{u}^{{\gamma^j}} & \textrm{if } (\beta^j)_k = 0  \\
      \bar{v_k}^{|\beta^j|-(m-|\eta^j|)}\bar{v}^{{\beta^j}'}\bar{u}^{{\gamma^j}} &  \textrm{if } (\beta^j)_k \not = 0 
    \end{cases},
\end{align*}
where $|{\beta^j}'|=|\beta^j|-|(\beta^j)_k|$. For $j :|\eta^j| \geq m$ the transform will be a monomial with factor $\bar{z}^{\eta^j}, |\eta^j| \geq m$. To compare orders it is also useful to consider the corresponding local picture.
\\As usual let $x' \in Y_i' \subset P_i'$ be a point in $(U_i')_k$ closed in a  $p$-fiber lying over $x \in Y_i \subset P_i$ that is contained in $Z$. Then we can choose a regular system of parameters 
$$
x_0=p, x_1, \dots, x_l, x_{l+1}, \dots, x_s, x_{s+1}, \dots,  x_r,x_{r+1}, \dots, x_c, \dots, x_n
$$
for $\mathcal{O}_{U_i, x}$ corresponding to the elements of $A_i$ above such that $ \mathcal{O}_{P_i, x} \cong \mathcal{O}_{N_i, x}/(x_{1}, \dots, x_l)$, $\mathcal{O}_{Y_i, x} \cong \mathcal{O}_{P_i, x}/(x_{l+1}, \dots, x_s)$ and $Z \cap P$ is defined in $\mathcal{O}_{P_i, x}$ by $(x_{l+1}, \dots, x_s, x_{s+1}, \dots, x_r)$. As the prime ideal  of $A_i$ corresponding to $x$ contains all $z,v$-variables but not necessarily all $u$-variables ,$(\mathcal{O}_{P_i}\mathcal{I}_i)_x$ is generated by the monomials $z^{\eta^j}v^{\beta^j}u^{{\gamma^j}'}$, where $u^{{\gamma^j}'}$ is the monomial given by $u^{\gamma^j}$ with any $u$-variables not contained in $x$ divided off. Similarly, $(\mathcal{O}_{P_i}(w^{\zeta^j}))_x$ is generated by the monomial $v^{\beta^j}u^{{\gamma^j}'} \subset \mathcal{O}_{P_i, x}$ and $(\mathcal{O}_{Y_i}(w^{\zeta^j}))_x$ is generated by the same monomial $v^{\beta^j}u^{{\gamma^j}} \subset \mathcal{O}_{Y_i, x}$. Then $(\mathcal{O}_{P_i'}\mathcal{I}_i')_{x'}$, $(\mathcal{O}_{P_i'}(w^{\zeta^j}))_{x'}$ and $(\mathcal{O}_{Y_i'}(w^{\zeta^j})')_{x'}$ are generated by the transform of the monomials as above but with ${\gamma^j}'$ instead of $\gamma^j$. 

In order to use these local transforms to compare supports we need 
$$
supp(\underline{\mathcal{I}}_i') \cap (U_i')_k \not = \emptyset \iff supp(\underline{\mathcal{C}}(\underline{\mathcal{I}}_i)') \cap (U_i')_k \not = \emptyset.
$$
Recall by \cref{lemma5235}
$$
supp(\underline{\mathcal{C}}(\underline{\mathcal{I}}_i)') = \bigcap_{j: m-|\eta^j| > 0} supp((U_i, N_i, Y_i, (w^{\zeta^j}), m-|\eta^j|, E_i)').
$$
Let $x',x$ be as above and suppose $x' \in supp(\underline{\mathcal{I}}_i') \cap (U_i')_k$ so that $x$ is contained in $supp(\underline{\mathcal{C}}(\underline{\mathcal{I}}_i)) = supp(\underline{\mathcal{I}}_i)$. We wish to show that for any $j: m-|\eta^j| > 0,$
$$
x' \in  supp((U_i, N_i, Y_i, (w^{\zeta^j}), E_i, m-|\eta^j|)').
$$
First suppose $(\beta^j)_k = 0$ for this $j$. Then 
\begin{align}
    &\textrm{ord}_{x'}\mathcal{O}_{P_i'} \nonumber (\bar{v_k}^{|\eta^j|+|\beta^j|-m}\bar{z}^{\eta^j}\bar{v}^{\beta^j}\bar{u}^{\gamma^j}) = |\eta^j|+\textrm{ord}_{x'}\mathcal{O}_{P_i'}(\bar{v_k}^{|\eta^j|+|\beta^j|-m}\bar{v}^{\beta^j}\bar{u}^{\gamma^j})\geq m 
    \\&\iff \textrm{ord}_{x'}\mathcal{O}_{Y_i'} \nonumber (\bar{v_k}^{|\eta^j|+|\beta^j|-m}\bar{v}^{\beta^j}\bar{u}^{\gamma^j}) \geq m-|\eta^j| 
    \\&\iff x' \in supp((U_i, N_i, Y_i, (w^{\zeta^j}), E_i, m-|\eta^j|)').   
    \label{eqn4.1}
\end{align}
If $(\beta^j)_k = 0$ then
\begin{align}
    &\textrm{ord}_{x'}\mathcal{O}_{P_i'}\bar{v_k}^{|\eta^j|+|\beta^j|-m}\bar{z}^{\eta^j}\bar{v}^{{\beta^j}'}\bar{u}^{\gamma^j} = |\eta^j|+\textrm{ord}_{x'}\mathcal{O}_{P_i'}(\bar{v_k}^{|\eta^j|+|\beta^j|-m}\bar{v}^{{\beta^j}'}\bar{u}^{\gamma^j})\geq m \nonumber 
    \\&\iff \textrm{ord}_{x'}\mathcal{O}_{Y_i'}(\bar{v_k}^{|\eta^j|+|\beta^j|-m}\bar{v}^{{\beta^j}'}\bar{u}^{\gamma^j}) \geq m-|\eta^j| \nonumber 
    \\&\iff x' \in supp((U_i, N_i, Y_i, (w^{\zeta^j}), E_i, m-|\eta^j|)').   
    \label{eqn4.2}
\end{align}
Conversely, let $x' \in supp((U_i, N_i, Y_i, (w^{\zeta^j}), E_i, m-|\eta^j|)')$ for all $j: m-|\eta^j| > 0$. We wish to show the order at $x'$ is $\geq m$ for each transform of a mononomial in $\mathcal{I}_i$. But by our computations it suffices to check the order of $\bar{v_k}^{|\eta^j|+|\beta^j|-m}\bar{z}^{\eta^j}\bar{v}^{{\beta^j}'}\bar{u}^{\gamma^j}$ at $x'$ is $\geq m$, which holds by (\ref{eqn4.1}) and (\ref{eqn4.2}) above. Thus to show the supports agree it suffices to consider the $(U_i')_k $ that intersect either and hence both supports.
\\Now when $supp(\underline{\mathcal{I}}_i') \cap (U_i')_k \not = \emptyset$, $\underline{\mathcal{I}}_i' \cap (U_i')_k$ will be of maximal order (by \cref{lemma4206}) and we can apply the $\underline{\mathcal{C}}$ construction. By the arguments of \cref{lemma5318} we have
$$
supp(\underline{\mathcal{I}}_i' \cap (U_i')_k) = supp(\underline{\mathcal{C}}(\underline{\mathcal{I}}_i' \cap (U_i')_k))
$$
so it suffices to show
\begin{align}
    supp(\underline{\mathcal{C}}(\underline{\mathcal{I}}_i' \cap (U_i')_k))  &= \bigcap_{j: m-|\eta^j| > 0} supp((U_i, N_i, Y_i, (w^{\zeta^j}), E_i, m-|\eta^j|)') \cap (U_i')_k
    \\&= supp(\underline{\mathcal{C}}(\underline{\mathcal{I}}_i)' \cap (U_i')_k).
\end{align} 

Now in constructing $\underline{\mathcal{C}}(\underline{\mathcal{I}}_i' \cap (U_i')_k)$ we need to consider which monomials in $\underline{\mathcal{I}}_i' \cap (U_i')_k$ have order $m$ to construct the "new" $z$-variables. As $\bar{z}^{\eta^j}$ with $|\eta^j|=m$ will be monomials of $\underline{\mathcal{I}}_i' \cap (U_i')_k$ the $\bar{z}$-variables will be a subset of these "new" $z$-variables. Denote by $\bar{t}$ the remaining "new" $z$-variables. Then each monomial can be written as $\bar{z}^{\eta^j}\bar{t}^{\xi^j}\bar{w}^{\delta^j}$, where $\bar{t}^{\xi^j}\bar{w}^{\delta^j}$ will come from the transform of $w^{\zeta^j}$ and we use the $\bar{w}^{\delta^j}$ with $|\xi^j|+|\eta^j| < m$ to construct $\underline{\mathcal{C}}(\underline{\mathcal{I}}_i' \cap (U_i')_k)$. It follows that for $j:|\eta^j| \geq m $, $\bar{z}^{\eta^j}$ will contribute degree $\geq m$ and so is only a monomial of degree $m$ if all other variables have $0$-power, i.e. $\xi^j = \delta^j = \bar{0}$. Hence, the new monomials of degree $m$ (those that contribute $\bar{t}$-variables) can only come from a $j: m-|\eta^j| > 0$ and be of the form
$$
\bar{v_k}^{|\eta^j|+|\beta^j|-m}\bar{z}^{\eta^j}\bar{v}^{{\beta^j}'}\bar{u}^{\gamma^j} = \bar{z}^{\eta^j}(w^{\zeta^j})'.
$$
If there are no new monomials of degree $m$ then there are no $\bar{t}$, $\delta^j=\zeta^j$ and $\underline{\mathcal{C}}(\underline{\mathcal{I}}_i' \cap (U_i')_k) = \underline{\mathcal{C}}(\underline{\mathcal{I}}_i)' \cap (U_i')_k$ and we are done. Note that if the order of any of the above monomials is $< m$ then $supp(\underline{\mathcal{I}}_i') \cap (U_i')_k = \emptyset$ so we may ignore this. Assume there is at least one new monomial of degree $m$. Let $j: m-|\eta^j| > 0$ correspond to this monomial. Then the degree of $(w^{\zeta^j})'$ is $|\beta^j|+|\gamma^j|+|{\beta^j}'|+|\eta^j|-m = m-|\eta^j|$ and it follows that $supp((U_i, N_i, Y_i, (w^{\zeta^j}), m-|\eta^j|, E_i)') \cap (U_i')_k$ is contained in the vanishing of the variables in $(w^{\zeta^j})'$ (it already contains all the $\bar{z}$-variables from the definition of $Y_i$). Thus
$$
 \bigcap_{j: m-|\eta^j| > 0} supp((U_i, N_i, Y_i, (w^{\zeta^j}), E_i, m-|\eta^j|)') \cap (U_i')_k \subset (Y_i)_k,
$$
where $(Y_i)_k$ is the closed subscheme of $(U_i')_k$ defined by the vanishing of the $\bar{z},\bar{t}$ variables that gives the closed subscheme of $(P_i)' \cap (U_i')_k$ in $\underline{\mathcal{C}}(\underline{\mathcal{I}}_i' \cap (U_i')_k)$. Let $x' \in  \bigcap_{j: m-|\eta^j| > 0} supp((U_i, N_i, Y_i, (w^{\zeta^j}), m-|\eta^j|, E_i)') \cap (U_i')_k$ and consider a $j: m-|\eta^j| > 0$ such that the monomial $\bar{z}^{\eta^j}(w^{\zeta^j})'$ has not become degree $m$ and $m-(|\eta^j|+|\xi^j|) > 0$, where $(w^{\zeta^j})' =  \bar{t}^{\xi^j}\bar{w}^{\delta^j}$. Then noting that $\textrm{ord}_{x'}\mathcal{O}_{Y_i'}(\bar{w}^{\delta^j}) = \textrm{ord}_{x'}\mathcal{O}_{(Y_i)_k}(\bar{w}^{\delta^j})$ as in observation \ref{observation5313},
\begin{align}
\label{eqn4.5}
 &\textrm{ord}_{x'} \mathcal{O}_{Y_i'}(w^{\zeta^j})' = \textrm{ord}_{x'} \mathcal{O}_{Y_i'}\bar{t}^{\xi^j} \bar{w}^{\delta^j} = |\xi^j|+\textrm{ord}_{x'}\mathcal{O}_{Y_i'} \bar{w}^{\delta^j} \geq m-|\eta^j| \nonumber
 \\&\iff \textrm{ord}_{x'} \mathcal{O}_{(Y_i)_k} \bar{w}^{\delta^j} \geq m-(|\eta^j|+|\xi^j|) \nonumber
 \\&\iff x' \in supp((U_i')_k, (N_i') \cap (U_i')_k, (Y_i)_k, \bar{w}^{\delta^j}, (E_i)'|_{(U_i')_k}, m-(|\eta^j|+|\xi^j|)),   
\end{align}
which implies
\begin{align*}
  &\bigcap_{j: m-|\eta^j| > 0} supp((U_i, N_i, Y_i, (w^{\zeta^j}), E_i, m-|\eta^j|)') \cap (U_i')_k 
  \\&\subset \bigcap_{j: m-(|\eta^j|+|\xi^j|) > 0} supp((U_i')_k, (N_i')_k, (Y_i)_k, \bar{w}^{\delta^j}, (E_i)'|_{(U_i')_k}, m-(|\eta^j|+|\xi^j|)) 
  \\&= supp(\underline{\mathcal{C}}(\underline{\mathcal{I}}_i' \cap (U_i')_k)).
\end{align*}
Conversely, let $x' \in supp(\underline{\mathcal{C}}(\underline{\mathcal{I}}_i' \cap (U_i')_k)) \subset (Y_i)_k$ and consider a $j:m-|\eta^j| > 0$. If $(w^{\zeta^j})'$ features in a monomial of degree $m$ then as argued before 
\\$supp((U_i, N_i, Y_i, (w^{\zeta^j}), E_i, m-|\eta^j|)') \cap (U_i')_k $ will be given by the vanishing of the $\bar{t}$-variables contained in $(w^{\zeta^j})'$ and hence $x' \in (Y_i)_k$ implies 
\\$x' \in supp((U_i, N_i, Y_i, (w^{\zeta^j}), E_i, m-|\eta^j|)') \cap (U_i')_k $. Now suppose $(w^{\zeta^j})'$ does not feature in a degree $m$ monomial but $m \leq (|\eta^j|+|\xi^j|)$. Then $|\xi^j| \geq m-|\eta^j|$ so $\textrm{ord}_{x'}\mathcal{O}_{Y_i'}(w^{\zeta^j})' = |\xi^j|+\textrm{ord}_{x'}\mathcal{O}_{(Y_i)_k}\bar{w}^{\delta^j} \geq m-|\eta^j|$ and 
\\$x' \in supp((U_i, N_i, Y_i, (w^{\zeta^j}), E_i, m-|\eta^j|)') \cap (U_i')_k $. Finally, if $m > (|\eta^j|+|\xi^j|)$ then we can use (\ref{eqn4.5}) above, which shows the reverse inclusion. \qedsymbol
\\
\\The argument for the previous proposition is long and tedious but we cannot think of a more elegant alternative. It is not, as one might expect, necessarily true that $\underline{\mathcal{C}}(\underline{\mathcal{I}}_i' \cap (U_i')_k) = \underline{\mathcal{C}}(\underline{\mathcal{I}}_i)' \cap (U_i')_k$ as illustrated by the following example.

\begin{example}
    Let $U=\textrm{Spec} (\mathbb{Z}[x, y, u ,v])$, $I = (x^2y^3, x^2v^6, y^4u^5) \subset \mathbb{Z}[x, y, u ,v]$, $(E^1, E^2, E^3, E^4)=(x, y, u, v)$ and $E=\sum_i E^i$. Then $(U, U, U, I, E, 5)$ is a $\Lambda$-marked monomial ideal of maximal order. The $z$-variables are $x, y$ and $w^{\zeta^1}=v^6,w^{\zeta^2}=u^5, m-|\eta^1|=3, m-|\eta^2|=1$. Consider the blow-up given by $Z=V(x, y, u, v)$. This is permissible and we can consider $U_{v}$ in $B_Z X$ given by the $v$-chart. Using the substitution rule and the coordinates $\bar{x}, \bar{y}, \bar{v}, \bar{u}$ we get that 
    $$
    I' \cap  U_{v}= (\bar{x}^2\bar{y}^3, \bar{x}^2\bar{v}^{=3}, \bar{y}^4\bar{u}^5\bar{v}^{4}).
    $$
    Applying the $\underline{\mathcal{C}}$-construction to $I' \cap  U_{v}$ the new $\bar{z}$-variables are $\bar{x}, \bar{y}, \bar{v}$ as the monomial $\bar{x}^2\bar{v}^3$ has order $5$ and we only have $({w_1}^{\zeta_1})=(\bar{u}^5)$. However, the transform $\underline{\mathcal{C}}(I)'$ will be given by the sum of the transforms of the $w^{\zeta^1}=v^6,w^{\zeta^2}=u^5$, which are $\bar{v}^{3}, \bar{u}^5\bar{v}^{4}$. This does not equal $(\bar{u}^5)$ and hence the marked ideals are not equal on $U_{v}$.
\end{example}

\begin{remark}
    The example above doesn't rely on the coefficients being in $\mathbb{Z}$ so we similarly get a counter-example to the claim in \cite[Section 8 Step 1]{BM1} that $\underline{\mathcal{C}}(\underline{\mathcal{H}}_{\sigma})' = \underline{\mathcal{C}}(\underline{\mathcal{H}}_{\sigma}')$. They do, however, have the same supports so this is not an issue for the purpose of order reduction.
\end{remark}

\begin{lemma}
    Any blow-up sequence for $U_i$ is permissible for $\underline{\mathcal{I}}_i$ if and only if it is permissible for $\underline{\mathcal{C}}(\underline{\mathcal{I}}_i)$.
\end{lemma}

\textit{Proof.} Follows by applying induction and using the arguments of the above proposition on the finitely many affine blow-up charts over an affine chart. \qedsymbol
\\
\\We must also consider what happens when we take the transforms when $\mathcal{C}(\mathcal{I}_i)$ does not make sense.

\begin{lemma}
\label[lemma]{lemma4.3.14} 
    Let $\underline{\mathcal{I}}_i$ be such that there exists no $j$ with $m-|\eta^j| > 0$. Then after a permissible blow-up sequence for $\underline{\mathcal{I}}_i$ either $supp\underline{\mathcal{I}}_i' = \emptyset$ or $supp\underline{\mathcal{I}}_i' =Y_i'$. Further, an order reduction of $\underline{\mathcal{I}}_i'$ can be achieved by blowing up a permissible $Z$ such that $ Z \cap N_i' = Y_i'$.
\end{lemma}

\textit{Proof.} We first show the lemma holds for a single blow-up. Let $Z=\cap E^k \cap U_i$ be a permissible blow-up for $\underline{\mathcal{I}}_i$ so either $Z \cap N_i =Y_i$ and we produce an order reduction by \cref{lemma5316} or $Z \cap N_i \subset Y_i$. Let notation be as in the proof of the above proposition. Now $supp\underline{\mathcal{I}}_i'  \subset Y_i'$ and locally at a closed point $x$ each monomial of $(\mathcal{O}_{P_i}\mathcal{I}_i)_x$ contains a factor $z^{\eta^j}$ with $|\eta^j| \geq m$. It follows that for a closed point $x'$ lying over $x$ in $(U_i')_k$ each monomial of $(\mathcal{O}_{P_i'}\mathcal{I}_i')_{x'}$ contains a factor $ \bar{z}^{\eta^j}$ with $|\eta^j| \geq m$ so $(Y_i)' \cap (U_i)_k \subset supp\underline{\mathcal{I}}_i'$ as all $\bar{z}$-variables vanish on $x' \in (Y_i)'$. For further blow-ups we can argue inductively over affine charts. Lastly, blowing up a permissible $Z$ such that $ Z \cap N_i' = Y_i'$ over an affine chart is an order reduction by the arguments of \cref{lemma5316} on affine covers. \qedsymbol

\begin{lemma}
    A permissible blow-up sequence for $\underline{\mathcal{I}}$ restricts to a permissible blow-up sequence for each $\underline{\mathcal{I}}_i$ and a permissible blow-up sequence for $\underline{\mathcal{I}}_i$ is the restriction over $U_i$ of a permissible blow-up sequence for $\underline{\mathcal{I}}$.

\end{lemma}

\textit{Proof.} We first show this for a single blow-up. Let $Z= \bigcap_{j \in J_0} (E^j)$ be such that $Z \cap U_i=\bigcap_{j \in J_0} (E^j \cap U_i) \not = \emptyset$. Note that every permissible centre for $\underline{\mathcal{I}}_i$ is of the form $Z \cap U_{i}$ for some such $Z$. Thus, it suffices to show
$$
Z \cap N_i  \subset supp(\underline{\mathcal{I}}_i) \iff Z \cap N \subset supp(\underline{\mathcal{I}}).
$$
If $Z \cap N \subset supp(\underline{\mathcal{I}})$ then clearly $Z \cap N \cap U_i = (Z \cap U_i) \cap N_i \subset supp(\underline{\mathcal{I}}) \cap U_i = supp(\underline{\mathcal{I}}_i)$. Conversely, suppose
$Z \cap N_i  \subset supp(\underline{\mathcal{I}}_i) \subset P_i$. Recall by \cref{lemma4503} there exists $J_1 \subset J_0, J_P \subset J_N$ such that $Z\cap N= \bigcap_{j \in J_1}(E^j)|_N, P=\bigcap_{j \in J_P} (E^j)|_N$. It follows that $J_P \subset J_1$ and so $Z \cap N \subset P$. Further, for each 
$$
\mathcal{O}_N{\mathcal{I}}_j = \mathcal{O}_X(-\sum_{i_j} a_{i_j} E^{i_j}|_N)
$$
at least one of the $E^{i_j}|_N$ is equal to a $E^j|_N, j \in J_1$ and $
\sum_{i_j, i_j 
\in J_1} a_{i_j} \geq m$ as this will be the degree of each "monomial" in $(\mathcal{O}_{N}\mathcal{I})(N_i)$ given by the image of the monomials in $\mathcal{I}(U_i).$ Thus for $x \in Z$ the monomials generating $(\mathcal{O}_P\mathcal{I})_x$ will have degree $\geq m$ and $x \in supp(\underline{\mathcal{I}})$. For further blow-ups we can argue inductively on finitely many affine opens. \qedsymbol
\\
\\\textbf{Applying induction.} We now have the relevant conditions to use the inductive hypothesis that order reductions exist for $\Lambda$-marked monomial ideals of $dim_{\mathbb{Z}}P < n$. We can apply the arguments of \cite[Proof of Theorem 8.5 Step 1]{BM1} using the analogous conditions developed above. Note we must be careful to also consider when there are no $j: m-|\eta^j| > 0$, which can be handled by \cref{lemma4.3.14}.
\\


\paragraph{\textbf{Monomial derivatives and $\underline{\mathcal{C}}(\underline{\mathcal{I}})$.}}
Our argument used the "by hand" construction of Bierstone and Milman (see \cite[Section 8.1 step 1]{BM1}) as opposed to the use of derivatives of ideals as in Kollar and other resolutions of singularities in characteristic $0$. The obstruction to using derivatives in characteristic $\not = 0$ and over $\mathbb{Z}$, is that usual derivatives produce monomials with coefficients that are not necessarily units, which leads to the failure of certain properties  (cf. \cite[Lemma 3.74]{10.2307/j.ctt7rptq}, \cite[Lemma 3.2]{bierstone2007functoriality}). For the ideals we consider, monomial derivatives can be used instead of usual derivatives and the "by hand" construction is recovered.

Let $X$ be a $\Lambda$-scheme smooth over $\mathbb{Z}$ with an affine $\Lambda$-cover and $E=\sum E^i$ a $\mathbb{Z}$-flat locally toric snc divisor. Let $\mathcal{I}=\sum \mathcal{I}_j$ be an ideal locally generated by monomials for $E$ and $r \geq 0$. We can define the ideal $\mathcal{D}^r_{mon}\mathcal{I}$, the $r$-th monomial derivative of $\mathcal{I}$, as the ideal generated by all $\beta$-th formal derivatives of $\mathcal{I}$ for $|\beta| \leq r$, as in \cite[lemma 7]{DifferentialOperatorsAlgorithmicWeightedResolution} but without coefficients. 
Let $m=\textrm{max-ord }\mathcal{I}$ and $\textrm{Spec} \, A_i=U_i$ be an affine open in the $\Lambda$-cover such that $supp(\mathcal{I}, m) \cap U_i \not = \emptyset$ and $\mathcal{I}(U_i)$ is generated by monomials in the elements $x_j \in A_i$, which can be split into $z$-variables and $w$-variables as in the $\underline{\mathcal{C}}$ construction. Then the $z$-variables that define the closed subscheme $Y_i \subset P_i$ are elements of $\mathcal{D}_{mon}^{\leq m-1}(\mathcal{O}_P\mathcal{I})$ such that for each $z_k $, $\mathcal{D}_{mon}^{\leq 1}(z_k) = (1)$, as in characteristic $0$ (See \cite[Corollary 4.2]{bierstone2007functoriality}). Moreover, $U_i$ is an affine open on which these $z_k$ define the smooth $\Lambda$-subscheme $Y_i$, i.e. monomial derivatives produce smooth hyperspaces of "maximal contact" compatible with $\Lambda$-structure. Consider the marked ideal $(\mathcal{I}, m)$ with $\textrm{max-ord }\mathcal{I}=m$. We define the monomial derivative of a marked ideal as in \cite[lemma 3.74]{10.2307/j.ctt7rptq}, \cite[Lemma 3.2]{bierstone2007functoriality} and define $(\mathcal{C}_{mon}(\mathcal{I}), (m-1)!)=\sum_{r=1}^{m-1}\mathcal{D}_{mon}^{\leq r}(\mathcal{I}, m)$, the \textit{monomial coefficient ideal} of $\mathcal{I}$ (cf. \cite[Section 3.4]{bierstone2007functoriality}). Then
$$
supp\underline{\mathcal{C}}(\underline{\mathcal{I}_i})=supp (\mathcal{C}_{mon}(\mathcal{O}_P\mathcal{I})|_{Y_i}, (m-1)!),
$$
where $\underline{\mathcal{C}}(\underline{\mathcal{I}_i})$ is as in 
\cref{C:5_Sec:Analysis:FM}. It follows that in the $\underline{\mathcal{C}}$ construction we are actually blowing up in the support of the local monomial coefficient ideal of $\mathcal{O}_P\mathcal{I}$ as is done in characteristic $0$. 
Using monomial derivatives in place of usual derivatives and following characteristic $0$ methods yields the same initial set up as the "by hand" construction but the way we argue for order reduction remains different. In characteristic $0$, rather than working directly with local equations of $\mathcal{I}$, order reduction is argued by certain properties of higher derivatives under birational transforms. An obstruction for this method in characteristic $p$ or mixed characteristic is again the failure of  \cite[Lemma 3.74]{10.2307/j.ctt7rptq}, \cite[Lemma 3.2]{bierstone2007functoriality}. As this property holds for monomial derivatives in our setting it seems likely that the same arguments as in characteristic $0$ could be used. However, the arguments would probably boil down to similar ones on monomial indices used in the "by hand" construction and so may not be very enlightening. Hence, we do not detail an analogue of the characteristic $0$ method though $\Lambda$-equivariance of this method may be worth investigating.


\subsubsection{Reduction to the monomial part}
\label{C5_S3_S2} We will show that assuming the existence of order reduction in maximal order we can reduce order reduction to the monomial case following the methods of Kollar \cite[Proof of 3.107 step 1]{10.2307/j.ctt7rptq} and Bierstone and Milman \cite[Step 2 of Section 8]{BM1}, \cite[Step II general case]{bierstone2007functoriality}. 

Let $\underline{\mathcal{I}} := (X, N, P, \mathcal{I}, E, m)$ be a $\Lambda$-marked monomial ideal so $\mathcal{I}$ is an ideal locally generated by monomials for $E$. We can collect all $E^k$ that divide every $\mathcal{I}_j$ into a single monomial ideal, the \textit{monomial part} of $\mathcal{I}$,
$$
\mathcal{M}(\mathcal{I}):=\sum_{k: \forall j \exists i_j=k } \mathcal{O}_X(-\sum_{k} a_{k} E^{k}),
$$
where $a_k=\textrm{min}_j \{a_{i_j}: i_j=k  \}$, i.e. we take out all common factors of $\mathcal{O}_{X}(-E^k)$. Define $\mathcal{N}(\mathcal{I})=\sum_j \mathcal{I}_j'$, the \textit{non-monomial} part, where $\mathcal{I}_j'=\mathcal{I}_j \cdot \mathcal{M}(\mathcal{I})^{-1}$ and no $E^i$ divide all $\mathcal{I}_j'$. Then
$$
\mathcal{I}=\mathcal{M}(\mathcal{I}) \cdot \mathcal{N}(\mathcal{I}) 
$$
and both $\mathcal{M}(\mathcal{I}),\mathcal{N}(I)$ are locally generated by monomials for $E$ and cannot further be factored into a monomial and non-monomial part. Note that $\mathcal{M}(\mathcal{I})$ is locally generated by a single monomial for $E$.

\begin{observation}[Fibers]
    For all prime ideals $(p) \subset \mathbb{Z}$ recall we have the fibers $\underline{\mathcal{I}}_p$. We can consider the same factorization of $\mathcal{I}_p \subset \mathcal{O}_{X_p}$ using $E_p$ and produce $\mathcal{M}(\mathcal{I}_p), \mathcal{N}(\mathcal{I}_p)$. It is clear that
    $$
    \mathcal{M}(\mathcal{I}_p)=(\mathcal{M}(\mathcal{I}))_p \textrm{ and } \mathcal{N}(\mathcal{I}_p)=(\mathcal{N}(\mathcal{I}))_p.
    $$
\end{observation}

Now that we have the necessary facts we proceed by applying the theorem for order reduction in maximal order to make $supp(\underline{\mathcal{I}})$ disjoint from $\mathcal{N}(I)$. We will use the companion of ideal of Bierstone and Milman's \cite{BM1}, \cite{bierstone2007functoriality}, though the companion ideal of Kollar would also work.
\\
\\Let $\underline{\mathcal{I}}$ be a $\Lambda$-marked monomial ideal with $\mathcal{I}=\mathcal{M}(\mathcal{I}) \cdot \mathcal{N}(\mathcal{I})$. Define
$$
\nu_{\underline{\mathcal{I}}}:= max\{\textrm{ord}_x(\mathcal{N}(\mathcal{I})): \space x \in supp(\underline{\mathcal{I}})\}.
$$
Define the $\Lambda$-marked monomial ideal
$$
\underline{\mathcal{N}}(\underline{\mathcal{I}}):=(X, N, P,  \mathcal{N}(\mathcal{I}),  E, \nu_{\underline{\mathcal{I}}}).
$$

\textbf{Reduction to $\nu_{\underline{\mathcal{I}}} < m$.} Assume that $\nu_{\underline{\mathcal{I}}} \geq m$. For any $x \in P$, $\textrm{ord}_x\mathcal{O}_P\mathcal{I}  \geq \textrm{ord}_x\mathcal{O}_P\mathcal{N}(\mathcal{I}) $ so $\nu_{\underline{\mathcal{I}}} = \textrm{max-ord }\mathcal{O}_P\mathcal{N}(\mathcal{I})$ and $\underline{\mathcal{N}}(\underline{\mathcal{I}})$ is of maximal order. Thus we may apply order reduction to $\underline{\mathcal{N}}(\underline{\mathcal{I}})$. We apply order reduction to $\underline{\mathcal{N}}(\underline{\mathcal{I}})$ as in the first step of Kollar \cite[proof of 3.107]{10.2307/j.ctt7rptq} until the maximal order of $\underline{\mathcal{N}}(\underline{\mathcal{I}})'$ is less than $m$, i.e. we have a permissible sequence for $\underline{\mathcal{I}}$ such that 
$$
\mathcal{I}' = \mathcal{M}(\mathcal{I}') \cdot \mathcal{N}(\mathcal{I}')
$$
with $\nu_{\underline{\mathcal{I}}'} < m$. Arguments for permissibility are similar to local blow-up coordinate arguments in \cref{C:5_Sec:Analysis:FM}. For convenience write $\underline{\mathcal{I}}$ for $\underline{\mathcal{I}}'$ and so we can assume $\nu_{\underline{\mathcal{I}}} < m$. We assume $\nu_{\underline{\mathcal{I}}} \not= 0$ and show that we can reduce to the $0$ case. Otherwise we can proceed straight to the $\nu_{\underline{\mathcal{I}}} = 0$ step below.
\\
\\\textbf{Reduction to $\nu_{\underline{\mathcal{I}}}=0$.} Define the two $\Lambda$-marked monomial ideals
\begin{align*}
    &\underline{\mathcal{N}}(\underline{\mathcal{I}}):=(X, N, P,  \mathcal{N}(\mathcal{I}), E, \nu_{\underline{\mathcal{I}}})
    \\&\underline{\mathcal{M}}(\underline{\mathcal{I}}):=(X, N, P,  \mathcal{M}(\mathcal{I}), E, m-\nu_{\underline{\mathcal{I}}})
\end{align*}
and the companion ideal
$$
\underline{\mathcal{G}}(\underline{\mathcal{I}}):= \underline{\mathcal{N}}(\underline{\mathcal{I}})+\underline{\mathcal{M}}(\underline{\mathcal{I}})=(X, N, P,  \mathcal{M}(\mathcal{I})^{\nu_{\underline{\mathcal{I}}}}+\mathcal{N}(\mathcal{I})^{m-\nu_{\underline{\mathcal{I}}}}, E, \nu_{\underline{\mathcal{I}}}(m-\nu_{\underline{\mathcal{I}}})).
$$
For convenience let $\mathcal{G}(\mathcal{I}): = \mathcal{M}(\mathcal{I})^{\nu_{\underline{\mathcal{I}}}}+\mathcal{N}(\mathcal{I})^{m-\nu_{\underline{\mathcal{I}}}}$.
It follows from \cref{lemma4404} that for $x \in X$,
$$
\textrm{ord}_x\mathcal{G}(\mathcal{I})= \textrm{min}\{(m-\nu_{\underline{\mathcal{I}}})\textrm{ord}_x\mathcal{N}(\mathcal{I}), {\nu_{\underline{\mathcal{I}}}}\textrm{ord}_x\mathcal{M}(\mathcal{I})\}.
$$
Then by \cref{lemma5104}, for $x \in P$
$$
\textrm{ord}_x\mathcal{O}_P\mathcal{G}(\mathcal{I}) = \textrm{min}\{(m-\nu_{\underline{\mathcal{I}}})\textrm{ord}_x(\mathcal{O}_P\mathcal{N}(\mathcal{I})), {\nu_{\underline{\mathcal{I}}}}\textrm{ord}_x(\mathcal{O}_P\mathcal{M}(\mathcal{I}))\}
$$
and $\underline{\mathcal{G}}(\underline{\mathcal{I}})$ is of maximal order. We can then apply order reduction to $\underline{\mathcal{G}}(\underline{\mathcal{I}})$ as in  \cite[Proof of Theorem 8.5 Step 2 (a)]{BM1} using properties of $\Lambda$-marked monomial ideals developed in \cref{C5_1} and local blow-up coordinates to argue for permissibility as in \cref{C:5_Sec:Analysis:FM}. Then we are either done or can assume $\nu_{\underline{\mathcal{I}}}=0$ and it suffices to consider order reduction in the monomial case.


\subsubsection{Monomial ideals and order reduction}
\label{C5_S3_S3}  
Let $\underline{\mathcal{I}} := (X, N, P, \mathcal{I}, E, m)$ be a $\Lambda$-marked ideal such that 
$$
\mathcal{I}=\mathcal{M}(\mathcal{I})= \mathcal{O}_X(-\sum_{j \in J} a_j E^{j}),
$$
where here $J \subset I$ for $I$ the totally ordered set indexing $E=\sum_{i \in I} E^i$. We show how we can do order reduction when $\mathcal{I}$ is monomial using the method of Kollar \cite[3.111 Step 3]{10.2307/j.ctt7rptq}. We outline the first step of the process to highlight that to ensure permissible blow-ups of $\Lambda$-marked monomial ideals we require centres as in \cite[3.111 Step 3]{10.2307/j.ctt7rptq} that intersect with $P$. Define $J_P:=\{j \in J: E^j \cap P \not =\emptyset   \}$ so that
$
\mathcal{O}_P\mathcal{I}= \mathcal{O}_P(-\sum_{j \in J_P} a_j (E^{j}|_P)).
$
\\
\\\textbf{Step 1.} Let $K_P \subset I$ be an indexing set as in \cref{lemma4503} such that $Z_P \cap N=P$ for $Z_P = \bigcap_{j \in K_P} E^j$. Find the smallest $j \in J_P$ such that $a_j \geq m$ is maximal. If no such $j$ can be found then proceed to the next step. Blow up $Z_P \cap E^j$, which is nonempty as $Z_P \subset P$ and $E^j \cap P \not = \emptyset$ and permissible as $Z_P \cap E^j \cap N =E^j \cap P \subset supp \underline{\mathcal{I}}$. Then $\bigcap_{j \in J_P} (E^j)' \cap (E^j)' =(Z_P)' \cap (E^j)' = \emptyset$. But $(Z_P)' \cap N'=P'$ by the proof of \cref{corollary4505} and so
$
\mathcal{O}_{P'}\mathcal{I}'= \mathcal{O}_{P'}(-\sum_{j \in J_P'} a_j ((E^{j})')|_{P'})).
$
Here $J_P'=J_P-j+j'$, where $j'$ corresponds to the exceptional divisor of the blow-up and $a_{j'}=a_j-m$, i.e. we remove $E^j$ as it is now disjoint from $P'$ and replace it with the exceptional divisor. Proceeding in this manner we have the transform $\underline{\mathcal{I}}'$ after a permissible sequence of blow-ups where we can assume $a_j < m$ for all $j \in J_{P'}$. 
\\
\\\textbf{Step $>$ 1.} We proceed exactly as in  \cite[3.111 Step 3]{10.2307/j.ctt7rptq} using lexicographic order. The algorithm terminates after step $n$ (where $n=dim_{\mathbb{Z}}P$) and $supp \underline{\mathcal{I}} = \emptyset$.

\begin{remark}[Algorithmic resolution]
    \label[remark]{remark5332} 
    As $K_P$ is not necessarily unique for each blow-up we make a choice of $K_P$ that gives the same $P', N'$ but may give different $X'$. In \cite{BM1} the proof of order reduction similarly encounters choices for blow-up centres, which can be made algorithmic by introducing an inductive invariant as in \cite[Section 9]{BM1} based on an ordering of codimension one orbit closures. As we have a total ordering on $E$ it seems we could similarly include an inductive invariant to make our order reduction algorithmic. Even without such an invariant in the monomial case we could choose $Z_P = \bigcap_{j \in K_P} E^j$ such that if we order the $j \in K_P$ as $j_1 < \dots < j_r $ then the tuple $(r, j_1, \dots, j_r)$ has smallest lexicographic order.
\end{remark}

\begin{remark}[Generic fiber]
    Let $\underline{\mathcal{I}}=(X, X, X, \mathcal{I}, E, m)$ be a monomial ideal with trivial embeddings. As blow-ups commute with fibers (see \cref{lemma5224}) the blow-up sequence (and hence order reduction) over the generic fiber $\underline{\mathcal{I}}_{0}$ will be precisely Kollar's arguments applied to $(X_{0}, (\mathcal{I}_{0}, m), E_{0})$.
\end{remark}

\section{Resolution from order reduction}
\label{C:6}

We conclude by describing how order reduction for $\Lambda$-marked monomial ideals proved in \cref{C:5} implies the existence of weak resolution for certain embedded $\Lambda$-schemes using principalization as in \cite{10.2307/j.ctt7rptq}.

\begin{theorem}[Principalization]
    \label[theorem]{theorem6101}
    (See \cite[Theorem 3.21]{10.2307/j.ctt7rptq}).
    Let $M$ be a smooth $\Lambda$-scheme over $\mathbb{Z}$ with an affine $\Lambda$-cover and a $\mathbb{Z}$-flat locally toric simple normal crossings divisor $E=\sum_{j \in J} E^j$ with $J$ totally ordered. Let $\mathcal{I}$ be an ideal sheaf that is locally generated by monomials of $E$ (see \cref{C:4_Locally monomial ideals}). Then there exists a sequence of blow-ups by smooth centres over $\mathbb{Z}$
    $$
    M_r \xrightarrow{\pi_{r}} M_{r-1} \xrightarrow{\pi_{r-1}} \dots \xrightarrow{\pi_1} M_0=M
    $$
    such that writing $\Pi_r$ as the composite of the blow-up morphisms, $(\Pi_r)^{-1}\mathcal{I}$ is \textit{locally principal}; generated by a single monomial (see \cref{C5_S3_S3}) for the total transform $E_r$. Moreover, $\Pi_r$ is an isomorphism over $M \backslash cosupp(\mathcal{I})$.
\end{theorem}

\textit{Proof.} Let $\underline{\mathcal{I}}:=(M, M, M, \mathcal{I}, E, 1)$ be the tuple with trivial embedding and marking $1$. By \cref{lemma3323} we can assume $E$ is intersection distinguished for the $\Lambda$-cover of $M$ and so $\underline{\mathcal{I}}$ is a $\Lambda$-marked monomial ideal. We can then apply order reduction to $\underline{\mathcal{I}}$ and use the arguments in the proof of \cite[Theorem 3.21]{10.2307/j.ctt7rptq}. \qedsymbol 

\begin{theorem}[Weak embedded resolution from principalization]
    \label[theorem]{theorem6102} 
    (See \cite[Corollary 3.22]{10.2307/j.ctt7rptq}).
    Let $M$ be a smooth $\Lambda$-scheme with an affine $\Lambda$-cover and a $\mathbb{Z}$-flat locally toric simple normal crossings divisor $E=\sum_{j \in J} E^j$ with $J$ totally ordered. Let $X \hookrightarrow M$ be an irreducible closed $\Lambda$-subscheme of codimension $\geq 2$ such that the ideal $\mathcal{I}$ defining $X$ is locally generated by monomials of $E$. There exists a morphism
    $$
    X' \xrightarrow[]{f} X,
    $$
    where $X'$ is a smooth $\Lambda$-scheme over $\mathbb{Z}$ and $f$ is $\Lambda$-equivariant, projective and birational.
\end{theorem}

Before proving the theorem we prove the following lemma, which we use to show $\Lambda$-equivariance of the resolution.

\begin{lemma}
    Let the following be a commuting diagram of $\Lambda$-schemes
    $$
    \begin{tikzcd}[cramped]
	   {X'} && {M'} \\
	   \\
	   X && M
	   \arrow["h", from=1-1, to=1-3]
	   \arrow["f", from=1-1, to=3-1]
	   \arrow["g", from=1-3, to=3-3]
	   \arrow["i"', hook, from=3-1, to=3-3]
    \end{tikzcd},
    $$
    where $i$ is a monomrphism in the category of schemes and $i, g, h$ are $\Lambda$-equivariant. Then $f$ is $\Lambda$-equivariant.
\end{lemma}

\textit{Proof.} Consider taking the diagram in \textbf{Sp} with morphisms from $W^{\ast}$ defining the $\Lambda$-structure. By $\Lambda$-equivariance of $i,g,h$ for all $n$,
$$
\begin{tikzcd}[cramped]
	{W^{\ast}_n(X')} && {X'} \\
	\\
	{W^{\ast}_n(X)} && X & M
	\arrow[from=1-1, to=1-3]
	\arrow[from=1-1, to=3-1]
	\arrow["f", from=1-3, to=3-3]
	\arrow[from=3-1, to=3-3]
	\arrow["i", hook, from=3-3, to=3-4]
\end{tikzcd}
$$
commutes. But $i$ a monomorphism implies the square itself must commute in \textbf{Sp} as $W^{\ast}_n(X), W^{\ast}_n(X'), X, X'$ are schemes. Taking the colimit gives $\Lambda$-equivariance of $f$. \qedsymbol

\textit{Proof of Theorem 5.0.2.} The arguments of \cite[Corollary 3.22]{10.2307/j.ctt7rptq} can be applied, using principalization above. This induces a morphism from strict transforms $\Pi_j: V_j \hookrightarrow X$, where $V_j$ is a $\mathbb{Z}$-smooth $\Lambda$-scheme by observation \ref{observation3316}. The preceding lemma shows that $\Pi_j$ is $\Lambda$-equivariant. \qedsymbol

\begin{corollary}[Fibers of resolutions]
    Let $(p)$ be a prime ideal of $\mathbb{Z}$. Then the fiber of $f$ over $(p)$,
    $$
    (X')_p \rightarrow X_p,
    $$
    is a projective morphism from $(X')_p$, a smooth $\mathbb{F}_p$-variety.
\end{corollary}

\begin{remark}
    In the previous corollary if $p=0$ then the fiber $f_0$ is very similar to the proof of \cite[Corollary 3.22]{10.2307/j.ctt7rptq} applied to $X_0 \rightarrow M_0$. In fact, if the generic point $\eta_X$ is contained in $X_0$ then $f_0$ is precisely the morphism produced by \cite[Corollary 3.22]{10.2307/j.ctt7rptq} and as blow-ups over $\mathbb{Z}$ commute with generic fibers $f_0$ is a composition of blow-ups of $(X_j)_0$ by $(Z_j \cap X_j)_0$.
\end{remark}

\begin{observation}[Limitations of principalization]
    \label{observation6201} 
    Theorem 5.0.2 has a number of limitations.
   \begin{enumerate}
       \item The resolution relies on the order reduction of \cref{C:5} and so is restricted to a limited class of closed subschemes.
       \item  For $p \not = 0$ the fiber of a weak resolution is not necessarily a blow-up sequence.
       \item  The intermediate strict transforms of $X$ need not have a $\Lambda$-structure despite the composition of blow-up morphisms being $\Lambda$-equivariant.
       \item Weak resolutions are not isomorphisms over the non-singular locus.
   \end{enumerate}

\end{observation}

\begin{example}[Monoid schemes]
    By restricting to monoid schemes we can address observations \ref{observation6201} (2) and (3). Let $M$ be a smooth $\Lambda$-scheme with an affine $\Lambda$-cover coming from the $\mathbb{Z}$-realization of a (separated) monoid scheme $M_{MSch}$ (see \cref{example2204}). Let $E=\sum_{j \in J} E^j$ be a $\mathbb{Z}$-flat locally toric simple normal crossings divisor on the cover of $M$ such that each $E^j$ is the $\mathbb{Z}$-realization of a principal equivariant ideal of $M_{MSch}$. Assume $J$ is totally ordered. Let $X \hookrightarrow M$ be an irreducible closed $\Lambda$-subscheme such that the ideal $\mathcal{I}$ defining $X$ is locally generated by monomials of $E$. It follows that $X$ is itself the $\mathbb{Z}$-realization of a closed equivariant monoid subscheme $X_{MSch}$ of $M_{MSch}.$ By \cref{theorem6102} we have 
    $
    X' \xrightarrow[]{f} X,
    $
    where $X'$ is a smooth $\Lambda$-scheme over $\mathbb{Z}$ and $f$ is $\Lambda$-equivariant, projective and birational. By \cref{lemma3315} the total transforms of $E$ will be given by $\mathbb{Z}$-realizations of monoid schemes and thus so will each blow-up centre $Z_k$. Thus every $Z_k \cap X_k$ will be the $\mathbb{Z}$-realization of an equivariant monoid closed subscheme of $X_k$ and each $X_k$ is a $\Lambda$-scheme with toric $\Lambda$-structure. This addresses observation \ref{observation6201} (2). Moreover, blowing up by $\mathbb{Z}$-realizations of monoid schemes commutes with fibers, which resolves observation \ref{observation6201} (3).
\end{example}

In upcoming work, again restricting to $\mathbb{Z}$-realizations of monoid schemes we will prove a stronger embedded resolution theorem based on work of Bierstone and Milman. The theorem will satisfy observation \ref{observation6201} (4) and observation \ref{observation6201} (2) and (3) will be automatic. We address observation \ref{observation6201} (1) by considering embeddings that are local "binomials". The proof of the theorem reduces to an application of order reduction for $\Lambda$-marked monomial ideals developed in \cref{C:5}.

\medskip
\section*{Acknowledgements}
The author would like to thank Prof. Christian Haesemeyer for his supervision, generosity and invaluable guidance throughout his doctoral studies and the writing of this paper. The author would also like to acknowledge the generosity of the David Lachlan Hay Memorial Fund, who's support through the University of Melbourne Faculty of Science Postgraduate Writing-Up Award made the writing of this paper possible.

\goodbreak
\label{Bibliography}


\bibliography{Bibliography}

\begin{thebibliography}{10}

\bibitem{BM1}
Edward Bierstone and Pierre Milman.
\newblock Desingularization of toric and binomial varieties.
\newblock {\em Journal of Algebraic Geometry}, 15, 12 2004.

\bibitem{bierstone2007functoriality}
Edward Bierstone and Pierre~D. Milman.
\newblock Functoriality in resolution of singularities.
\newblock {\em Publications of the Research Institute of Mathematical Sciences}, 44(2):609--639, 2008.

\bibitem{borger2009lambdarings}
James Borger.
\newblock Lambda-rings and the field with one element, 2009.
\newblock preprint.

\bibitem{borger2015basicII}
James Borger.
\newblock The basic geometry of {W}itt vectors. {II}: Spaces.
\newblock {\em Mathematische Annalen}, 351:877--933, 2010.

\bibitem{borger2015basicI}
James Borger.
\newblock The basic geometry of {W}itt vectors, {I}: The affine case.
\newblock {\em Algebra and Number Theory}, 5(2):231--285, 2011.

\bibitem{Bruns_Herzog_1998}
Winfried Bruns and H.~Jürgen Herzog.
\newblock {\em Cohen-Macaulay Rings}.
\newblock Cambridge Studies in Advanced Mathematics. Cambridge University Press, 2 edition, 1998.

\bibitem{monoid}
Guillermo Cortiñas, Christian Haesemeyer, Mark Walker, and Charles Weibel.
\newblock Toric varieties, monoid schemes and $cdh$ descent.
\newblock {\em Journal für die reine und angewandte Mathematik (Crelles Journal)}, 2015, 06 2011.

\bibitem{cutkoskyresolution}
S.D. Cutkosky.
\newblock {\em Resolution of Singularities}.
\newblock Graduate studies in mathematics. American Mathematical Soc., 2004.

\bibitem{10.2307/j.ctt1b7x7vc}
William Fulton.
\newblock {\em Introduction to Toric Varieties. (AM-131)}.
\newblock Princeton University Press, 1993.

\bibitem{PMIHES_1967__32__5_0}
Alexander Grothendieck.
\newblock {\'E}l\'ements de g\'eom\'etrie alg\'ebrique : {IV.} {{\'E}tude} locale des sch\'emas et des morphismes de sch\'emas, {Quatri\`eme} partie.
\newblock {\em Publications Math\'ematiques de l'IH\'ES}, 32:5--361, 1967.

\bibitem{2a777e2f-398f-34da-9fbc-ecc37bb2c7e1}
Heisuke Hironaka.
\newblock Resolution of singularities of an algebraic variety over a field of characteristic zero: {II}.
\newblock {\em Annals of Mathematics}, 79(2):205--326, 1964.

\bibitem{10.2307/j.ctt7rptq}
János Kollár.
\newblock {\em Lectures on Resolution of Singularities (AM-166)}.
\newblock Princeton University Press, 2007.

\bibitem{DifferentialOperatorsAlgorithmicWeightedResolution}
Jonghyun Lee.
\newblock Differential operators and algorithmic weighted resolution.
\newblock Honours thesis, Brown University, 04 2020.

\bibitem{stacks-project}
{Stacks project authors}.
\newblock The stacks project.
\newblock \url{https://stacks.math.columbia.edu}, 2024.
\newblock accessed 2024-08-18.

\bibitem{FOAG}
Ravi Vakil.
\newblock {MATH 216: Foundations of Algebraic Geometry}, 2017.
\newblock November 18, 2017 draft accessed 2021-02-03.

\end{thebibliography}

\end{document}